\documentclass[12pt,oneside,openany,article]{memoir}
\usepackage{mempatch}

\nouppercaseheads 
\usepackage[amsthm,thmmarks,hyperref]{ntheorem}
\usepackage[noTeX]{mmap}
\usepackage{etex}

\usepackage[english]{babel}
\usepackage[leqno]{mathtools}

\usepackage{mathrsfs}
\usepackage{upgreek}
\usepackage[compress,square,comma,numbers]{natbib}
\usepackage{hypernat} 

\usepackage{sfmath}

\usepackage{microtype}

\usepackage{subfig}
\usepackage{wrapfig}
\usepackage{array}

\usepackage{graphicx,color}
\definecolor{gray}{gray}{0}
\pagecolor{white}

\usepackage{hyperxmp}
\usepackage{xr-hyper}
\usepackage{nameref}
\usepackage[pdftex,bookmarks,pdfnewwindow,plainpages=false,unicode]{hyperref}

\usepackage{bookmark}

\usepackage{enumitem}

\usepackage{amssymb} 

\hypersetup{
colorlinks=true,
linkcolor=black,
citecolor=black,
urlcolor=blue,
pdfauthor={Victor Ivrii},
pdftitle={Short Loops and Pointwise Spectral Asymptotics},
pdfsubject={Sharp Spectral Asymptotics},
pdfkeywords={Microlocal Analysis,  Sharp Spectral Asymptotics, Schort loops},
bookmarksdepth={4}
}

\numberwithin{equation}{chapter}

\theoremstyle{plain}
\newtheorem{theorem}{Theorem}[chapter]

\newtheorem{proposition}[theorem]{Proposition}
\newtheorem{corollary}[theorem]{Corollary}

\theoremstyle{definition}

\newtheorem{definition}[theorem]{Definition}
\newtheorem{Problem}[theorem]{Problem}
\theoremstyle{remark}
\newtheorem{remark}[theorem]{Remark}

\newtheorem{problem}[theorem]{Problem}

\numberwithin{equation}{chapter}

\DeclareMathAlphabet{\mathpzc}{OT1}{pzc}{m}{it}

\newcommand{\cA}{\mathcal{A}}
 \newcommand{\cB}{\mathcal{B}}
 \newcommand{\cC}{\mathcal{C}}
 
 \newcommand{\cE}{\mathcal{E}}
 
 \newcommand{\cG}{\mathcal{G}}

 \newcommand{\cL}{\mathcal{L}}
 \newcommand{\cM}{\mathcal{M}}
 \newcommand{\cN}{\mathcal{N}}

 \newcommand{\cR}{\mathcal{R}}
 \newcommand{\cQ}{\mathcal{Q}}
 \newcommand{\cS}{\mathcal{S}}
 \newcommand{\cT}{\mathcal{T}}

 \newcommand{\cX}{\mathcal{X}}
 \newcommand{\cY}{\mathcal{Y}}
 \newcommand{\cZ}{\mathcal{Z}}

 \newcommand{\sC}{\mathscr{C}}

 \newcommand{\sF}{\mathscr{F}}

 \newcommand{\sL}{\mathscr{L}}

\newcommand{\inter}{{\mathsf{int}}}

\newcommand{\st}{{\mathsf{st}}}

\newcommand{\corr}{{\mathsf{corr}}}

\newcommand{\D}{{\mathsf{D}}}
\newcommand{\ext}{{\mathsf{ext}}}

\newcommand{\MW}{{\mathsf{MW}}}
\newcommand{\N}{{\mathsf{N}}}

\newcommand{\R}{{\mathsf{R}}}

\newcommand{\T}{{\mathsf{T}}}
\newcommand{\trans}{{\mathsf{trans}}}

\newcommand{\W}{{{\mathsf{W}}}}
\newcommand{\w}{{\mathsf{w}}}

\newcommand{\const}{{\mathsf{const}}}

\newcommand{\new}{{{\mathsf{new}}}}

\newcommand{\bN}{{\mathbb{N}}}

\newcommand{\bR}{{\mathbb{R}}}

\newcommand{\bZ}{{\mathbb{Z}}}

\newcommand{\fz}{{\mathfrak{z}}}

\def\1{\boldsymbol {|}}
\newcommand{\brangle}{{\boldsymbol{\rangle}}}

\newcommand{\Def}{\mathrel{\mathop:}=}

\newcommand{\Hess}{\operatorname{Hess}}

\newcommand{\mes}{\operatorname{mes}}

\newcommand{\rank}{\operatorname{rank}}
\newcommand{\supp}{\operatorname{supp}}

\newcommand{\Tr}{\operatorname{Tr}}

\newenvironment{claim}[1][{\textup{(\theequation)}}]{\refstepcounter{equation}\vglue10pt
\begin{trivlist}
\item[{\hskip\labelsep#1}]}{\vglue10pt\end{trivlist}}

\newenvironment{claim*}[1][{}]{\vglue10pt
\begin{trivlist}
\item[{\hskip\labelsep#1}]}{\vglue10pt\end{trivlist}}

\newenvironment{phantomequation}[1][]{\refstepcounter{equation}}{}
\newcounter{note}

\newenvironment{Note}[1]{\refstepcounter{note}{$^{\thenote}$}\reversemarginpar
\marginpar{\hskip-40pt\fbox{\begin{minipage}{79pt} {$^{\thenote}$\  #1}\end{minipage}}}}

\DeclareTextCommand{\textinfty}{PU}{\9042\036}

\DeclareTextCommand{\textge}{PU}{\9042\145}
\DeclareTextCommand{\textle}{PU}{\9042\144}
\DeclareTextCommand{\texthat}{PD1}{\136}

\begin{document}
\title{$2\D$- and $3\D$-magnetic Schr\"odinger Operator with Irregular Coefficients}
\author{Victor Ivrii}

\maketitle
{\abstract%
Now we start analysis of the Schr\"odinger operator with the strong magnetic
field. Keeping in mind the crucial role played by the microlocal canonical form in Chapter~\ref{book_new-sect-13} of \cite{futurebook}, devoted to such operators in the smooth case, the lack of the ultimate smoothness here seems to be the catastrophic because complete microlocal canonical form requires such smoothness. However, combining ideas of Chapter~\ref{book_new-sect-13} and Sections~\ref{book_new-sect-4-5}  and  \ref{book_new-sect-5-3} of \cite{futurebook} with an idea of the incomplete reduction, with a remainder which is just small enough for our needs, we manage under very moderate smoothness conditions to derive sharp remainder estimates: namely, under non-degeneracy assumption $\nabla (V/F)\ne 0$, assuming a bit more than $\sC ^1$ and $\sC ^2$ for $d=3$ and $d=2$ respectively \footnote{Actually, the metrics and the magnetic field should belong to $C^4$ as $d=2$ and $\mu \gg h^{-1}$.}, we recover remainder
estimates $O(h ^{-2}+\mu h ^{-1})$ and $O(\mu ^{-1}h ^{-1}+1)$ respectively
exactly as in the smooth case. Further, for $d=3$ without this non-degeneracy
condition we recover remainder estimate $O(h ^{-2}+ \mu h ^{-1-2/ (l+2)})$ in $\sC ^l$-case, $1\le l\le 2$. Here $h\ll 1$ and $\mu\gg 1$ are Planck and coupling constants respectively.

In comparison with version 1 of 5.5 year ago this version contains more results and some minor corrections.
\endabstract}

\setcounter{chapter}{-1}
\chapter{Introduction}
\label{book_new-sect-18-0}

\chapter{Introduction}
\label{book_new-sect-18-1}

\section{Preliminary}
\label{book_new-sect-18-1-1}

This Chapter is a continuation of Chapter~\ref{book_new-sect-13} and Sections~\ref{book_new-sect-4-5} and~\ref{book_new-sect-5-3}. Here we deal with the Schr\"odinger operator with the strong magnetic field in dimensions $2$ and $3$ under weak smoothness assumptions.

Keeping in mind the crucial role played by the microlocal canonical form in Chapter~\ref{book_new-sect-13}, devoted to such operators in the smooth case, the lack of the ultimate smoothness here seems to be the catastrophic because complete microlocal canonical form requires such smoothness. Instead, in this Chapter we will make a \emph{partial reduction\/}\index{canonical form!partial reduction} with a remainder which is just small enough for our needs but not excessively small, and such reduction requires much lesser smoothness.

This reduction serves two purposes.

First, if the magnetic field is not very strong $\mu \le h^{-\frac{1}{3}}$\,\footnote{\label{foot-18-1} We skip logarithmic factors $|\log h|^\beta$ in this preliminary part.}, we will use it only to prove that the singularities leave the diagonal and propagate away from it and thus the time interval for which we have the standard decomposition of the propagator trace could be extended from $\{|t|\le T_*=\bar{T}\Def \epsilon \mu ^{-1}\}$\,\footnote{\label{foot-18-2} This result is just a rescaling of the standard theory (i.e. the theory with $\mu =1$).} to a larger one thus improving remainder estimate $O(\mu h^{1-d})$ which is due to this rescaling as well.

Namely, for $d=2$ we use non-degeneracy assumption $|\nabla (V/F)|\ne 0$ and
extend time to $\{|t|\le T^*=\epsilon \mu \}$ thus improving estimate to
$O(\mu ^{-1}h^{-1})$. We also consider a weaker non-degeneracy assumption $|\nabla (V/F)|+|\det \Hess (V/F)|\ne 0$ but it will require at the end more smoothness.

For $d=3$ no non-degeneracy assumption at this stage is needed: we consider zone $\{|\xi_3| \ge C\mu^{-1}\}$ where $\xi _3$ is a velocity along magnetic field and we use this movement to prove that the singularities leave the diagonal, increasing time interval to 
$\{|t| \le T^*=T^*(\xi_3)=\epsilon |\xi _3|\}$\,\footnote{\label{foot-18-3} Actually, a bit better but we skip logarithmic factors here.}
and another zone $\{|\xi_3| \le C\mu^{-1}\}$ which is small enough to estimate it contribution to the remainder by $O(\mu ^{\frac{3}{2}}h^{-\frac{3}{2}})= O(h^{-2})$.

On the other hand, if the magnetic field is strong enough, we use this canonical form to derive spectral asymptotics as well. For $d=3$ and $\mu h\le 1$ we do this only in the zone $\{|\xi_3| \le C_0(\mu h)^{\frac{1}{2}}\}$, treating zone
$\{|\xi_3| \ge C_0(\mu h)^{\frac{1}{2}}\}$ in the same way as before. Only at this stage the qualified estimate of the remainder in canonical form is needed.

Combining this idea of this ``poor man canonical form'' with the ideas of Chapter~\ref{book_new-sect-13}, and sections~\ref{book_new-sect-4-5} and~\ref{book_new-sect-5-3} we manage to derive sharp spectral asymptotics under very moderate smoothness conditions: namely, 
\begin{enumerate}[label=(\roman*), leftmargin=*]
\item For $d=2$  assuming 
$\sC ^{2,1}$-regularity  under non-degeneracy assumption $\nabla (V/F)\ne 0$ we recover remainder estimate $O\bigl(\mu ^{-1}h ^{-1}+1\bigr)$  (and the same or almost the same remainder estimate under weaker non-degeneracy assumption) -- exactly as in the smooth case; under weaker smoothness assumption remainder estimate would be larger;

\item For $d=3$, assuming  $\sC ^{1,2}$-regularity  we recover remainder estimates $O\bigl(h ^{-2}+\mu h ^{-1}\bigr)$ (and a the same remainder estimate under weaker non-degeneracy assumption)  -- exactly as in the smooth case; 

\item For $d=3$ without any non-degeneracy assumption we recover remainder estimate  $O\bigl(h ^{-2}+ \mu h ^{-1- 2/(l+2)}\bigr)$ in
$\sC ^{l,\sigma}$-case, with appropriate $\sigma $.
\end{enumerate}

\begin{remark}\label{rem-18-1-1}
(i) Actually in the case of stronger magnetic field we always assume that the metrics and the magnetic field should belong to $\sC^{2,1}$;

\medskip\noindent
(ii) To recover remainder estimate $O(1)$ for $d=2$ and $\mu \gg h^{-1}$ we need to assume that the metrics and the magnetic field should belong to $\sC^4$.
\end{remark}

Recall that here and below we use 

\begin{definition}\label{def-18-1-2}
$\sC ^{l,\sigma }$ denotes a space of functions which 
$\lfloor l \rfloor$-derivatives are continuous with the continuity modulus
$|x-y| ^{l-\lfloor l \rfloor}\bigl|\log |x-y|\bigr| ^{-\sigma }$ if
$l \notin \bN$ or $l \in \bN$ and $\sigma >0$; for $l\in \bN$ and $\sigma \le 0$
one should take $l-1$ instead of $\lfloor l\rfloor$.
\end{definition}

Also $\sigma $ will be positive and large enough in the most of the statements\footnote{\label{foot-18-4} In contrast to non-magnetic case we are not obsessed here to making $\sigma $ as small as possible but we are still pushing hard $l$ down.}.

In this Chapter we consider $d=2$ and $d=3$ only. For $d\ge 4$ the sharp remainder estimate was not derived so far in this book  even in the smooth case due to lack of the canonical form\footnote{\label{foot-18-5} Because of variable multiplicities of the eigenvalues of the matrix magnetic intensity and their ``high-order resonances'' when one of their linear combination with integer coefficients vanishes.}. However, in the next few chapters we will to cover multidimensional case as well, using the idea of the partial canonical form introduced here\Note{Pending. See my papers \cite{ivrii:IRO4} for non-degenerate case like $d=2$, \cite{ivrii:IRO5} for uniformly degenerate case like $d=3$ and \cite{ivrii:IRO8, ivrii:IRO9} for $4$-dimensional degenerating case.}. However, a higher smoothness may be required.

So, let us consider Schr\"odinger operator in the strong magnetic field (\ref{book_new-13-1-1}) i.e. 
\begin{multline}
A=A_0+V(x),\qquad A_0=\sum_{j,k\le d}P_jg ^{jk}(x)P_k, \\
P_j=hD_j- \mu V_j(x),
\quad \mu \ge 1
\label{18-1-1}
\end{multline}
with real symmetric positive metric tensor $(g ^{jk})$, real-valued $V_j$, $V$ and small parameter $h$ and large parameter $\mu $.

We assume that the corresponding \emph{magnetic field intensity\/} 
\begin{equation}
F ^l(x)= \frac{1}{2}\sum _{jk}\upvarepsilon ^{jkl}
\bigl(\partial _{x_j}V_k-\partial _{x_k}V_j\bigr)
\label{18-1-2}
\end{equation}
with the standard absolutely skew-symmetric tensor 
$\upvarepsilon ^{jkl}$\!\footnote {\label{foot-18-6} With non-zero components equal $\pm g^{\frac{1}{2}}$ where here and
in what follows $g=\det (g^{jk})$ and $(g_{jk})=(g^{jk})^{-1}$.}
does not vanish in domain in question. Note that for $d=2$ automatically
$F ^1=F ^2=0$ and thus magnetic field is reduced to $F ^3$.

\begin{remark}\label{rem-18-1-3}
For $d=3$ we can straighten vector field $\mathbf{F}=(F^1,F^2,F^3)$ and direct
it along $x_3$ by an appropriate choice of variables $x_1,x_2$ and make $V_3=0$ by an appropriate \emph{gauge transformation\/}
$A \mapsto e ^{-i\mu h ^{-1}\Phi (x)}Ae ^{i\mu h ^{-1}\Phi }$.

We assume that this is has been done and impose smoothness conditions to the reduced operator.
\end{remark}

Thus we assume that for this reduced operator
\begin{equation}
V_3=0,\qquad \partial _{x_3}V_k=0 \quad k=1,2
\label{18-1-3}
\end{equation}
where the second condition is due to the first one and assumption $F^1=F^2=0$;
then one can reduce it further to $V_2=0$ and then
\begin{equation}
\mathbf{F}=(0,0,F ^3),\quad F ^3=
g ^{-{\frac{1}{2}}}(\partial _{x_1}V_2-\partial _{x_2}V_1),\qquad
g=\det (g ^{jk}) ^{-1}
\label{18-1-4}
\end{equation}
and the scalar intensity of magnetic field is
\begin{multline}
F=g_{33}^{\;{\frac{1}{2}}} |F ^3| =
(g ^{11}g ^{22}-g ^{12}g ^{21})^{\frac{1}{2}}
|(\partial _{x_1}V_2-\partial _{x_2}V_1)|,\\
(g_{jk})=(g ^{jk}) ^{-1}.
\label{18-1-5}
\end{multline}
In particular for $d=2$ one can assume that $g ^{33}=1$, $g ^{3j}=0$ for
$j=1,2$; without taking absolute value we get \emph{pseudo-scalar intensity\/}
then.

Further, we assume that
\begin{equation}
F\ge c^{-1}.
\label{18-1-6}
\end{equation}

Multiplying operator by $F^{-\frac{1}{2}}$ from both sides one can further reduce to the case
\begin{equation}
F=1;
\tag*{$\textup{(\ref*{18-1-6})}^*$}
\label{18-1-6-*}
\end{equation}
then operator will be 
\begin{gather}
\sum_{j,k}P_jF^{-1}g^{jk}P_k+F^{-1}(V -h^2V')
\label{18-1-7}\\
\shortintertext{with}
V'=F^{\frac{1}{2}}\sum_{j,k} 
\partial_j\bigl(g^{jk}\partial_k F^{-\frac{1}{2}}\bigr)
\label{18-1-8}
\end{gather}
and this perturbation will not affect our results unless $d=2$, $\mu \ge h^{-1}$ (in which case $V$ should be corrected).

\section{Smooth case: survey}
\label{book_new-sect-18-1-2}

From Chapter~\ref{book_new-sect-13} we know that for $d=2$ in the smooth case one can reduce operator further to the microlocal canonical form
\begin{gather}
\mu ^2 \sum_{m+n+p\ge 1} \omega_{m,n,p} (x_2,\mu ^{-1}hD_2)
\bigl(x_1 ^2+\mu^{-2} h ^2D_1 ^2\bigr)^m \mu ^{-2n-p} h^p\label{18-1-9}\\
\intertext{with the principal part}
F(x_2,\mu ^{-1}hD_2) \bigl(\mu ^2x_1 ^2+h ^2D_1 ^2\bigr)+V(x_2,\mu ^{-1}hD_2) \tag*{$\textup{(\ref*{18-1-9})}_0$}\label{18-1-9-0}
\end{gather}
and for $d=3$ in the smooth case one can reduce operator to the microlocal
canonical form
\begin{multline}
\mu ^2 \sum_{m+n+p+q\ge 1} \omega_{m,q,k,p}(x_2,x_3, \mu ^{-1}hD_2) \times\\
\bigl(x_1 ^2+\mu ^{-2} h ^2D_1 ^2\bigr)^m (hD_3)^q \mu ^{-2m-2q-2k-p}
h^p\label{18-1-10}
\end{multline}
with the principal part
\begin{equation}
F(x_2,\mu ^{-1}hD_2) \bigl(\mu ^2x_1 ^2+h ^2D_1 ^2\bigr)+h ^2D_3 ^2+V(x_2,x_3,\mu ^{-1}hD_2)
\tag*{$\textup{(\ref*{18-1-10})}_0$}\label{18-1-10-0}.
\end{equation}

Then decomposing into Hermitian functions of 
$\mu ^{\frac{1}{2}}h^{-{\frac{1}{2}}}x_1$ one can replace harmonic oscillator $(h ^2D_1 ^2+\mu ^2x_1 ^2)$ by its eigenvalues
$(2j+1)\mu h$ with $j\in \bZ^+$. So, basically for $d=2$ our operator is a kind of the mixture of the usual second order operator and a family of 
$\mu ^{-1}h$-pseudo-differential operators while for $d=3$ it is a kind of the mixture of the usual second order operator and a family of $h$-differential $1$-dimensional Schr\"odinger operators with respect to $x_3$ which are also $\mu^{-1}h$-pseudo-differential operators with respect to $x_2$. This explains the major differences between $d=2$ and $d=3$:

\medskip
For $d=2$ there are sharp Landau levels $(2j+1)\mu h$ (the eigenvalues of the infinite multiplicity as $g ^{ik}$ and $V$ are constant and domain is $\bR^2$); without some non-degeneracy assumption one cannot expect the remainder estimate better than $O(\mu h ^{-1})$ and with an appropriate non-degeneracy assumption one can improve it up to $O(\mu ^{-1}h ^{-1}+1)$ because only $\asymp (\mu h)^{-1}$ of Landau levels are not classically forbidden.

\medskip
Meanwhile for $d=3$ Landau levels are rather bottoms of the branches of the continuous
spectrum and the non-degeneracy assumption is not that crucial: the best possible remainder estimate is $O\bigl(h ^{-2}+\mu h ^{-1}\bigr)$ while the worst possible remainder estimate depends on the smoothness of the reduced operator (presuming nothing wrong happened during reduction) and is 
$O\bigl(h ^{-2}+\mu h ^{-1- 2/(l+2)}\bigr)$ in the heuristic accordance with section~\ref{book_new-sect-5-3}.

\medskip
Finally, for $d\ge 4$ Landau levels are $\sum_{1\le k \le r}(2j_k+1) f_k \mu h $ where $f_k > 0$, $\pm if _k$ are non-zero eigenvalues of the matrix magnetic
intensity $F$, $2r=\rank F$ and $j \in \bZ^{+\,r}$, and one can think operator
to be similar to $2$-dimensional as $\rank F = d$ and to $3$-dimensional as
$\rank F <d$ but there is no complete canonical form even if $\rank F$ and
multiplicities of $f_k$ are constant because $\sum _k j_k f_k$ can vanish for
some $j\in \bZ^r$ at some points (resonances). This leads to the difficulties in our analysis which we tackle in few next chapters.\Note{Pending}

\section{Non-smooth case: heuristics as $d=2$}
\label{book_new-sect-18-1-3}

In the non-smooth case we do not have the complete canonical form but actually we do not need it. Further, if magnetic field is not very strong we do not need even logarithmic uncertainty principle\footnote{\label{foot-18-7} Recall that in the non-smooth case we approximate operator in question by operator with $\varepsilon$-smooth coefficients.} (\ref{18-1-12}) below to justify quantization of $V(x_2,\xi_2)$ because we do not quantize it in this case. 

Consider $\xi_2$-shift for time $T$ in the propagation. Note that it does not exceed $CT$ and under non-degeneracy assumption $|\partial_{x_2} V|\ge \epsilon$ the shift is exactly of magnitude $\asymp T$ as $T\le \epsilon \mu$. This shift is ``observable'' on the ``quantum level'' and the propagator trace vanishes provided logarithmic uncertainty principle $T \times \varepsilon \ge C h |\log h|$ holds; here  $x_2$- scale is $\varepsilon$. 

In order to be able to take $T$ as small as $\bar{T}=\epsilon \mu ^{-1}$ we need to pick up
\begin{equation}
\varepsilon = C\mu h|\log h|
\label{18-1-11}
\end{equation}
which in addition to the semiclassical remainder estimate $O(\mu ^{-1}h^{-1})$ brings an approximation error $O(h^{-2}\varepsilon ^l |\log h|^{-\sigma})$ for coefficients belonging to $\sC^{l,\sigma}$. This error will be small for 
$\mu \le h^{1-\delta}$ and $\delta = 2/(l+2)$ which would take care of $\mu \le h^{1-\delta}$ in the smooth case but we do not have such luxury now.

As $\mu$ grows larger we need to treat our operator as a family of
$\mu^{-1}h$-pseudo-differential operators and then we need the logarithmic uncertainty principle $\varepsilon \times \varepsilon \ge C\mu^{-1}h|\log \mu|$
\begin{equation}
\varepsilon = C(\mu ^{-1} h |\log \mu|)^{\frac{1}{2}}
\label{18-1-12}
\end{equation}
which in addition to the semiclassical remainder estimate $O(\mu ^{-1}h^{-1})$ brings an approximation error 
$O\bigl(h^{-2}\varepsilon ^l |\log h|^{-\sigma}\bigr)$ as $\mu \lesssim h^{-1}$.

On the other hand, for $\mu \gtrsim h^{-1}$ and we consider Schr\"odinger-Pauli operator, the semiclassical remainder estimate is $O(1)$ and an approximation error is $O\bigl(\mu h ^{-1}\varepsilon ^l |\log h|^{-\sigma} \bigr)$ from approximation of $V$ but also 
$O\bigl(\mu ^2\varepsilon ^{\bar{l}} |\log h|^{-\bar{\sigma}} \bigr)$ from approximation of $g^{jk},V_j$ where we assume that $g^{jk},F\in \sC^{{\bar{l}},\bar{\sigma}}$. 

An actual reduction and calculations are the very significant part of our analysis.

For the best results we need to pick $\varepsilon$ to be the minimum of (\ref{18-1-11}) and (\ref{18-1-12}) i.e.
\begin{equation}
\varepsilon \Def 
C\left\{\begin{aligned}
&\mu h|\log h| \qquad 
&&\text{as\ \ } \mu \le (h|\log h|)^{- \frac{1}{3}}\\
&(\mu ^{-1} h |\log \mu|)^{\frac{1}{2}} \qquad 
&&\text{as\ \ } \mu \ge (h|\log h|)^{- \frac{1}{3}}
\end{aligned}\right.
\label{18-1-13}
\end{equation}
and the threshold is $\mu^*_1\Def (h|\log h|)^{- \frac{1}{3}}$. Below it magnetic field is considered to be \emph{weak\/} and above it \emph{strong\/} (and we treat these two cases differently). The second threshold is 
$\mu^*_2 \Def h^{-1}$; for $\mu$ above it we modify our conditions and only finite number of Landau levels matter; we call magnetic field \emph{superstrong\/} then.

\section{Non-smooth case: heuristics as $d=3$}
\label{book_new-sect-18-1-4}

In contrast to $2$-dimensional case, there is a movement along $x_3$ with the speed $|\xi_3|$ and for time $T$ shift with respect to $x_3$ is
$\asymp T\cdot |\xi_3|$ (provided $|\xi_3|\ge C \mu ^{-1}$ and
$T\le |\xi_3|$). This shift is ``observable'' on the ``quantum level'' and the
propagator trace vanishes provided the logarithmic uncertainty principle
$T |\xi_3| \times |\xi_3| \ge Ch|\log h|$ holds; also we need to assume that
$\varepsilon \times |\xi_3|\ge Ch|\log h|$. Here scale with respect to $\xi_3$
is $\epsilon |\xi_3|$. 

Plugging $T=\bar{T}=\epsilon \mu ^{-1}$ we arrive to the inequality
\begin{equation}
|\xi_3|\ge \bar{\rho}\Def 
C\max \bigl(\mu ^{-1}, (\mu h |\log h|)^{\frac{1}{2}}\bigr)
\label{18-1-14}
\end{equation}
and also we need to take $\varepsilon$ depending on $\xi_3$: 
$\varepsilon = Ch|\log h| /|\xi_3|$. Then
the contribution of the \emph{outer zone\/} $\{|\xi_3| \ge \bar{\rho}\}$ to the approximation error will be $O(h^{-2})$ under $\sC^{1,2}$ regularity of the coefficients.

Further, in this case one can replace $T^* (\xi_3) = \epsilon |\xi_3|$ by a bit larger $T^*(\xi_3)= \epsilon |\xi_3| \bigl|\log |\xi_3|\bigr|^2$ and then the contribution of this zone to the remainder does not exceed 
$Ch^{-2}\int T^{*\,-1}(\xi_3)\, d\xi_3$ which
does not exceed $Ch^{-2}$ for this increased $T^*(\xi_3)$.

On the other hand, the contribution of the \emph{inner zone\/}
$\{|\xi_3|\le \bar{\rho}\}$ to the remainder does not exceed
$C\mu \bar{\rho} h^{-2}$ which in turn does not exceed
$Ch^{-2}$ as $\mu \le \mu^*_1\Def (h |\log h|)^{-\frac{1}{3}}$. In this case we take $\varepsilon = Ch\log h \bar{\rho}^{-1}$ in the inner zone and then its contribution to the approximation error does not exceed $Ch^{-2}$ as well. In fact we will increase this threshold to  $\mu^*_1\Def h^{-\frac{1}{3}}$.

\medskip
For $\mu \ge h^{-\frac{1}{3}}$ we will consider the outer zone in the same way as above but in the inner zone $\{|\xi_3|\le \bar{\rho}\}$ we will treat our operator as a family of $1$-dimensional Schrödinger operators. Then one needs to take $\varepsilon = Ch|\log h|/|\xi_3|$ in this zone as well. Now theory breaks in two cases.

\medskip\noindent 
(a) Without non-degeneracy assumption. Here we apply in rather
straightforward manner arguments of subsection~\ref{book_new-sect-5-3-1}.

\medskip\noindent
(b) Under non-degeneracy assumption. It is still not very useful to use shift with respect to $x_3,x_2,\xi_2$ because it lands us with 
$\varepsilon = C\mu h |\log h|$ which leads to the total approximation error
$O(h^{-2})$ only under excessive smoothness assumptions. The role of this
condition is more subtle: after $|\xi_3|^2$-partition with respect to
$x_2,x_3,\xi_2$ this condition assures that the ellipticity condition fails
only for elements with the total relative measure
$\asymp \min \bigl( |\xi_3|^2/\mu h ,1\bigr)$ and this leads us in the end to
both remainder and approximation error not exceeding $Ch^{-2}$. 

The logarithmic uncertainty principle requires condition
$|\xi_3|\ge Ch^{\frac{1}{3}}|\log h|^{\frac{1}{3}}$ which can be weaken to
$|\xi_3|\ge Ch^{\frac{1}{3}}$. Zone $\{|\xi_3|\le Ch^{\frac{1}{3}} \}$ is small enough and its total contribution does not exceed $Ch^{-2}$ as well.

Finally, for $\mu \ge h^{-1}|\log h|^{-1}$ the outer zone disappears.

\section{Plan of the Chapter}
\label{book_new-sect-18-1-5}

Let us discuss the plan of this chapter. Sections~\ref{book_new-sect-18-2}--\ref{book_new-sect-18-5} are devoted to $2$-dimensional case and Sections~\ref{book_new-sect-18-5}--\ref{book_new-sect-18-8} are devoted to $3$-dimensional case.

In section~\ref{book_new-sect-18-2} we consider case $d=2$ when magnetic field is relatively weak. Using ``precanonical form'' we prove our first main theorem \ref{thm-18-2-9}.

In section~\ref{book_new-sect-18-3} we consider case $d=2$ when magnetic field is relatively strong. We derive ``poor man canonical form''. The true (complete) canonical form does not exist because we lack the infinite smoothness. However
precanonical form exists with a rather small error but we need to plug corrected potential $W$ instead of $V$ and their difference is the source of correction term. This corrected potential $W$ is an original potential $V$ averaged along magnetrons and as $(l,\sigma)\succeq (2,1)$ their difference  is sufficiently small and therefore the correction term is sufficiently small as well. 

In section~\ref{book_new-sect-18-4} we prove major Tauberian estimates. Then in section~
\ref{18-4-5} we make calculations, prove statements our main
theorem \ref{thm-18-5-4}  and  finish $2$-dimensional case.

In section~\ref{book_new-sect-18-6} we consider case $d=3$ when magnetic field is relatively weak. We analyze the outer zone defined by (\ref{18-1-14})
(if $(x_3,\xi _3)$ are ``free''(``non-magnetic'') variables) and prove that the
contribution of this zone to the remainder is $O(h^{-2})$; then we prove instantly theorem~\ref{thm-18-6-4} and we improve it to theorem~\ref{thm-18-6-12} to cover the case $\mu\le h^{-\frac{1}{3}}$ rather than $\mu\le (h|\log h|)^{-\frac{1}{3}}$.

In section~\ref{book_new-sect-18-7} we consider case $d=3$ when magnetic field is relatively strong. We derive ``poor man canonical form''.

In section~\ref{book_new-sect-18-8} we analyze the inner zone defined by inequality opposite to (\ref{18-1-14}) as $h^{-\frac{1}{3}}\le \mu \le Ch^{-1}$ and
 and prove major Tauberian estimates. 
 
Then in section~\ref{18-4-9} we make calculations, prove statements our main
theorem \ref{thm-18-9-4}  and \ref{thm-18-9-6} and finish $2$-dimensional case.

\medskip
In what follows we consider operator (\ref{book_new-13-1-1}) assumingly that conditions (\ref{book_new-13-1-1}), (\ref{book_new-13-1-4}), (\ref{book_new-13-1-5}) and (\ref{book_new-13-3-45}).

\section{Standard results rescaled}
\label{book_new-sect-18-1-6}

By \emph{standard results\/} we mean results for $\mu =1$. While rescaling final spectral asymptotics by $x\mapsto \mu x$, $h \mapsto \mu h$ we would get the standard remainder estimate $O(h^{1-d})$ converted into
$O\bigl((\mu h)^{1-d}\cdot \mu ^d\bigr)=O(\mu h^{1-d})$, something which
we want to improve, rescaling of the intermediate results will be much more fruitful.

Let us consider some $(\gamma ,\rho )$-admissible partition in $(x,\xi )$ and $\varepsilon $-mollification with respect to $x$ with $\varepsilon $ depending on $(x,\xi )$ such that\footnote{\label{foot-18-8} In fact both $\gamma $ and especially $\rho $ could be vectors and then assumption (\ref{18-1-15}) should hold for their components.}
\begin{equation}
C_s\rho ^{-1}h|\log h| \le \varepsilon \le \gamma .
\label{18-1-15}
\end{equation}
\emph{At first we assume that $\varepsilon$ is a numerical small parameter\/}.

Let $\tilde{A}$ be such approximate operator and $\tilde{e}(x,y,\tau )$ be its spectral projector. Let $U(x,y,t)$ be the Schwartz kernel of the corresponding propagator $e ^{-ih ^{-1}t\tilde{A}}$. We do not write $\tilde{U}$ as we will never consider propagator for the non-mollified (i.e. original) operator.

\begin{proposition}\label{prop-18-1-4}
Let $\varepsilon$ be a small numerical parameter such that
\begin{gather}
\varepsilon \ge C_s h|\log h|\label{18-1-16},\\
\shortintertext{and}
\mu \le C_s ^{-1}h ^{-1}|\log h| ^{-1}.
\label{18-1-17}
\end{gather}
(i) Then
\begin{equation}
|\int \bar{\phi}_\zeta (\tau )\Bigl(d_\tau \tilde{e}(x,x,\tau )-
F_{t\to h ^{-1}\tau }\bigl(\bar{\chi}_T(t)U(x,x,t)\bigr)\,d\tau
\Bigr)|\le Ch ^s
\label{18-1-18}
\end{equation}
provided $\bar{\phi}(\tau )=1$ for $\tau <-1$,
$\bar{\phi} (\tau )=0$ for $\tau >- \frac{1}{2} $,
$\bar{\chi}(t)=1$ for $|t|\le \frac{1}{2} $,
$\bar{\chi}(t)=0$ for $|t|\ge 1$\,\footnote{\label{foot-18-9} We pick up all auxiliary functions $\phi $, $\chi $ etc satisfying assumptions
of Section~\ref{book_new-sect-2-3}:
\begin{equation}
|D ^\alpha f|\le (CN) ^{(|\alpha |-m)_+}\qquad
\forall \alpha :|\alpha |\le N =\epsilon _s |\log h|
\label{18-1-19}
\end{equation}
with large enough exponent $m$.},
$T=\bar{T}\Def \epsilon_0 \mu ^{-1}$ with small enough constant $\epsilon_0 >0$,
$\zeta \ge C_s\mu h|\log h|$\,\footnote{\label{foot-18-10} So $\zeta$ and $T$ also satisfy logarithmic uncertainty principle $\zeta \cdot T\ge C_sh|\log h|$.}
and $\bar{\chi}_T(t)=\bar{\chi}(t/T)$,
$\bar{\phi}_\zeta (\tau )=\bar{\phi}(\tau /\zeta )$;

\medskip\noindent
(ii) Moreover as $|\tau|\le \epsilon$
\begin{gather}
|F_{t\to h ^{-1}\tau }\bigl(\bar{\chi}_T(t)-
\bar{\chi}_{T_*}(t) \bigr) U(x,x,t) | \le Ch^s\label{18-1-20}\\
\intertext {with $T=\bar{T}$ and $T_*\Def C_s h|\log h|$,}
|F_{t\to h ^{-1}\tau }\bar{\chi}_{T_*}(t)
\bigl(U(x,x,t) - U^0(x,x,t)\bigr) | \le 
Ch^{3-d}+Ch^{1-d} \vartheta(h) \label{18-1-21}\\
\shortintertext{with}
U^0(x,x,t)\Def
h^{-d}\int e^{ih^{-1}t\tau} \,d_\tau \tilde{\cN}^\MW_x (\tau)\label{18-1-22}
\end{gather}
and
\begin{multline}
|F_{t\to h ^{-1}\tau }\bar{\chi}_{T_*}(t)\bigl(U(x,x,t)- U_\st^0(x,x,t)\bigr) | \le \\
C\mu ^2 h^{3-d}+ Ch^{1-d}\vartheta (h) \label{18-1-23}
\end{multline}
with
\begin{equation}
U_\st ^0(x,x,t)\Def
h^{-d}\int e^{ih^{-1}t\tau} \,d_\tau \tilde{\cN}^\W_x (\tau)
\label{18-1-24}
\end{equation}
where both $\tilde{\cN}^\W_x(\tau)$ and $h^{-d}\tilde{\cN}^\MW_x (\tau)$
are constructed for the mollified operator $\tilde{A}$ and
$\vartheta (t) = t ^l |\log t |^{-\sigma}$ as $g^{jk}, F_{jk}, V\in \sC^{l,\sigma}$.
\end{proposition}

\begin{proof} 
(i) Proof of estimate (\ref{18-1-18}) is standard because in the scale
$x\mapsto \mu x$, $t\mapsto \mu t$, $\tau \mapsto \mu ^{-1}\tau $,
$\mu \mapsto 1$, $h \mapsto \mu h$ we have the standard case (with $\mu =1$).
Note that
\begin{multline}
\int \bar{\phi}_\zeta (\tau ) F_{t\to h ^{-1}\tau }
\Bigl(\bar{\chi}_T(t)U(x,x,t)\Bigr)\,d\tau =\\
\bar{\phi}_\zeta (hD_t)\Bigl(\bar{\chi}_T(t) U(x,x,t)\Bigr)\Bigr|_{t=0}.
\label{18-1-25}
\end{multline}
Inequality (\ref{18-1-20}) also follows from the same rescaling arguments applied to the standard results of section~\ref{book_new-sect-4-5} as we take in account the $\xi$-microhyperbolicity which is equivalent to (\ref{18-1-26}): 
\begin{equation}
V\le -\epsilon_0.
\label{18-1-26}
\end{equation}

\medskip\noindent
(ii) To prove (\ref{18-1-21})--(\ref{18-1-22}) and (\ref{18-1-23})--(\ref{18-1-24}) we also apply rescaling arguments and recall construction of $U(x,y,t)$ by the method of successive approximations as $|t|\le C_sh |\log h|$.

However, as an unperturbed operator one can also take
\begin{gather}
\bar{A}(x,y,hD_x)\Def \sum_{j,k} g ^{jk}(y)\bar{P}_j\bar{P}_k + V(y)
\label{18-1-27}\\
\shortintertext{with}
\bar{P}_j=hD_j-\mu V_j (y;x),\qquad 
V_j(y;x)= V_j(y)+\langle \nabla V_j(y), x-y\rangle
\label{18-1-28}
\end{gather}
rather than more standard $V_j(y)$ (where the former includes magnetic field and the latter does not) and here we must plug mollified $\tilde{g} ^{jk}$, 
$\tilde{V}_j$, and ${\tilde V}$ instead of the original (non-mollified) coefficients.

One can see easily that for such operator the restriction of the spectral kernel
to the diagonal $\bar{e}(y,y,\tau )$ is exactly 
$h^{1-d}\tilde{\cN}^\MW_y(\tau )$ and that all other terms of the successive approximation contribute no more than $O(h^{3-d})$.

Really, the contribution of the $k$-th term ($k=1,2,3,\ldots $) does not exceed $Ch^{-d-(k-1)}T^{2k-1}$ with $T=h|\log h|$ as long as $g^{jk}, V\in \sC^{1,1}$ and $V_j\in \sC^{2,1}$; then each next term acquires factor $T^2h^{-1}$. 

Thus contribution of any term with $k\ge 4$ is clearly $O(h^{3-d})$ while contribution of the third term is $O(h^{3-d}|\log h|^5)$. Further, analysis of subsection~\ref{book_new-sect-4-5-4} (based on the rescaling technique) shows that one can take $T=h$ in this estimate and then the contribution of the second term does not exceed $Ch^{3-d}$ as well. 

Consider the second term 
\begin{equation}
\frac {1} {2\pi i} \Bigl(\bar{G} ^+ (\tilde{A}-\bar {A})\bar{G} ^+ -
\bar{G} ^- (\tilde{A}-\bar{A}) \bar{G} ^-\Bigr)\updelta (x-y)\updelta(t)
\label{18-1-29}
\end{equation}
where $\bar{G} ^{\pm}$ are respectively forward and backward parametrices of the evolution operator $hD_t-\bar{A}$.

Consider $\tilde{A}-\bar{A}=R_1+R_2$ where $R_1$, $R_2$ are second order operators with the coefficients
$O\bigl(\mu |x-y|\vartheta (|x-y|)\bigr)$ and $O\bigl(\vartheta (|x-y|)\bigr)$ provided $\vartheta (z)\ge Cz^2$. The same rescaling procedure implies that if we replace $\tilde{A}-\bar{A}$ in (\ref{18-1-29}) by $R_1,R_2$, the final error will be $O\bigl(h^{1-d} \mu h \vartheta (h)\bigr)$ and 
$O\bigl(h^{1-d} \vartheta (h)\bigr)$ respectively. On the other hand, for $\vartheta(z)=O(z^2)$ we get $O(h^{3-d})$.

It implies estimates (\ref{18-1-21})--(\ref{18-1-22}).

\medskip\noindent
(iii) To prove estimates (\ref{18-1-23})--(\ref{18-1-24}) we need to use the same construction but with the standard unperturbed operator $\bar{A}(y,hD_x)$ instead; it is given by (\ref{18-1-27}) with $\bar{P}_j=hD_j-V_j(y)$. Then one can estimate the contribution of the $k$-th term in the successive approximation by $Ch^{-d+k}\mu ^{k-1}$ and thus all the terms but the first two contribute less than $C\mu ^2 h^{2-d}$. The contribution of the first term is exactly $h^{-d}\tilde{\cN}^\W_y(\tau)$ and the only term we need to consider is the second one. 

The same argument as above work but now 
\begin{equation*}
\tilde{A}-\bar{A}=R_0+ R_1+R_2,\qquad
R_0= \mu \langle x-y, \nabla_y \tilde{A}(y,hD_x) \rangle.
\end{equation*}
Then the contributions of the term given by (\ref{18-1-29}) with $\tilde{A}-\bar{A}$ replaced by $R_0$ will be identically $0$ we get a term odd with respect to $(x-y,D_x)$ and therefore its restriction to $\{x=y\}$ vanishes.

It implies estimates (\ref{18-1-23})--(\ref{18-1-24}).
\end{proof}

\begin{corollary}\label{cor-18-1-5}
In frames of proposition~\ref{prop-18-1-4}(ii) 
\begin{multline}
h^{-1} |\int ^\tau F_{t\to h ^{-1}\tau }\bar{\chi}_{T_*}(t)
\bigl(U(x,x,t) - U^0(x,x,t)\bigr) \,d\tau | \le \\
Ch^{2-d}+Ch^{-d} \vartheta(h) \label{18-1-30}
\end{multline}
and
\begin{multline}
h^{-1}\int^\tau |F_{t\to h ^{-1}\tau }\bar{\chi}_{T_*}(t)\bigl(U(x,x,t)-
U_\st^0(x,x,t)\bigr)\,d\tau | \le\\
C\mu ^2 h^{2-d}+ Ch^{-d}\vartheta (h). \label{18-1-31}
\end{multline}
\end{corollary}

\chapter{$d=2$: Weak magnetic field}
\label{book_new-sect-18-2}

We start from $2$-dimensional case which is more transparent and clean.

\section{Heuristics}
\label{book_new-sect-18-2-1}

\subsection{Pilot-model}
\label{book_new-sect-18-2-1-1}

The pilot-model of such operator is $h^2D_1^2(hD_2-\mu x_1)^2+V(x)$ and metaplectic map
$hD_1\mapsto hD_1$, $hD_2 -\mu x_1\mapsto -\mu x_1$,
$x_1 \mapsto x_1+ \mu ^{-1}hD_2$, $x_2 \mapsto x_2 +\mu ^{-1}hD_1$ transforms
this operator into
\begin{equation}
h^2D_1^2 +\mu ^2 x_1^2 + V(x_1+ \mu ^{-1}hD_2, x_2 +\mu ^{-1}hD_1).
\label{18-2-1}
\end{equation}
In some sense this operator is modelled by
\begin{equation}
h^2D_1^2 +\mu ^2 x_1^2 + V( \mu ^{-1}hD_2, x_2 ),
\label{18-2-2}
\end{equation}
Then along coordinates $x_2$ and $\xi_2\sim \mu ^{-1}hD_2$ singularities propagate with the speed $O(\mu ^{-1})$ and under assumption 
$|\partial_{\xi_2} V |\asymp \nu $ for time $T$ the shift with respect to $x_2$ is $\asymp \mu ^{-1}T\nu $. To satisfy logarithmic uncertainty principle we should assume that 
\begin{equation}
\mu ^{-1}T \nu \cdot \varepsilon \ge C\mu ^{-1}h |\log (\mu ^{-1}h )|
\label{18-2-3}
\end{equation}
because $V( \mu ^{-1}hD_2, x_2 )$ is $\mu ^{-1}h$-pseudo-differential operator with respect to $x_2$ and $\varepsilon $ is a characteristic size in $\xi_2$-variable. Therefore,
$\Gamma U \Def \int U(x,x,t)dx $ is negligible for $T_*\le |t| \le T^*$ with \begin{equation}
T_* \Def C\varepsilon ^{-1}\nu^{-1} h |\log (\mu ^{-1}h)|
\label{18-2-4}
\end{equation}
and $T^*=\epsilon \mu $ (as we will show).

On the other hand, from the standard results rescaled we know $\Gamma U $ for
$|t|\le \bar{T} \Def \epsilon \mu ^{-1}$. In order these intervals to overlap we
need  
\begin{equation*}
T_* =C\varepsilon^{-1}\nu^{-1}h|\log (\mu ^{-1}h)|\le \bar{T} =
\epsilon \mu^{-1}
\end{equation*}
i.e. $\varepsilon \ge C\mu h |\log h| \nu^{-1}$ and therefore
$\mu \le h^{-1}|\log h|^{-1}$. So, let us pick up the minimal possible value
\begin{equation}
\varepsilon = C\mu h \nu^{-1}|\log h|
\label{18-2-5}
\end{equation}
and then the approximation error will be 
$Ch^{-2}\vartheta (\mu h \nu^{-1}|\log h|)$.

This construction works for pretty large $\mu $; however, these arguments are
optimal only for
\begin{equation}
\mu \le \bar{\mu}_1 \Def (h  |\log h|)^{-\frac{1}{3}}.
\label{18-2-6}
\end{equation}

\begin{remark}\label{rem-18-2-1}
In the next section we will pick up 
$\varepsilon = C\bigl(\mu ^{-1}h|\log (\mu ^{-1} h) |\bigr)^{\frac{1}{2}}$ which (for $\nu\asymp 1$) is less than  $C\mu h |\log h|$ as $\mu \ge \bar{\mu}_1$.

Note that we need inequality
$\varepsilon \ge C\bigl(\mu ^{-1}h|\log (\mu ^{-1} h) |\bigr)^{\frac{1}{2}}$ even to
consider $V(\mu^{-1}h D_2, x_2)$ as a legitimate $\mu^{-1}h$-pseudo-differential operator. So, in the frames of this section $V(\mu^{-1}h D_2, x_2)$ is not even defined properly. The short answer to this obstacle is that we should work directly with operator $V(x_1+ \mu ^{-1}hD_2, x_2 +\mu ^{-1}hD_1)$ or just not to make any reduction but still to be able to reproduce the above arguments.
\end{remark}

\subsection{Propagation of singularities}
\label{book_new-sect-18-2-1-2}

The goal of the this subsubsection is to prove that at the energy level $0$ and below singularities propagate with respect to $x_1,x_2$ with the finite speed, and that with respect to ``variables'' 
$Q_j\Def x_j + \sum _k \beta _{jk} \mu ^{-1} P_k$ (with appropriate coefficients $\beta _{jk}$) singularities propagate with the speed not exceeding 
$C_0\mu ^{-1}$; further, the same is true for
$Q'_j \Def Q_j + \sum_{ki} \mu ^{-2} \beta '_{jki}P_kP_l$ with appropriate coefficients $\beta '_{jki}$; finally, in the latter case singularities propagate along trajectories of the ``Liouvillian'' field of
$(V- \tau )/F$\,\footnote{\label{foot-18-11} Recall that
\begin{equation}
\frac {dx} {dt}=\mu ^{-1} (\nabla (V-\tau )/F))^\perp.
\label{18-2-7}
\end{equation}}.

We define $\beta _{jk}$ to have Poisson brackets vanish:
$\{Q_j, P_k\} \equiv 0 \mod O(\mu ^{-1})$, or simply
\begin{equation}
Q_1= x_1- F_{12}^{-1}P_2, \quad Q_2 = x_2 +F_{12}^{-1} P_1, \qquad \{P_1,P_2\}=-\mu F_{12}
\label{18-2-8}
\end{equation}
as we would do in the smooth case; then
\begin{equation}
\{ A_0, Q_j\} =\mu ^{-1} \sum _{ki} \alpha _{jki}P_kP_i;
\label{18-2-9}
\end{equation}
where $A_0=A-V$. In the smooth case we would be able to eliminate
$\alpha _{jkl}$ modulo $\omega _j A_0$ and $O(\mu ^{-2})$ by an appropriate
choice of coefficients $\beta '_{jki}$ above but now $\alpha_{jki}$ have
smoothness $(l-1,\sigma)$ and then $\beta '_{jki}$ will have the same smoothness which is not enough to plug them into Poisson brackets unless 
$(l,\sigma)\succeq (2,0)$. 

\begin{claim}\label{18-2-10}
As $(1,1)\preceq (l,\sigma) \preceq (2,0)$ let us replace coefficients 
$\beta '_{jki}$ by their $\mu^{-1}$-mollifications\footnote{\label{foot-18-12} Which makes sense only as $\mu^{-1}\ge \varepsilon $ i.e. 
$\mu \ge (h\nu^{-1}|\log h|)^{\frac{1}{2}}$.} which provides $\mu \vartheta (\mu^{-1})$-approximation to them. 
\end{claim}

Then for these new coefficients we have
$\{ A_0, Q'_j\}\equiv \mu ^{-1} \omega _j A^0$ and therefore
\begin{multline}
\{ A, Q'_ j\} \equiv
\mu ^{-1}\bigl(\omega _j (\tau +V) + F_{12}^{-1}\cL _j V -
\omega _j(\tau - A)\bigr) \\
\mod O\bigl(\vartheta (\mu^{-1})+\mu^{-2}\bigr)
\label{18-2-11}
\end{multline}
with $\cL _j \Def (-1)^{j-1}\partial _{x_{2-j}}$ (as $g^{jk}=\updelta_{jk}$ at the given point). 

In the smooth case we know from Chapter~\ref{book_new-sect-13} that $\omega _j = \cL _j F^{-1}$ and so it should be in our case as well:
\begin{equation}
\{ A, Q'_j\} \equiv \mu ^{-1} \Bigl( \cL _j \bigl({\frac {V - \tau } F}\bigr)+
\omega _j(\tau - A) \Bigr) \mod O\bigl(\vartheta (\mu^{-1})+\mu^{-2}\bigr)
\label{18-2-12}
\end{equation}
and this equality will be sufficient to prove that singularities propagate along
trajectories of (\ref{18-2-7}).

\begin{remark}\label{rem-18-2-2}
(i) The propagation speed of (\ref{18-2-7}) is $\asymp \nu \mu^{-1}$. The right-hand side of (\ref{18-2-12}) is larger than the error as 
\begin{equation}
\nu \ge C_1 \mu \vartheta (\mu^{-1})+C_1\mu^{-1}.
\label{18-2-13}
\end{equation}
In particular, as $(l,\sigma)\succeq (2,0)$ we need no $\mu^{-1}$-mollification and $\nu\ge C\mu^{-1}$.

\medskip\noindent
(ii) As there is a mollification with respect to $\varepsilon$ we also need to assume that 
\begin{equation}
\nu \ge C_1 \varepsilon^{-1} \vartheta (\varepsilon)+C_1\varepsilon
\label{18-2-14}
\end{equation}
but it is more important only as $\varepsilon \ge \mu^{-1}$ i.e. when 
$\mu \ge (h|\log h|)^{-\frac{1}{2}}$.

\medskip\noindent 
(iii) In what follows we assume that \underline{either} 
$(l,\sigma)\succeq (1,1)$ and non-degeneration condition
\begin{equation}
|\nabla V/F|\ge \epsilon_0
\label{18-2-15}
\end{equation} 
is fulfilled a \underline{or} $(l,\sigma)\succeq (2,0)$ and non-degeneracy assumption 
\begin{equation}
|\nabla V/F|+|\det \Hess V/F| \ge \epsilon_0
\label{18-2-16}
\end{equation} 
is fulfilled leaving to the reader problem~\ref{problem-18-2-3} below.
\end{remark}

\begin{Problem}\label{problem-18-2-3}
Investigate the case when $(1,1)\prec (l,\sigma)\prec (2,0)$ and assumption (\ref{18-2-16}) is replaced by the estimate to $\mes (\{x: |\nabla V/F|\le \nu\})$ as $\nu \to 0$; in this case $\epsilon \eta$-vicinity should be replaced by $\gamma$-vicinity with $\gamma^{-1}\vartheta (\gamma)=\epsilon \nu$. Here $\gamma \ge C_2\mu^{-1}$ due to (\ref{18-2-13}).
\end{Problem}

\section{Propagation of singularities: rigorous results}
\label{book_new-sect-18-2-2}

First let us prove the finite speed of the propagation:

\begin{proposition} \label{prop-18-2-4} 
Let $d\ge 2$, $\mu \le h^{-1+\delta }$ and let
\begin{equation}
M \ge \sup _{x:\sum _{jk}g^{jk}\xi _j \xi _k +V = 0} 2\sum _k g^{jk}\xi _k
+\epsilon\qquad \forall j
\label{18-2-17}
\end{equation}
with arbitrarily small constant $\epsilon >0$.

Let $\phi _1$ be supported in $B(0,1)$, $\phi _2=1$ in $B(0,2)$, $\chi $ be
supported in $[-1,1]$\,\footnote{\label{foot-18-13} Recall that all such auxiliary functions are appropriate in the sense of section~\ref{book_new-sect-2-3}.}. Let $\bar{T}_* = Ch |\log h| \le T \le T^{*\prime}=\epsilon _0$.

Then
\begin{multline}
|F_{t\to h^{-1}\tau } \chi _T(t) \bigl(1- \phi _{2, MT} (x-\bar{x})\bigr)
\phi_{1, MT} (y-\bar{x}) U (x,y,t) | \le Ch^s \\[2pt]
\forall \tau \le \epsilon_1
\label{18-2-18}
\end{multline}
where here and below $\epsilon_1>0 $ is a small enough constant.
\end{proposition}

\begin{remark}\label{rem-18-2-5}
Surely results will be much more precise in the distance according to $g_{jk}V^{-1}$ metrics.
\end{remark}

\begin{proof}
Proof follows the proof of Theorem~\ref{book_new-thm-2-3-2}. Namely, let $\upchi$ be the same function as there; recall that it was supported in $(-\infty, 0]$ and satisfying certain ``regularity'' conditions. Then as a main auxiliary symbol we pick up
\begin{multline}
\upchi \Bigl(f(x,t) -\epsilon _3\Bigr)\qquad \text{with}\quad
f(x,t)=\frac{1}{T}\bigl(|x-\bar{x}|^2+\epsilon _2\bigr)^{\frac{1}{2}}
\mp M \frac{t}{T},
\label{18-2-19}
\end{multline}
We pick signs ``$\mp$'' analyzing cases $\pm t >0$.
\end{proof}

Now let us investigate magnetic drift. 

\begin{proposition} \label{prop-18-2-6} 
Let $d=2$ and condition \textup{(\ref{18-1-6})} be fulfilled. Let 
$\mu _0\le \mu \le h^{-1+\delta }$ where here and below $\mu _0>0$, $C_1$, and $1/\epsilon_1$ are large enough constants depending on all the other constants, exponents and the choice of axillary functions below. Let \underline{either} $(l,\sigma)\succeq (1,1)$ and $\nu=1$ \underline{or} $(l,\sigma)\succeq (2,0)$ and 
\begin{equation}
\nu \ge C_0 \max(\mu^{-1},\varepsilon).
\label{18-2-20}
\end{equation}
Consider ball $B(\bar{x},\epsilon \nu)$.

\medskip\noindent
(i) Assume that $|\cL _{V/F}(\bar{y})|\le \nu $. Let
\begin{equation}
T_* \Def C \varepsilon ^{-1}\nu^{-1} h|\log h| + 
C\nu^{-1} (\mu h |\log h|)^{\frac{1}{2}} \le T \le T^* = \epsilon _1\mu.
\label{18-2-21}
\end{equation}
Then\footnote{\label{foot-18-14} Recall that $b^\w$ means Weyl quantization of symbol $b$ and due to condition (\ref{18-2-21}) logarithmic uncertainty principle holds and this quantization of the symbols involved is justified, $\,^t\!B$ means the dual operator.} 
\begin{multline}
|F_{t\to h^{-1}\tau } \bar{\chi} _T(t)
\bigl(1-\phi _{2, \mu ^{-1}\nu T} (Q_1-\bar{x}_1, Q_2-\bar{x}_2) \bigr)^\w U(x,y,t)\times
\\[3pt]
\, ^t\!\bigl(\phi _{1, \mu ^{-1}\nu T} (Q_1-\bar{x}_1, Q_2-\bar{x}_2) \bigr)^\w|
\le Ch^s\qquad \forall \tau , |\tau |\le \epsilon_1\nu;
\label{18-2-22}
\end{multline}
(ii) Assume now that $|\cL _{V/F}(\bar{y})|\asymp \nu$. Let $\chi $ be supported in $[1-\epsilon _0, 1+\epsilon_0]$.

Further, let $\bar{x} = \bar{y} + \mu ^{-1} \cL _{V/F}(\bar{y})T$ with $T$
satisfying \textup{(\ref{18-2-21})} and $T\le \epsilon_1 \mu$.

Then 
\begin{multline}
|F_{t\to h^{-1}\tau } \chi _T(t)
\bigl(1-\phi _{2, M\epsilon \mu ^{-1}\nu T} 
(Q_1-\bar{x}_1, Q_2-\bar{x}_2)\bigr)^\w U(x,y,t)\times
\\[3pt] 
^t\!\bigl(\phi_{1,\epsilon \mu ^{-1}T}(Q_1-\bar{y}_1, Q_2-\bar{y}_2)\bigr)^\w|
\le C h^s\qquad \forall \tau , |\tau |\le \epsilon_1\nu.
\label{18-2-23}
\end{multline}
\end{proposition}

\begin{proof}
Proofs of both statements follow the proof of Theorem~\ref{book_new-thm-2-3-2}. Namely, let $\upchi$ be the same function as there; recall that it was supported in $(-\infty, 0]$ and satisfying certain ``regularity'' conditions. Then as a main auxiliary symbol we pick up
\begin{multline}
\upchi \Bigl(f\bigl( Q_1(x,\xi), Q_2(x,\xi),t\bigr)-\epsilon _2\Bigr) \qquad
\text{with}
\\
f(Q,t)=\frac{1}{T}
\bigl(\mu^2\nu^{-2}|Q-{\bar Q}|^2+\epsilon _3\bigr)^{\frac{1}{2}} \mp 
M \frac{t}{T}
\label{18-2-24}
\end{multline}
and
\begin{multline}
\upchi \Bigl(f \bigl( Q'_1(x,\xi), Q'_2(x,\xi),t\bigr)-\epsilon _2\Bigr) \qquad
\text{with} \\
f(Q',t)=\frac{1}{T}
\bigl(\mu^2\nu^{-2}|Q'-\bar{Q}'-\mu ^{-1}\ell t|^2+
\epsilon _3\bigr)^{\frac{1}{2}} \pm \epsilon_4 \frac{t}{T}
\label{18-2-25}
\end{multline}
in the proofs of statements (i) and (ii) respectively where $\epsilon _k>0$ are small constants. We pick signs ``$\mp$'' analyzing cases $\pm t >0$. Here 
$\ell = \cL _{V/F}(\bar{y})$ in the proof of statement (ii). 
\end{proof}

The following corollary follows immediately from proposition \ref{prop-18-2-6}(ii):

\begin{corollary} \label{cor-18-2-7} 
(i) Statement (ii) of proposition \ref{prop-18-2-6} remains true for all
$T$ satisfying \textup{(\ref{18-2-21})} if we redefine $\bar{x}$ as
$\bar{x}= \Psi _T (\bar{y})$ where $\Psi _t$ is Liouvillian flow defined by
\textup{(\ref{18-2-7})}.

\medskip\noindent
(ii) As $C_0\nu^{-1} \le T \le \epsilon_1 \mu$ 
\begin{equation}
|F_{t\to h^{-1}\tau } \Bigl( \chi _T(t) \Gamma_x \bigl(\psi (x) U\bigr)\Bigr)
|\le Ch^s
\qquad \forall \tau: |\tau |\le \epsilon_1\nu .
\label{18-2-26}
\end{equation}
Recall that $\chi $ is supported in $[-1,-\frac{1}{2}]\cup[\frac{1}{2},1]$.
\end{corollary}

Recall that under condition (\ref{book_new-13-3-45}) i.e. 
\begin{equation}
V\le -\epsilon_0
\label{18-2-27}
\end{equation}
estimate 
\begin{equation}
|F_{t\to h^{-1}\tau }\Bigl(\chi _T(t)\Gamma\bigl(\psi (x)U\bigr)\Bigr)| \le Ch^s
\qquad \forall \tau: |\tau |\le \epsilon_1.
\label{18-2-28}
\end{equation}
holds for $Ch|\log h|\le T\le \bar{T}= \epsilon_1\mu^{-1}$. Note that $T_*\le \bar{T}$ iff 
\begin{equation}
\nu \ge C(\mu^3h|\log h|)^{\frac{1}{2}}\label{18-2-29}
\end{equation}
and $\varepsilon$ is defined by (\ref{18-2-5}) where we picked the smallest possible $\varepsilon$. 

Therefore 
\begin{claim}\label{18-2-30}
Under assumptions (\ref{18-2-27}), (\ref{18-2-29})\,\footnote{\label{foot-18-15} Which we replace by \ref{18-2-29-'}.} and (\ref{18-2-5}) estimate (\ref{18-2-28}) holds for $T_*\Def Ch|\log h|\le T\le T^*\Def \epsilon _1\mu$. 
\end{claim}

Out of two conditions (\ref{18-2-29}) and (\ref{18-2-5}) the first is the nastiest. To weaken it we need a special analysis; however note that $\bar{T}\le \epsilon \nu^{-1}$ so we can restrict ourselves to $T\le \epsilon \nu^{-1}$ when evolution is restricted to $C_0\mu^{-1}$ vicinity of the original point. 

\begin{proposition}\label{prop-18-2-8}
Let assumptions of proposition~\ref{prop-18-2-6} be fulfilled.

\medskip\noindent
(i) Let $g^{jk}=\updelta_{jk}$, $F=1$. Then \textup{(\ref{18-2-28})} holds for $T_*\le T\le T^*$ with
\begin{equation}
T_*\Def C\varepsilon ^{-1}\nu^{-1}h|\log h|.
\label{18-2-31}
\end{equation}
(ii) In the general case \textup{(\ref{18-2-28})} holds for $T_*\le T\le T^*$ with
\begin{equation}
T_*\Def C\varepsilon ^{-1}\nu^{-1}h|\log h|+ C\nu^{-1} \bigl(\bar{\vartheta}(\mu^{-1}) h|\log h|\bigr)^{\frac{1}{2}}
\label{18-2-32}
\end{equation}
where $\bar{\vartheta}$ is a modulus of continuity of $g^{jk},F$.
\end{proposition}

\begin{proof}
(i) In this case we can assume without any loss of the generality that 
$A_0= h^2D_1^2 + (hD_2-\mu x_1)^2$. Let us consider propagation with respect to $\xi_2$. We can use then function $f=T^{-1}\mu^{-1}\nu^{-1}\xi_2$ and as for time $T$ shift with respect to $\xi_2$ is $\asymp \nu T $ and logarithmic uncertainty principle means exactly that $T\ge T_*$ with $T_*$ defined by (\ref{18-2-31}).

\medskip\noindent
(ii) Consider the general case. Then without any loss of the generality we can assume that $F=1$ (otherwise we can divide by it) and $g^{jk}(\bar{x})=\updelta_{jk}$, $\nabla g^{jk}(\bar{x})=0$, $V_1(\bar{x})=V_2(\bar{x})=0$, $\nabla V_1(\bar{x})=0$. In this case 
$\mu Q_1= \xi_2 + \mu \phi (x)$ with 
$\nabla \phi = O\bigl(\mu \bar{\vartheta}(\mu^{-1})\bigr)$ and quantization conditions is $T\ge T_*$ with $T_*$ defined by (\ref{18-2-32}).
\end{proof}

Now condition $T_*\le \bar{T}$ becomes (\ref{18-2-5}) plus 
\begin{equation}
\nu \ge C\bigl(\mu^3 \bar{\vartheta}(\mu^{-1})h|\log h|\bigr)^{\frac{1}{2}}.
\tag*{$\textup{(\ref*{18-2-29})}'$}\label{18-2-29-'}
\end{equation}

\section{Main theorem}
\label{book_new-sect-18-2-3}

Now we can prove main results of this section. Under assumption (\ref{18-2-15}) we pick up $\varepsilon=C\mu h|\log h|$ provided $\nu=1$ satisfies \ref{18-2-29-'} i.e. 
\begin{equation}
C\mu^3\bar{\vartheta} (\mu^{-1}) h|\log h|\le 1;
\label{18-2-33}
\end{equation}
as $(\bar{l},\bar{\sigma})=(1,1)$ it is equivalent to 
$\mu \le \epsilon_1h^{-\frac{1}{2}}$ and as $(\bar{l},\bar{\sigma})=(2,0)$ it is no restriction at all.

Under assumptions (\ref{18-2-16}) (and thus $(l,\sigma)\succeq (2,0)$), (\ref{18-2-33}) we introduce variable
\begin{equation}
\nu =\nu (x)= \epsilon |\nabla V/F| + \bar{\nu},\qquad
\bar{\nu}\Def C_0(\mu h|\log h|)^{\frac{1}{2}}+C_0\mu^{-1}
\label{18-2-34}
\end{equation}
and variable $\varepsilon=\varepsilon (x)$ by (\ref{18-2-5}) where $\bar{\nu}$ ensures that (\ref{18-2-20}) and \ref{18-2-29-'} are fulfilled.

Let $\tilde{A}$ be $\varepsilon$-approximation\footnote{\label{foot-18-16} One can take two copies $\tilde{A}^\pm$ of approximation such that the corresponding estimate  (\ref{18-2-35}) or (\ref{18-2-37})
holds for both and $\tilde{A} ^-\le A\le \tilde{A} ^+$ in the operator
sense. We call them \underbar{framing approximations}. This remark holds for
all our theorems \ref{thm-18-2-9}, \ref{thm-18-6-10}, \ref{thm-18-9-4}, \ref{thm-18-9-6} in full.} of $A$.

\begin{theorem} \label{thm-18-2-9} 
Let $d=2$ and $A$ be a self-adjoint in $\sL ^2(X)$ operator defined by \textup{(\ref{18-1-1})}. Let conditions \textup{(\ref{18-1-6})} and \textup{(\ref{18-2-27})} be fulfilled in $B(0,1)\subset X\subset \bR ^2$. 

\medskip\noindent
(i) Further, let $(l,\sigma)\succeq (1,1)$, condition \textup{(\ref{18-2-33})} be fulfilled and condition \textup{(\ref{18-2-15})} be fulfilled in $B(0,1)$. Then estimate
\begin{multline}
\R^\MW \Def |\int \Bigl( \tilde{e}(x,x,0)-
h^{-2}\cN _{2} ^\MW (x,0)\Bigr)\psi (x)dx|\le\\ 
C\mu ^{-1}h ^{-1}+ Ch^{-2}\vartheta (\mu h|\log h|)
\label{18-2-35}
\end{multline}
holds; this estimate also holds with the standard Weyl expression albeit with infinite\footnote{\label{foot-18-17} Actually finite, due to finite smoothness.} number of terms generated by the principal part:
\begin{multline}
\R^\W_\infty \Def |\int \Bigl( \tilde{e}(x,x,0)-
h^{-2}\cN _{2,\infty} ^\W (x,0)\Bigr)\psi (x)dx|\le\\ 
C\mu ^{-1}h ^{-1}+ Ch^{-2}\vartheta (\mu h|\log h|)
\label{18-2-36}
\end{multline}

\medskip\noindent
(ii) On the other hand, let $(l,\sigma)\succeq (2,1)$, condition \textup{(\ref{18-2-33})} be fulfilled and condition \textup{(\ref{18-2-16})} be fulfilled in $B(0,1)$. Then estimate
\begin{equation}
\R^\MW \le C\mu ^{-1}h ^{-1}+ C\mu ^2 |\log h|
\label{18-2-37}
\end{equation}
holds; this estimate also holds with the standard Weyl expression $\R^\W$.
\end{theorem}

\begin{proof}
(i) Due to proposition \ref{prop-18-2-8} estimate (\ref{18-1-20}) holds 
with $T=T^*=\epsilon \mu $ and $T_*= Ch|\log h|$ and therefore due to theory of section~\ref{book_new-sect-4-5}
\begin{equation}
|F_{t\to h ^{-1}\tau }\Bigl(\bar{\chi}_T(t) \Gamma \bigl( \psi (x)U\bigr)
\Bigr) | \le Ch^{-1}.
\label{18-2-38}
\end{equation}
Therefore due to the standard Tauberian procedure
\begin{multline}
|\Gamma \bigl(\tilde{e} (.,.,0)\psi\bigr) -
h^{-1}\int _{-\infty}^0 \Bigl( F_{t\to h ^{-1}\tau }
\bigl( \bar{\chi}_T(t) \Gamma (U \psi )\bigr)\Bigr)\, d\tau | \le\\
\frac {C} {T^*} h^{-1}\asymp C\mu^{-1}h^{-1} 
\label{18-2-39}
\end{multline}
with $T=T^*$ and due to (\ref{18-1-20}) with $T=Ch|\log h|$; but then 
\begin{multline}
|\int \Bigl( \tilde{e}(x,x,0)-
h^{-2}\tilde{\cN} _{2} ^\MW (x,0)\Bigr)\psi (x)dx|\le 
C\mu ^{-1}h ^{-1}.
\label{18-2-40}
\end{multline}
Finally, one can see easily that replacing $\tilde{\cN}^\W_{2,\infty}$ by $\cN^\MW_{2}$ we make an error not exceeding 
$Ch^{-2}\vartheta (\varepsilon)$.

In these arguments one can replace $\cN^\MW_{2}$ by $\cN^\W_{2,\infty}$.

\medskip\noindent
(ii) The same arguments work in the \emph{exterior zone\/}\index{zone!exterior} 
\begin{equation}
\cX_\ext= \{x: |(\nabla V/F)(x) |\ge C\bar{\nu}\};
\label{18-2-41}
\end{equation}
therefore contribution of this zone to both the Tauberian remainder and to the left-hand expression (\ref{18-2-40}) do not exceed $C\mu^{-1}h^{-1}$. 

Further, contributions of the \emph{interior zone\/}\index{zone!interior}
\begin{equation}
\cX_\inter=\{x: |(\nabla V/F)(x) |\le C\bar{\nu}
\label{18-2-42}
\end{equation}
to both the Tauberian remainder and to the left-hand expression (\ref{18-2-40}) do not exceed 
$C\mu h^{-1}\mes (\cX') \asymp C\mu h^{-1}\bar{\nu}^2\asymp C\mu^{-1}h^{-1}+ C\mu ^2|\log h| $.

Finally, one can see easily that replacing $\tilde{\cN}^\MW_{2}$ by $\cN^\MW_{2}$ we make an error not exceeding 
\begin{equation*}
Ch^{-2}\int_{\cX_\ext} \vartheta \bigl(\mu h|\log h|/\nu^2\bigr)\,\nu d\nu +
Ch^{-2}\vartheta (\bar{\varepsilon}) \mes (\cX') 
\end{equation*}
with $\bar{\varepsilon}=C_0\mu h \bar{\nu}^{-1} |\log h|\asymp 
C_0\min \bigl( (\mu h|\log h|)^{\frac{1}{2}}, \mu^2h|\log h|\bigr)$;
as $(l,\sigma)\succeq (2,1)$ we arrive to estimate (\ref{18-2-37}).
\end{proof}

In these arguments one can replace $\cN^\MW_{2}$ by $\cN^\W_{2,\infty}$ but all extra terms do not exceed $C\mu^2$.

\begin{corollary}\label{cor-18-2-10}
In the framework of theorem~\ref{thm-18-2-9} let $(l,\sigma)\succeq (2,1)$ and one of conditions~\textup{(\ref{18-2-15})}, \textup{(\ref{18-2-16})} be fulfilled. Then estimate 
\begin{equation}
\R^\MW = O(\mu^{-1}h^{-1})
\label{18-2-43}
\end{equation}
holds under assumption \textup{(\ref{18-2-6})} 
(i.e. $\mu \le \bar{\mu}_1= (h|\log h|)^{\frac{1}{3}}$). The same estimate holds for $\R^\W$.
\end{corollary}

\begin{remark}\label{rem-18-2-11}
Under stronger smoothness assumptions this estimate (\ref{18-2-43}) holds for larger $\mu$ under assumption (\ref{18-2-15}) but it does not hold for $\R^\W$ under assumption (\ref{18-2-16}).
\end{remark}

\chapter{$d=2$: Canonical form}
\label{book_new-sect-18-3}

Now we analyze the case
\begin{equation}
\bar{\mu} _1 \Def (h|\log h|)^{-\frac{1}{3}}
\le \mu \le \bar{\mu}_2\Def (h|\log h|)^{-1};
\label{18-3-1}
\end{equation}
the second restriction will be removed in subsection~\ref{book_new-sect-18-3-3}.

Our arguments will not change much for larger $\mu $ as well but some of them
become much simpler while other become a bit more complicated. We will take
\begin{equation}
\mu^{-1}\ge \varepsilon \ge C(\mu ^{-1}h |\log h|)^{\frac{1}{2}}
\label{18-3-2}
\end{equation}
which makes it possible to use $\mu^{-1}h$-pseudo-differential operators.

\section{Reduction}
\label{book_new-sect-18-3-1}

\subsection{Main part of operator. I}
\label{book_new-sect-18-3-1-1}

Now our goal is to reduce our operator to a some kind of the canonical form. It would be easier for constant or at least smooth $g^{jk}$ and $F$. However, we can overcome technical difficulties in a more general case as well. Instead of normal Fourier integral operators we will use ``poor man Fourier integral operators'' i.e. propagators of ``rough'' pseudo-differential operators.

\begin{problem}\label{problem-18-3-1}
In this section for a sake of simplicity of some arguments we assume that $\sigma \ge 0$ but one can get rid of this assumption easily. We leave it to the reader.
\end{problem}

In this subsection we will reduce $A_0$ to pre-canonical form which is really easy for constant metrics $g^{jk}$ and $F$ and is not needed at all for Euclidean metrics and constant $F$.

First of all, we can assume without any loss of the generality that
\begin{equation}
F^3= F=1.
\label{18-3-3}
\end{equation}
Really, changing orientation, if necessary, we can make $F^3>0$. Then
multiplying mollified operator by $F^{-{\frac{1}{2}}}$ both from the right and from the left and commuting it with $P_j$ we will get operator of the same form with $g^{jk}$, $V$ replaced by $F^{-1}g^{jk}$, $F^{-1}V$ respectively, modulo
operator $h^2 V'(x)$ where $V'(x)$ is a function, linear with respect to second
derivatives of $g^{jk}$ and $F$ and thus $V'$ does not exceed
$C\bigl(1+\vartheta (\varepsilon )\varepsilon ^{-2}\bigr)$; one
can see easily that after multiplication by $h^2$ it does not exceed
$C\vartheta (\varepsilon )$ which is an approximation error anyway (plus $Ch^2$
which is less than $C\mu ^{-1}h$). Surely, such transformation would affect
$e(x,x,0)$ but we will return to the original operator later.

In Chapter~\ref{book_new-sect-13} on this step of reduction we applied $\mu^{-1}h$-Fourier integral operator and let us try the same now. Namely, let us transform our operator by operator $T(1)$, where
\begin{equation}
T(t) \Def e^{-i\mu ^{-1}h^{-1}tL^\w}, \qquad
L^\w= {\frac{1}{2}}\sum_{jk}P_j L ^{jk}(x) P_k
\label{18-3-4}
\end{equation}
with $L^{jk} \in \sF ^{\bar{l} ,\bar{\sigma} }$. Obviously $L^\w$ is a Weyl
$\mu^{-1}h$-quantization of the symbol
\begin{equation}
L (x,\xi) \Def \mu ^2 \sum_{j,k} L^{jk}(x) p_j(x,\xi) p_k(x,\xi),
\qquad p_j\Def \xi _j -V_j(x).
\label{18-3-5}
\end{equation}

\begin{remark}\label{rem-18-3-2}
In this subsection $L,T,t$ do not denote the same things as everywhere else and
$\sF ^{l,\sigma}$ denotes space of functions functions satisfying
\begin{equation*}
|\partial _{x,\xi }^\nu f|\le
C_\nu \bigl(\varepsilon ^{l-|\nu|}|\log \varepsilon |^{-\sigma} +1\bigr)
\end{equation*}
with the same constants $C_\nu$ as in section~\ref{book_new-sect-2-3}.
\end{remark}

Let us consider ``Heisenberg evolution'' $Q(t)=T(-t){\bar Q}T(t)$ for
$\mu^{-1}h$-pseudo-differential operator ${\bar Q}$; then
\begin{equation}
\partial _t Q(t) = i\mu ^{-1}h^{-1}[L^\w,Q(t)], \qquad Q(0)={\bar Q}.
\label{18-3-6}
\end{equation}
Let us define first symbol $q(t)$ as a solution to
\begin{equation}
\partial _t q(t) = \{\mu ^{-2}L ,q\}, \qquad q(0)=\bar{q}\implies
q(t) = \bar{q}\circ \Phi _t
\label{18-3-7}
\end{equation}
with the standard Poison brackets where $\phi_t$ is a corresponding Hamiltonian
flow. Let us consider differential equations defining Hamiltonian flow $\phi_t$ in terms of $p_j$ and $x_k$:
\begin{equation}
{\frac{d\ }{dt}} x_j= \sum_k L^{jk}p_k, \qquad
{\frac{d\ }{dt}}p_j= \sum_k \Lambda ^{jk}p_k + \sum _{k, m} \beta ^{km}_jp_k p_m
\label{18-3-8}
\end{equation}
with $\Lambda \Def f(x) \ell J$, where $\Lambda =\bigl(\Lambda ^{jk}\bigr)$,
$\ell =\bigl(L^{jk}\bigr)$ and $J=\left(\begin{smallmatrix} 0& -1\\ 1& 0\end{smallmatrix}\right)$ are $2\times 2$-matrices, here and below
$\beta_* \in \sF ^{\bar{l}-1,\bar{\sigma}}$,
$f(x)= \partial_{x_1}V_2-\partial_{x_2}V_1$.

One can see easily that for $|t|\le c$
\begin{gather}
|p|\le c \implies |p\circ \phi _t|\le C|p|,\label{18-3-9}\\[2pt]
\phi_t \in \sF^{ \bar{l}, \bar{\sigma}},\label{18-3-10}\\[2pt]
p_j\circ \phi _t = e^{\Lambda (x)t}p+\sum _{j,k}\beta '_{jk}p_jp_k
\qquad \text{on\ \ } \{|p|\le c\mu ^{-1}\} ,\label{18-3-11}\\
x_j\circ \phi_t - x_j = \sum_j \beta '' _{jk} p_k\label{18-3-12}
\end{gather}
where $p=(p_1,p_2)$, $x=(x_1,x_2)$.

Further, one can see easily that
\begin{multline}
a_0\circ \phi_t = \sum g_t ^{jk}p_jp_k + M_t, \quad
(g_t)=e^{t \,^t\!\Lambda (x) }(g)e^{t\Lambda (x)}\\
\text{with\ \ }M_t=\sum_{i,j,k}\beta _{ijk;\,t} p_ip_jp_k
\label{18-3-13}
\end{multline}
where $a_0=\sum_{j,k}g^{jk}p_jp_k$ and $(g)$, $(g_t)$ denote corresponding
matrices.

Now, let us consider in this region $\{|p|\le c\mu ^{-1}\}$ the third and higher
derivatives of the symbol $\alpha_{jk} p_jp_k$ with 
$\alpha_* \in \sF^{\bar{l},\bar{\sigma}}$; then
\begin{equation}
|\partial _{x,\xi } ^\nu (\alpha_{jk} p_jp_k)| \le
C_\nu \varepsilon ^{-|\nu| +l}|\log \mu |^{-\bar{\sigma}} \mu ^{-2}+ C
\label{18-3-14}
\end{equation}
for $3\le |\nu| \le c|\log \mu |$ where $C_\nu $ are standard constants of section~\ref{book_new-sect-3-2} even if $V_j\in \sF^{\bar{l},\bar{\sigma}}$ only. 

Similarly for $|\nu| \ge 3$ we have
\begin{equation}
|\partial _{x,\xi } ^\nu (\beta_{ijk} p_ip_jp_k)| \le
C_\nu \varepsilon ^{-|\nu| +\bar{l}-1}|\log \mu |^{-\bar{\sigma}} \mu ^{-3}+C;
\label{18-3-15}
\end{equation}
and therefore $\mu^2 M_t\in \sF^{\bar{l} ,\bar{\sigma}}\cap
\mu ^{-1}\sF^{\bar{l}-1,\bar{\sigma} }$.

Now let us pass from symbols to their quantizations. Let us recall that
$(\bar{l},\bar{\sigma})\succeq (2,0)$ and also that in the standard smooth situations and for the standard quantizations $[a^\w,b^\w]= -ih\{a,b\}^\w + O(h^3)$; then for the standard quantization of $\mu^2 a_0\circ \phi _t$ equation
(\ref{18-3-6}) will be fulfilled modulo operators with $A$-bound\footnote{\label{foot-18-18} Here in contrast to the standard definition we use everywhere we say that $A$ bound of $B$ is $R$ if
\begin{equation}
\|Bu\|\le R \bigl(\|u\|+\|A^m u\|\bigr)
\label{18-3-16}
\end{equation}
with some exponent $m$ which is fixed in advance.}
not exceeding $C\mu ^2 (\mu^{-1}h)^2 \times \varepsilon ^{-1}\mu ^{-1} \times
(\varepsilon ^{-1}\mu ^{-1} + \varepsilon ^{-2}\mu ^{-3})$
which in turn does not exceed $C\mu ^{-1}h$.

Therefore, 
\begin{claim}\label{18-3-17}
Transformed operator $T(-t)A_0T(t)$ differs from the quantization of the transformed symbol $\mu^2 a_0\circ \phi_t$ by an operator with $A$-bound not exceeding $C\mu^{-1}h$.
\end{claim}

Recall that $A_0$ is a Weyl $\mu^{-1}h$-quantization of $a_0$. This is as good as we need.

Now we can pick up $L^{jk}\in \sF^{\bar{l} ,\bar{\sigma}}$ such that
$e^{^t\!\Lambda } (g) e^\Lambda =\left(\begin{smallmatrix} 1& 0\\ 0& g^{-1} \end{smallmatrix}\right) $ with
$g=\det (g)^{-1}$ because $e^\Lambda $ could be any matrix with determinant
equal to 1. Therefore we get that modulo operator with $A$-bound not
exceeding $\mu^{-1}h$
\begin{gather}
T(-1)A_0T(1) \equiv P_1^2 +P_2(g^{-1})^\w P_2 + \mu ^2 M^\w, \label{18-3-18}\\[3pt]
M \Def a_0 \circ \phi _1 -p_1^2-g^{-1}p_2^2,\qquad
g'\Def g\circ \phi_1.\label{18-3-19}
\end{gather}
Recall that according to (\ref{18-3-13}) $M=\sum _{i,j,k}\beta _{ijk}p_ip_jp_k$ with $\beta _{ijk}\in \sF^{\bar{l}-1,\bar{\sigma}}$.
\medskip

Consider now how the above transformation affects $V$. Note that
\begin{equation*}
|\partial_{x,\xi}^\nu (V\circ \phi_t)|\le
C_\nu \varepsilon ^{l-|\nu|}|\log h|^{-\sigma}
\end{equation*}
in the same region $\{|p|\le c\mu ^{-1}\}$. Then one can see easily that 
\begin{equation}
T(-1)VT(1)\equiv (V\circ \phi_1)^\w
\label{18-3-20}
\end{equation}
modulo operator with $A$-bound not exceeding 
$C(\mu ^{-1}h)^2 \times \varepsilon ^{-1}\mu ^{-1}\times
\varepsilon ^{l-3}|\log h|^{-\sigma}$ (which is obviously less than
$C\vartheta (\varepsilon )$).

Also note that we used only that $V_1,V_2\in \sF^{\bar{l},\bar{\sigma}}$. Therefore without any loss of the generality one can assume that
\begin{equation}
V_1=0,\quad P_1=hD_1, \quad
V_2,\partial_{x_1} V_2=g^{-{\frac{1}{2}}}\in \sF^{\bar{l},\bar{\sigma}}
\label{18-3-21}
\end{equation}
where the third assertion follows from the first and the second ones and
\begin{equation}
p_2 = \alpha (x,\xi _2)\bigl(x_1 - \lambda (x_2,\xi _2)\bigr),
\qquad \alpha, \lambda \in \sF^{\bar{l},\bar{\sigma}}.
\label{18-3-22}
\end{equation}

\subsection{Main part of operator. II}
\label{book_new-sect-18-3-1-2}

In this subsection we will reduce $A_0$ to canonical form which is really easy for constant $g^{jk}$ and $F$.

Our next transformation is $T'(1)$ with
\begin{equation}
T'(t)=e^{-it\lambda ^\w D_1}\qquad
\text{where\ \ }
\lambda ^\w =\lambda ^\w(x_2,\mu ^{-1}hD_2);
\label{18-3-23}
\end{equation}
then
\begin{multline}
T'(-t)D_1T'(t)=D_1,\quad T'(-t)\lambda ^\w T'(t)=\lambda ^\w, \\
T'(-t)x_1T'(t)=x_1+t \lambda ^\w
\label{18-3-24}
\end{multline}
precisely.

To calculate transformations of other operators we need to introduce
the corresponding Hamiltonian flow $\phi'_t$:
\begin{gather}
\xi_1=\const,\;\lambda (x_2,\xi_2)=\const,\; x_1= x_1(0)+t \lambda (x_2,\xi_2),
\label{18-3-25}\\[2pt]
{\frac{d\ }{dt}}x_2 = \xi _1 \bigl(\partial_{\xi_2}\lambda (x_2,\xi_2)\bigr),
\qquad
{\frac{d\ }{dt}}\xi_2 = -\xi _1 \bigl(\partial_{x_2}\lambda (x_2,\xi_2)\bigr).
\label{18-3-26}
\end{gather}

This flow is less regular than $\phi_t$:
\begin{equation}
|\partial_{x,\xi}^\nu \phi'_t|\le C_\nu \mu ^{-1}
\bigl(1+\varepsilon ^{\bar{l}-1-|\nu |}|\log \mu |^{-\bar{\sigma}}\bigr).
\label{18-3-27}
\end{equation}
However, using again above arguments we can estimate the difference between operators
$T'(-1)\alpha ^\w T'(1)$ and $(\alpha \circ \phi'_1)^\w$ with
$\alpha \in \sF^{\bar{l} ,\bar{\sigma}}$: its $A$-bound does not exceed
$(\mu ^{-1}h)^2 \times \varepsilon ^{-1}\mu ^{-1} \times
\varepsilon ^{-2}\mu ^{-1}$ which in turn does not exceed $C\mu^{-1}h$.
The same statement is true for $\mu ^{-1}\beta ^\w$ with
$\beta \in \sF^{\bar{l}-1,\bar{\sigma} }$.

Therefore, \emph{modulo operator with $\bar{A}_0$-bound not exceeding
$C\mu ^{-1}h$\/}
\begin{gather}
T'(-1)T(-1)A_0T(1)T'(1) \equiv \bar{A}_0 + \mu^2 (M')^\w,\label{18-3-28}\\
\shortintertext{where}
\bar{A}_0\Def \bar{P}_1^2 + \bar{P}_2^2 , \quad \bar{P}_1=hD_1,
\bar{P}_2=-\mu x_1,\label{18-3-29}\\
M'= a_0\circ \Phi - \bar{a}_0=
\sum _{i,j,k}\beta '_{ijk}\bar{p}_i \bar{p}_j \bar{p}_k\label{18-3-30}
\end{gather}
with $\bar{p}_1\Def \xi_1$, $\bar{p}_2\Def -x_1$,
$\Phi\Def \phi _1\circ \phi'_1$.

Similarly, \emph{modulo operator with $\bar{A}_0$-bound not exceeding $C\vartheta(\varepsilon)$\/}
\begin{equation}
T'(-1)T(-1)VT(1)T'(1) \equiv (V\circ \Phi)^\w.
\label{18-3-31}
\end{equation}

\subsection{Potential}
\label{book_new-sect-18-3-1-3}

Now we have $\mu ^{-1}h$-pseudo-differential operator $\bar{A}$, with the symbol
$\mu^2 \bar{a}_0 = \mu^2 (x_1^2+\xi_1^2)$, perturbed by $\mu ^2 (M')^\w$ and
$(V')^\w$ with $V'=(V\circ \Phi)^\w$. The first perturbing operator would vanish
for constant $g^{jk}$ and $F$ but the second one would still cause trouble even in this case because it's symbol depends on $x_1,\xi_1$ as well\footnote{\label{foot-18-19} In this case $\Phi $ is a linear symplectomorphism.}.

Our goal is to make a perturbation a $\mu^{-1}h$-pseudo-differential operator with the symbol depending on $x_2,\xi_2$, and $ r ^2\Def x_1^2+\xi_1^2$ only.

Let us consider the smooth case first as a pilot model; then we can decompose
$V\circ \Phi$ into asymptotic series with respect to $x_1,\xi_1$. To get rid of
linear terms we need to make a shift in $x_1,\xi_1$ of a magnitude $\mu ^{-2}$;
then the increment of $\mu^2 \bar{a}_0$ will absorb these terms. 

We can continue further\footnote{\label{foot-18-20} We cannot eliminate terms of even degrees $2m$ completely, but we can reduce them to 
$\mu^2 \omega _m(x_2,\xi_2) (x_1^2 +\xi _1^2)^m$; see Chapter~\ref{book_new-sect-13}}
but what we need really is an error in magnetic Weyl expression not exceeding
$C\bigl(\mu ^{-1}h +\vartheta (\varepsilon)\bigr)h^{-2}$. In $3$-dimensional case we will need much less than this (an error in operator $O(\mu ^{-2})$ will be almost sufficiently small) but now we should request $O(\mu ^{-1}h)$ error
which can be as large as $O(\mu ^{-4})$ now and to run things in a ``smooth''
manner we would need $l=4$ while our objective is to get the best estimate as
$l=2$\,\footnote{\label{foot-18-21} Actually we need $l=3$ because the ``weak magnetic field'' approach leads to $O\bigl(h^{-2}(\mu h|\log h|)^l\bigr)$ error and the ``strong magnetic field'' approach in the discussed above smooth version produces $O\bigl(h^{-2} \mu ^{-l}\bigr)$ error; threshold is as 
$\mu =h^{-{\frac{1}{2}}}|\log h|^{-{\frac{1}{2}}}$ and to keep error below $\mu ^{-1}h$ we need $(l,\sigma)= (3,\frac{3}{2})$ and we will proceed in this assumption in our next Chapter~19
devoted to the higher dimensions.}.

Let us return to the non-smooth case and try to repeat the above arguments. To make the shift described above we need to use transformation
$T''(-t)=e^{-i \mu h^{-1}t S^\w}$ with $\mu^{-1}h$-pseudo-differential operator $S^\w$ with the symbol belonging to $\mu ^{-2}\sF^{l,\sigma}$; on the other hand, to accommodate $M'$ we will need to add $\mu^{-1}h$-pseudo-differential operator with the symbol satisfying improved estimate (\ref{18-3-15}). Namely, one can see easily that
\begin{equation}
|\partial _{x,\xi } ^\nu M' | \le
C_\nu \varepsilon ^{-|\nu| +1}|\log h|^{-1}\mu ^{-3}+
C\mu ^{(-3+\nu_1)_-}
\label{18-3-32}
\end{equation}
with $\nu_1$ counting derivatives with respect to $(x_1,\xi_1)$ only.

So in fact
\begin{equation}
S=S'+S'', \qquad S'\in \mu^{-2}\sF^{l,\sigma},\qquad
S''\in \mu^{-2}\sF^{1,0}
\label{18-3-33}
\end{equation}
and $S''$ satisfies (\ref{18-3-32}); some adjustments to be done later.

\begin{remark}\label{rem-18-3-3} 
One can check easily that such operator $T''(t)$ would be
$\mu^{-1}h$-pseudo-differential operator iff
\begin{equation}
\mu \ge \bar{\mu}^*_1 \Def Ch^{-\frac{1}{3}}|\log h|^{\frac{1}{3}}.
\label{18-3-34}
\end{equation}
It is much simpler to analyze transformation by $\mu^{-1}h$-pseudo-differential operators than by $\mu^{-1}h$-Fourier integral operators with not very regular symbols and one can avoid some hassle if (\ref{18-3-34}) holds. This restriction is not smoothness related: one just need to take away from $\mu^2 M'+V'$ a term of magnitude $\mu ^{-1}$.

One can cover case $\mu \le \bar{\mu}^*_1$ by results of the previous section;
this would require $(l,\sigma)= (\bar{l},\bar{\sigma})=(2,3)$ \underline{and} assumption (\ref{18-2-15}) only to achieve the best possible estimate $O(\mu^{-1}h)$. The second restriction would be the most unfortunate.

However we will use more sophisticated arguments to cover
$\mu \in [\bar{\mu}_1, \bar{\mu}^*_1]$ here and prove the best possible
estimate for $\sigma=\bar{\sigma}=1$ instead.
\end{remark}

Let us consider corresponding Hamiltonian flow $\psi _t$. One can see easily
that
\begin{equation}
|\partial_{x,\xi} ^\nu (\psi _t -I)|\le
C_\nu \mu ^{-2}\varepsilon ^{l-1-|\nu| }|\log h|^{-\sigma} +
C_\nu \mu^{-3}\varepsilon ^{-|\nu|}
\label{18-3-35}
\end{equation}
where $I$ is an identity map.

Using arguments of the previous subsection one can prove easily that 
$T''(-t){\bar P}_jT''(t)\equiv \mu ({\bar p}_j\circ \psi_t)^\w$ modulo operator
with the upper $\bar{A}_0$-bound
\begin{multline}
\mu (\mu ^{-1}h)^2 \times
\bigl(1+\mu ^{-2}\varepsilon ^{l-3}|\log \mu|^{-\sigma} +
\mu ^{-3} \varepsilon ^{-2}|\log \mu |^{-1}\bigr) \times \\[2pt]
\bigl(\mu ^{-2}\varepsilon ^{l-4}|\log h|^{-\sigma} +
\mu ^{-3} \varepsilon ^{-3}|\log \mu |^{-1}\bigr);
\label{18-3-36}
\end{multline}
then the same is true for
$T''(-t)\bar{A}_0T''(t)\equiv \mu^2 (\bar{a}_0\circ \psi_t)^\w$. 

The bad news is that this expression is not necessarily less than $C\mu^{-1}h$ even if $(l,\sigma)=(2,1)$; the trouble appears as $\mu $ is close to $h^{-\frac{1}{3}}$ in which case expression (\ref{18-3-36}) is close to $h$. Also, expression (\ref{18-3-36}) is not less than $C\vartheta (\varepsilon)$ for $l< \frac{3}{2}$ and $\mu $ close to $h^{- \frac{1}{3}} $ again.

We could avoid these problems by slightly increasing $\bar{l}$ and assuming
that $l \ge \frac{3}{2}$ but there is a better way. Note that (\ref{18-3-35}) with $\varepsilon $ replaced by $\bar{\epsilon}\Def \mu ^{-1-\delta }$ with small enough $\delta >0$ is less than 
$C\bigl(\mu ^{-1}h + \vartheta (\varepsilon)\bigr)$ (with the original parameter 
$\varepsilon$ in the last expression). Then let us assume that

\begin{claim}\label{18-3-37}
$S'\in \mu^{-2} \bar{\sF}^{l,\sigma}$ and $S''$ satisfies (\ref{18-3-32})
with $\varepsilon $ replaced by $\bar{\varepsilon}\Def \mu ^{-1-\delta}$;
here and below $\bar{\sF}^*$ means that in the definition of the class
$\varepsilon $ replaced by $\bar{\varepsilon}$ while $\sF^*$ denotes the
original class. 
\end{claim}
Then the same arguments show that
\begin{equation*}
T''(-t)(\mu ^2 M'+V')^\w T''(t)\equiv
\bigl((\mu ^2 M'+V')\circ \psi _t\bigr)^\w
\end{equation*}
with the same error. So, modulo operator with $A$-bound not exceeding 
$C\bigl(\mu^{-1}h+ \vartheta (\varepsilon)\bigr)$ we need to consider the quantization of
\begin{equation}
\bigl( \mu^2 \bar{a_0} +\mu ^2 M' +V'\bigr)\circ \psi _1.
\label{18-3-38}
\end{equation}
One can see easily that
$\bar{p}_j\circ \psi_t - \bar{p}_j- t\{ \bar{p}_j, S'+S''\} \in
\mu ^{-4}\bar{\sF}^{l-2,\sigma}$ and therefore
\begin{equation*}
\mu ^2 \bigl( \bar{a}_0\circ \psi_t - t \bar{a}_0- \{ \bar{a}_0, S'+S''\}\bigr) \in \mu ^{-2} \bar{\sF}^{l-1,\sigma};
\end{equation*}
we conclude that if $V'\in \bar{\sF}^{l,\sigma}$ and $M'$ satisfied (\ref{18-3-32}) with $\varepsilon $ replaced by $\bar{\varepsilon}$ then 
$(\mu ^2 M'+ V')\circ \psi _t - (\mu ^2 M'+ V')$ would belong to the same class.
In this case ``the main part'' of (\ref{18-3-38}) would be 
$\mu ^2 \bar{a}_0 + \mu^2 \{ \bar{a}_0, S'+S''\} + \mu ^2 M'+V'$ and we would define $S',S''$ from equations
\begin{gather}
\{ \bar{a}_0,S'\}= -\mu ^{-2}(V'-W'),\label{18-3-39}\\[2pt]
\{ \bar{a}_0,S''\}= -(M'-\mu ^{-2}W'').\label{18-3-40}
\end{gather}
Since $\{ \bar{a}_0, q\}=(\xi_1\partial _{x_1}-x_1\partial _{\xi_1})q$ these
equations would be solvable and solutions would satisfy (\ref{18-3-37}) if and only if
\begin{gather}
W' (x_2,\xi_2, r )=
\frac{1}{2\pi} \int V'( r \cos t , r \sin t; x_2,\xi _2)\,dt,\label{18-3-41}\\
W'' (x_2,\xi_2, r )=
\frac{1}{2\pi}\mu ^2 \int M'( r \cos t , r \sin t; x_2,\xi _2)\,dt;\label{18-3-42}
\end{gather}
where we recall that $ r ^2= x_1^2+\xi_1^2$.

One can calculate easily the difference between
$\bigl( \mu^2 \bar{a_0} +\mu ^2 M' +V'\bigr)\circ \psi _1$ and
$\mu ^2 \bar{a}_0 + \mu^2 \{ \bar{a}_0, S'+S''\} + \mu ^2 M'+V'$; it will be
$O(\mu^{-2})$ which is small enough to be taken care of by $\mu^{-1}h$-pseudo-differential operator transformations in what follows. More precisely
\begin{multline}
N\Def \Bigl( \mu^2 \bar{a}_0 +\mu ^2 M' +V'\Bigr)\circ \psi _1-\\
\Bigl(\mu ^2\bar{a}_0 + \mu^2 \{ \bar{a}_0, S'+S''\} + \mu ^2 M'+V'\Bigl)
\in \mu^{-2} \bar{\sF}^{l-1,\sigma}.
\label{18-3-43}
\end{multline}
Since our symbols belong to $\sF^*$ rather than $ \bar{\sF}^*$, we just
replace $V',M'$ by their $ \bar{\varepsilon}$-mollifications
$\bar{V}', \bar{M}'$ and define $\bar{S}',\bar{S}'',\bar{W}',\bar{W}''$ by
(\ref{18-3-39})--(\ref{18-3-42}) with $\bar{V}',\bar{M'}$ instead of $V',M'$.

Then we have unaccounted term $K^\w$ in the operator \emph{before\/} the last transformation with $K=(\mu^2 \tilde{M}'+ \tilde{V}')$, where 
$\tilde{M}'=M'- \bar{M}'$, $\tilde{V'}=V'-\bar{V}'$ and then
$K \in \vartheta( \bar{\varepsilon})\sF^{0,0}$
and it also belongs to the same class as $V'+\mu^2 M'$ did.

To finish this part we need to understand how operator $K^\w$ is transformed by
$T''(1)$. One can see easily that $K= K_0+K_1+\ldots+K_m$ where $K_j$ is the
difference between $\varepsilon_j$ and $\varepsilon_{j+1}$-mollifications,
$\varepsilon_j = 2^j\varepsilon$ and $\varepsilon_m =\bar{\varepsilon}$.
Obviously $K_j \in \vartheta (\varepsilon _j)\sF_j ^{0,0}$ where
\begin{claim}\label{18-3-44}
$\sF_j^*$ denotes corresponding class with $\varepsilon$ replaced by
$\varepsilon_j$.
\end{claim}

Then one can see easily that $T''(-1)K_j^\w T''(1) \equiv K_j^\w$ modulo
operator with the norm not exceeding
$\mu^{-2}\varepsilon_j^{l-2}|\log h|^{-\sigma}$ and then
$T''(-1)K^\w T''(1) \equiv K^\w$ modulo operator with the norm not exceeding
$\mu^{-2}\varepsilon ^{l-2}|\log h|^{-\sigma}$ as $l<2$ or
$\mu^{-2}|\log h|^{1-\sigma}$ as $l=2$ and in both cases we got
$O\bigl(\vartheta(\varepsilon)\bigr)$ estimate unless
$l=2, \mu \le \bar{\mu}^*_1$.

The similar calculations imply that 
$T''(-1)K^\w T''(1) \equiv (K')^\w$ modulo operator with the norm not exceeding
$C\vartheta (\varepsilon)$ with $K'=K+\{K,{\bar{S}}'+{\bar{S}}''\}$ as $l=2$ and $\mu\in [\bar{\mu} _1, \bar{\mu}^*_1]$.

\subsection{Final reduction}
\label{book_new-sect-18-3-1-4}

Thus we arrived to operator $\bar{A}_0 + \bar{W}^\w + (L+K')^\w$ where
$\bar{W}=W' + W''$, symbol $L$ is defined by (\ref{18-3-43}) and
\begin{equation}
K'=\sum _{0\le j\le m} K_j, \qquad K_j\in
\vartheta (\varepsilon _j)\sF_j ^{0,0}\subset \sF_j ^{l,\sigma}.
\label{18-3-45}
\end{equation}
Fortunately, both operators $L^\w$ and $K^{\prime\w}$ are small enough to be reduced by $\mu^{-1}h$-pseudo-differential operator.

Namely, let us consider transformation by
$U = \bigl(e^{-i\mu ^{-1}h^{-1}(S'''+S^{IV}}\bigr)^\w$
where $S'''$ and $S^{IV}$ are symbols of the same class as $L$ and $K$
respectively; in particular, $S^{IV}$ admits decomposition of type (\ref{18-3-45}).

One can check easily that
\begin{equation*}
e^{-i\mu ^{-1}h^{-1}(S'''+S^{IV})} \in \sF^{0,0}
\end{equation*}
as long as $\mu ^{-l}h^{-1}|\log h|^{-\sigma}\le C\varepsilon ^{-1}$ which is always the case because $\mu \ge C(h|\log h|)^{-\frac{1}{3}}$ and we assume that 
\begin{equation}
(l,\sigma)\succeq (1,2).
\label{18-3-46}
\end{equation}
Also one can see easily that modulo operator with $\bar{A}_0$-bound not exceeding $C\bigl(\mu^{-1}h + \vartheta (\varepsilon) \bigr)$
\begin{multline}
U^{-1}\Bigl( \bar{A}_0+\bigl( \bar{W}+ K' +L \bigr)^\w \Bigr) U \equiv\\
\bar{A}_0+\bigl( \{\bar{a}_0, S'''+S^{IV}\} + \bar{W}+ K +L +R \bigr)^\w
\label{18-3-47}
\end{multline}
where \emph{one needs to include $R\in \mu^{-4} \sF^{l-2,\sigma}$ only
for $(l,\sigma)\succeq (2,0)$\/}.

Now we need just to define $S''',S^{IV},W''',W^{IV}$ from equations similar to
(\ref{18-3-39})--(\ref{18-3-42}):
\begin{gather}
\{ \bar{a}_0, S'''\} = -(L-W'''),\label{18-3-48}\\[3pt]
\{ \bar{a}_0, S^{IV}\} = -(K'-W^{IV}),\label{18-3-49}\\[2pt]
W''' (x_2,\xi_2, r )=
{\frac{1}{2\pi}}\int L( r \cos t , r \sin t; x_2,\xi _2)\,dt,\label{18-3-50}\\
W^{IV} (x_2,\xi_2, r )=
{\frac{1}{2\pi}}\int K'( r \cos t , r \sin t; x_2,\xi _2)\,dt.\label{18-3-51}
\end{gather}

So, \emph{modulo operator with $\bar{A}_0$-bound not exceeding
$C\bigl(\mu^{-1}h+ \vartheta (\varepsilon)+\mu^{-4}\bigr)$ we
reduced $A$ to $\bar{A}_0+W^\w$ where $W= \bar{W} +W''+W'''$ depends
on $x_2,\xi_2, r$ only\/}.

This is almost the end of the story: for
$(l,\sigma)\succeq (2,0)$ and $\mu\in [ \bar{\mu}_1,h^{-\frac{1}{3}}]$ one needs to make one step more: namely, there is operator $R^\w$ in the right-hand expression of (\ref{18-3-47}) and we should make transformation
$U '= \bigl(e^{-i\mu ^{-1}h^{-1}S^V}\bigr)^\w$ with
$S^V\in \mu^{-4} \sF^{l-2,\sigma}$ defined by
\begin{gather}
\{ \bar{a}_0, S'''\} = -(R-W^V) \label{18-3-52}\\
\shortintertext{with}
W^V (x_2,\xi_2, r )=
{\frac{1}{2\pi}}\int R( r \cos t , r \sin t; x_2,\xi _2)\,dt.\label{18-3-53}
\end{gather}

So, \emph{modulo operator with $\bar{A}_0$-bound not exceeding
$C\bigl(\mu^{-1}h+ \vartheta (\varepsilon)\bigr)$ we
reduced $A$ to $\bar{A}_0+W^\w$ where for $(l,\sigma)\succeq (2,0)$ one needs to include $W^V$ into $W$ as well\/}.

\section{Calculations}
\label{book_new-sect-18-3-2}
\subsection{Calculations. I}
\label{book_new-sect-18-3-2-1}
Now we will need to calculate $W$ more explicitly. To do this we repeat our
construction, paying attention to explicit expressions we get rather than their
smoothness. Also we will consider Weyl symbols rather than operators themselves. Let us start from the beginning of subsubsection~\ref{book_new-sect-18-3-1-3}.3 ``\nameref{book_new-sect-18-3-1-3}''. There we started from symbol 
$\bigl(\mu^2\bar{a}_0 +b\bigr)$ with
$\bigl(b=\mu^2 M'+V'\bigr)$. Then we applied $\psi _1$ transformation and arrived to the symbol
\begin{multline}
\bigl(\mu^2\bar{a}_0 +b\bigr)\circ \psi_1 =\mu^2\bar{a}_0\ +\\
\shoveright{\mu^2 \{ \bar{S} , \bar{a}_0\} +
\frac{1}{2} \mu^2\{ \bar{S} , \{ \bar{S}, \bar{a}_0\}\}\ +
\frac{1}{6} \mu^2\{ \bar{S} , \{ \bar{S} , \{ \bar{S} , \bar{a}_0\}\}\}+} \\
b \quad + \qquad
\{ \bar{S}, b\}\qquad +\qquad
{\frac{1}{2}} \{ \bar{S}, \{ \bar{S} , b\}\}\quad+\ \dots \label{18-3-54}
\end{multline}
where $\bar{S}=\bar{S}'+\bar{S}''$ and we skipped terms not exceeding
$C\bigl(\mu^{-1}h+ \vartheta (\varepsilon)\bigr)$. 

Plugging $\mu^2 \{ \bar{S}, \bar{a}_0\} = -\bar{b}+ \bar{W}_0$ where $\bar{b}$ is a $\bar{\varepsilon}$-mollification of $b$ and
\begin{equation}
\bar{W}_0= \cM _ r (\bar{b}) \Def
\frac{1}{2\pi} \int_0^{2\pi} b( r\cos t, r \sin t;x_2,\xi_2)\,dt,
\label{18-3-55}
\end{equation}
we get $\mu^2\bar{a}_0 + \bar{W}_0 +b_1$ with
\begin{multline}
b_1={\frac{1}{2}}
\{ \bar{S} , - \bar{b}+\bar{W}_0\}+(b- \bar{b})+
\{ \bar{S}, b\}+\\
\frac{1}{6} \{ \bar{S} , \{ \bar{S} , - \bar{b}+\bar{W}_0 \}\}+
 \frac{1}{2} \{ \bar{S} , \{ \bar{S}, b\}\}
\label{18-3-56}
\end{multline}
as a result of subsubsection \ref{book_new-sect-18-3-1-3}.3 ``\nameref{book_new-sect-18-3-1-3}''. 

To calculate $W$ as a result subsubsection \ref{book_new-sect-18-3-1-4}.4 ``\nameref{book_new-sect-18-3-1-4}''
(without the last step) we need to apply $\cM_ r$ to this expression:
\begin{multline}
W = \bar{W}_0+ (\cM_r b_1 ) = \cM_r( b ) +\\[3pt]
\cM_r\bigl( \{ \bar{S}, -\frac{1}{2}\bar{b}+ \frac{1}{2} \bar{W}_0+b\}\bigr) +
\frac{1}{2} \cM_r\bigl( \{\bar{S},\{\bar{S},-\frac{1}{3}\bar{b} + \frac{1}{3}\bar{W}_0 +b\}\}\bigr).
\label{18-3-57}
\end{multline}

Finally, if $l=2$, $\sigma >0$ and we want to eliminate $O(\mu^{-4})$ from an
error, we need to make a last step and define $S_1$ from the equation
\begin{equation}
\mu^2\{S_1,\bar{a}_0\}= -b_1 + (W-\bar{W}_0),
\label{18-3-58}
\end{equation}
calculate
\begin{align}b_2={\frac{1}{2}} \mu^2\{ S_1, \{ S_1, \bar{a}_0\}\}+ \{S_1, \bar{W}_0\} +
b_1 - (W-\bar{W}_0)&=\label{18-3-59}\\
\{S_1, -{\frac{1}{2}}b_1+{\frac{1}{2}}W + {\frac{1}{2}}\bar{W}_0\}+
b_1 - (W-\bar{W}_0)&\notag
\end{align}
and add to $W$ term $\cM_r(b_2)$.

Let $W_0\Def \cM_r (b)$ and $b'= \frac{1}{2}\bar{W}_0- \frac{1}{2} \bar{b} +b$.
Then
\begin{equation}
\partial _\phi \bar{S} = \frac{1}{2} \mu ^{-2}\bar{\beta}, \qquad \text{with}
\quad\beta = b-W_0, \quad \bar{\beta}= \bar{b}- \bar{W}_0.
\label{18-3-60}
\end{equation}

One can see easily that $W_0\in \sF^{l,\sigma}$ and that
$W_0 \equiv b \mod O\bigl(\vartheta ({\bar\varepsilon})+\mu ^{-2})$.

Note that we can take ${\bar{S}}=0$ as $x_1=\xi_1=0$ and in this case
\begin{equation}
|\partial_{x_2,\xi_2} \bar{S} |\le
C\mu ^{-2}\bigl(\mu ^{-1} + \bar{\varepsilon} ^{l-1}|\log \mu |^{-\sigma}\bigr).
\label{18-3-61}
\end{equation}

One can check easily that with the final error
$O\bigl(\mu^{-1}+\vartheta(\varepsilon)\bigr)$ in all the above calculations one
can replace the standard Poisson brackets $\{.,.\}$ by the ``short'' ones
$\{.,.\}'$, involving only $x_1,\xi_1$
\begin{equation}
\{f,g\}'= \partial_{\xi_1} f \cdot \partial_ {x_1} g -
\partial_{\xi_1} g \cdot \partial_ {x_1} f 
\label{18-3-62}
\end{equation}
which in polar coordinates $( r,\phi)$ at $(x_1,\xi_1)$-plane are 
\begin{equation}
\{f,g\}'= r^{-1}\bigl((\partial_\phi f)\cdot (\partial_ r g)-
(\partial_\phi g)\cdot (\partial_ r f)\bigr).
\label{18-3-63}
\end{equation}
Then integrating by parts if
$b'$ is differentiated by $\phi$ we get easily that
\begin{multline}
\cM_r \bigl(\{ \bar{S} , b'\}'\bigr)=
r ^{-1} \partial_r \cM_r (b'\partial _\phi \bar{S} )=
 \frac{1}{2} \mu^{-2} r ^{-1}\partial_r \cM_ r
\bigl(b'{\bar\beta}\bigr)\equiv \\
W_1\Def \frac{1}{4} \mu^{-2} r ^{-1}\partial_r \cM_r \bigl( (W_0-b)^2\bigr).
\label{18-3-64}
\end{multline}

Further, for $(l,\sigma)\succeq (2,0)$ let us consider the last terms in (\ref{18-3-55}) and $\cM_r (b_2)$\,\footnote{\label{foot-18-22} One can prove easily that they are not essential as $(l,\sigma)\preceq (2,0)$
and we will see that they are $\asymp \mu^{-4}$ as $(l,\sigma)\succeq (2,0)$.}.
One can see easily that only ``short'' Poisson brackets (\ref{18-3-62}) should be considered.

Then using the same arguments as above one can prove easily that
modulo
$O\bigl(\mu^{-1}h + \vartheta (\varepsilon)\bigr)$ these terms are equal
$ \frac{1}{8} \mu ^{-4} r^{-1}
\partial_ r \cM_ r \bigl( r^{-1}(b-W_0)^2 W \bigr)$ and\linebreak
$\frac{1}{4} \mu^{-2} r^{-1}\partial_ r \cM_ r
\bigl( (b_1-W_1)^2\bigr)$ respectively with $b_1$ \emph{redefined as\/}
\begin{equation}
b_1 \Def \frac{1}{2} \{S, b+W_0\}, \quad W_1 = \cM_r (b_1).
\tag*{$\textup{(\ref*{18-3-56})}^*$}\label{18-3-56-*}
\end{equation}

\subsection{Calculations. II}
\label{book_new-sect-18-3-2-2}

Now we want more explicit expressions for $W$, $W_1$ via $g^{jk}$, $V$.

First of all, note that for a constant metrics $g^{jk}$ we have
$b=V'= V\circ \Phi$, where $\Phi: T^*\bR^2\to T^*\bR^2$ is a linear map. Then
$W_0= \cM_ r V'$ and since under map $\Phi $ a circle
\begin{gather}
\cC_ r=\{(x_2,\xi_2)=\const, x_1= r \cos t,\xi_1 = r \sin t\}
\label{18-3-65}\\
\intertext{translates into ellipse with $x$-projection} 
\cE_{r}(y)=
\{\sum_{j,k}g_{jk}(x_j-y_j)(x_k-y_k)= r ^2\}
\label{18-3-66}
\end{gather}
where $y=\Phi (\bar{x}_2,\bar{\xi}_2)$ is its center and $g_{jk}$ are elements of the inverse matrix $(g^{jk})^{-1}$, we get that 
\begin{equation}
W_0=\cM_{\cE_r }(\Phi (\bar{x}_2,\bar{\xi}_2))
\label{18-3-67}
\end{equation}
which is a corresponding average of $V$ along $\cE _{ r}(y)$.

Similarly, we see that 
\begin{multline}
W_1\equiv - \frac{1}{4}\mu^{-2}|\nabla V'|^2\equiv\\
-\frac{1}{4}\mu^{-2}\sum g^{jk}(\partial _{x_j}V)(\partial _{x_k}V)\circ \Phi 
\mod O(\mu ^{-4}).
\label{18-3-68}
\end{multline}

\begin{Problem}\label{problem-18-3-4}
In the general case, the deviation of $\cC_ r \circ \Phi $ from the ellipse $\cE_r$ (with $g_{jk}$ calculated in the center) is $O(\mu ^{-2})$ and thus
\begin{equation}
W\equiv \cM_{\cE _ r }(\Phi (\bar{x}_2,\bar{\xi}_2)) \mod O(\mu^{-2}).
\label{18-3-69}
\end{equation}
This does not require $l>1$. It would be nice to get more
precise answer in the general case as well.
\end{Problem}

\section{Strong magnetic field $\mu \ge C_0^{-1}(h|\log h|)^{-1}$}
\label{book_new-sect-18-3-3}

In this case we cannot claim anymore that $|p_j|\le C\mu ^{-1}$ in the
microlocal sense even if $\mu \le ch^{-1}$ because
$f\bigl(\mu p_1(x,\xi), \mu p_2x,\xi), x\bigr)$ is
no more quantizable symbol. However, now $\mu ^{-1}\le \varepsilon$ and
$f\bigl(\varepsilon ^{-1}p_1(x,\xi),\varepsilon ^{-1}p_2(x,\xi), x\bigr)$ is
quantizable symbol and we can claim that $|p_j|\le C\varepsilon$ in the
microlocal sense in the zone of $\tau \le C_0(1+ \mu h)$ due to the standard ellipticity arguments. Basically, this makes our construction much simpler; however, we need to remember that $h^2$ is no more non-essential term.

Construction of subsubsections~\ref{book_new-sect-18-3-1-1}.1, \ref{book_new-sect-18-3-1-2}.2 remains the same but we need to take in account correction (\ref{18-1-8})
appearing when multiplying $A$ by $F^{-\frac{1}{2}}$ from the left and right
and moving both copies of $F^{-\frac{1}{2}}$ inside to $g^{jk}$.

Further, since the principal part of operator $A$ is of magnitude $\mu h$ rather than of magnitude $1$ as it was before, the  error due to approximation of the main part of operator is
$O\bigl(\mu ^2 \varepsilon ^{\bar{l}} |\log \mu|^{-\bar{\sigma}}\bigr)$ and it is non-essential (i.e. $O(\mu ^{-1}h+1)$) for $(\bar{l},\bar{\sigma})=(4,2)$ only.

In the very input of subsubsection~\ref{book_new-sect-18-3-1-2}.2 we can take a Taylor decomposition with respect to $(x_1,\xi_1)$ (as we did in the smooth case) and then modulo ``non-essential symbol'' (which belongs to
$\bigl(\mu^{-1}h + \vartheta (\varepsilon)\bigr)\sF^{0,0}$) we get
\begin{multline}
\mu ^2 (x_1^2+\xi_1^2) + V(x_2,\xi_2) + \\
\mu^2 \sum_{j=0}^3 \beta _j(x_2,\xi _2) x_1^j \xi_1^{3-j} +
\sum _{j=0}^1 \alpha _j (x_2,\xi_2) x_1^j \xi_1^{1-j}+\\
\mu^2 \sum_{j=0}^4 \beta '_j(x_2,\xi _2) x_1^j \xi_1^{4-j} +
\sum _{j=0}^2 \alpha '_j (x_2,\xi_2) x_1^j \xi_1^{2-j}
\label{18-3-70}
\end{multline}
with $\beta _j \in \sF^{1,1}$, $\beta _j \in \sF^{0,1}$,
$\alpha _j \in \sF^{l-1,\sigma}$, $\beta _j \in \sF^{l-2,\sigma}$
and we need the last term only if $l\ge 2$.

However, we need to remember the terms in decomposition of the main part;
we can eliminate the third order term but we will reduce the fourth-order term.
More precisely, to get rid of the third and the forth terms in (\ref{18-3-70}) we use again the construction of subsubsection~\ref{book_new-sect-18-3-1-3}.3 but we look explicitly for
\begin{equation}
S'=\sum _{j=0}^1 \hat{\alpha}_j (x_2,\xi_2) x_1^j \xi_1^{1-j},\qquad
S''= \mu^2 \sum_{j=0}^3 \hat{\beta}_j(x_2,\xi _2) x_1^j \xi_1^{3-j}
\label{18-3-71}
\end{equation}
with $\hat{\alpha}_j, \hat{\beta}_j$ of the same regularity 
as $\alpha _j,\beta _j$.

Now we have expression (\ref{18-3-70}) without the third and the fourth terms; the fifth and the sixth terms are modified but retain $\sF^{l-2,\sigma}$-regularity provided $(l,\sigma )\succ (2,0)$; for 
$(l,\sigma )\preceq (2,0)$ we get non-essential terms. For $\mu \ge c^{-1}h^{-1}$ these terms are non-essential as well. So, \emph{in these two cases our construction is complete\/}.

To complete reduction for $\mu \in [(h|\log h|)^{-1}, ch^{-1}]$, 
$(l,\sigma )\succ (2,0)$ we apply again ``poor man Fourier integral operator'', this time with
\begin{equation}
S'=\sum _{j=0}^2 \hat{\alpha}'_j (x_2,\xi_2) x_1^j \xi_1^{2-j},\qquad
S''= \mu^2 \sum_{j=0}^4 \hat{\beta}'_j(x_2,\xi _2) x_1^j \xi_1^{4-j}
\label{18-3-72}
\end{equation}
with $\hat{\alpha}'_j, \hat{\beta}'_j\in \sF^{l-2,\sigma}$. However,
we cannot eliminate these terms completely, we can only to reduce them to
\begin{phantomequation}\label{18-3-73}\end{phantomequation}
\begin{equation}
\bar{\alpha} (x_2,\xi_2) \bigl(x_1^2+\xi_1^2\bigr), \qquad
\bar{\beta} (x_2,\xi_2) \bigl(x_1^2+\xi_1^2\bigr) ^2
\tag*{$\textup{(\ref*{18-3-73})}_{1,2}$}
\end{equation}
respectively with $\bar{\alpha},\bar{\beta}\in \sF^{l-2,\sigma}$.

One can calculate easily that the final expression is resembling one from subsubsections~\ref{book_new-sect-18-3-2-1}.1,~\ref{book_new-sect-18-3-2-2}.2 because in our situation $\cM_\cC f=f + {\frac{1}{2}} \mu ^{-2}\Delta _{x_1,\xi_1}f$ modulo
non-essential term and one can recalculate this way $W_1$ too. However,
construction here reminds one of the smooth theory.

So, we have proven

\begin{proposition} \label{prop-18-3-5} 
Let $d=2$, $\mu \ge (h|\log h|)^{-\frac{1}{3}}$ and $(\bar{l},\bar{\sigma})\succeq (2,0)$, $(l,\sigma)\succeq (1,2)$.

Then modulo operator with $\bar{A}_0$-bound not exceeding
$C\bigl(\mu^{-1}h+ \vartheta (\varepsilon)\bigr)$ with $\varepsilon = C(\mu^{-1}h|\log h|)^{\frac{1}{2}}$ one can reduce operator $A$ to operator $\bar{A}_0+W^\w$ where $W= W_0+W_1+W_2$ depends on $x_2,\xi_2$ and 
$\mu ^{-1}A_0^{\frac{1}{2}}$ (or $\mu ^{-1}\bar{A}_0^{\frac{1}{2}}$) only.
\end{proposition}

\chapter{$d=2$: Tauberian theory}
\label{book_new-sect-18-4}

\section{Decomposition}
\label{book_new-sect-18-4-1}

We again consider more difficult case $\mu \le \epsilon (h|\log h|)^{-1}$
first; the opposite case can be considered in the same way with rather obvious
modifications and simplifications.

There is not much what is left: removing non-essential operator and
compensating this by adding $\mp C_0\vartheta (\varepsilon )$ to the lower/upper
approximation (which in turn will not violate our estimates; see below), we get
a $1$-dimensional operator-valued $\mu ^{-1}h$-pseudo-differential operator
\begin{equation}
\bar{A}_0+ W(x_2,\mu^{-1}hD_2, \mu ^{-1}\bar{A}_0)
\label{18-4-1}
\end{equation}
with an operator-valued symbol, acting in auxiliary space $\sL^2(\bR)$. 

Taking as in Section~\ref{book_new-sect-13-4} decomposition into 
\begin{equation}
\Upsilon _n (x_1)\Def (\mu ^{-1}h)^{- \frac{1}{4}}
\upsilon _n\bigl( (\mu ^{-1}h)^{- \frac{1}{2}} x_1\bigr)
\label{18-4-2}
\end{equation}
with Hermitian functions $\upsilon _n$ and coefficients which are functions of
$x_2$ we get a family of $1$-dimensional $\mu^{-1}h$-pseudo-differential operators
\begin{equation}
\cA_n\Def r_n^2 + W(x_2,\mu^{-1}hD_2,\mu^{-1}r_n),\qquad
r_n \Def \bigl((2n+1)\mu h\bigr)^{\frac{1}{2}}
\label{18-4-3}
\end{equation}
as $\mu h n\le C_0$ .

Even if we cannot apply directly spectral projectors associated with
$\cA_n$ we can consider propagator $e^{ih^{-1}t\tilde{A}}$ and its Schwartz kernel $u(x,y,t)$ and prove easily that for 
$t, |t|\le T_1\Def \epsilon \mu $
\begin{equation}
U(x,y,t) \equiv \cT
\biggl(\sum _{n\ge 0} U_n(x_2,y_2,t) \Upsilon _n (x_1)\Upsilon _n (y_1)\biggr)
\cT^{-1}
\label{18-4-4}
\end{equation}
modulo Schwartz kernel of an operator which becomes negligible multiplied by $\bigl(\psi (x) q (\mu p_1,\mu p_2)\bigr)^\w$ where $\psi $, $q$ are smooth and
supported in $B(0,{\frac{1}{2}})$ and $B(0,C_0)$ respectively; recall that we call operator \emph{negligible\/} if its norm does not exceed $C(\mu ^{-1}h)^s$ with large enough exponent $s$.

Here $U$ is the Schwartz kernel of the propagator $e^{ih^{-1}t\tilde{A}}$ of
the transformed back reduced operator\footnote{\label{foot-18-23} Apart of just reduction we removed some ``non-essential'' terms.}\footnote{\label{foot-18-24} Actually we have two of them: $\tilde{A}^\pm$ such that
$\tilde{A}^-\le A\le \tilde{A}^+$ and $A$-bounds of operators
$A-\tilde{A}^\pm$ do not exceed
$C\bigl(\mu ^{-1}h + \vartheta (\varepsilon)\bigr)$ (so they are operators
described in footnote \footref{foot-18-16}).}
\begin{equation}
\tilde{A}\Def \cT
\Bigl(\bar{A}_0+ W(x_2,\mu^{-1}hD_2, \mu ^{-1}\bar{A}_0)\Bigr)\cT^{-1}
\label{18-4-5}
\end{equation}
which is $1$-dimensional $\mu^{-1}h$-pseudo-differential operator with operator-valued symbol, $\cT$ is a full transformation constructed in subsection~\ref{book_new-sect-18-3-1}, and 
\begin{claim}\label{18-4-6}
$U_n$ are Schwartz kernels of the propagators $e^{ih^{-1}t\cA _n}$ for $\cA_n$. 
\end{claim}
Further, by definition $W=W_0+W_1$ because $W_2$ is not essential in our sense.

\section{Estimate under assumption (\ref{18-2-15})}
\label{book_new-sect-18-4-2}

So, for $T=T_1$, $\tau \le c$
\begin{multline}
F_{t\to h^{-1}\tau }\Bigl(\chi_T(t) \Gamma \bigl(\psi (x)
u(x,y,t)\bigr)\Bigr)\equiv \\
F_{t\to h^{-1}\tau }\biggl(\chi_T(t)
\sum _{n\ge 0} \Gamma \Bigl(\tilde{\psi}
\bigl(U_n(x_2,y_2,t) \Upsilon _n(x_1)\Upsilon_n(y_1) \Bigr)\biggr)
\label{18-4-7}
\end{multline}
modulo $O(\mu ^{-s})$ where as before $\Gamma $ is a trace for Schwartz kernel
(restriction to the diagonal and integration) and 
\begin{equation}
\tilde{\psi}=\cT^{-1}\psi \cT
\label{18-4-8}
\end{equation}
is $\mu ^{-1}h$-pseudo-differential operator rather than a function; its symbol is not very regular (see below).

In this decomposition we can consider only $n\le C_1(\mu h)^{-1}$ because for any other $n$ interval $(-\infty , C_0)$ is a classically forbidden zone for
operator $\cA_n$ and one can prove easily that their contributions do not
exceed $C\mu ^{-s}(\mu h n +1)^{-s} $ with the arbitrarily large $s$ and thus
their total contribution is negligible.

Note that all operators $\cA_n$ are non-degenerate at level $0$ in the sense that
\begin{equation}
|\cA_n | \le \epsilon_0\implies |\nabla _{x_2,\xi_2}\cA_n| \asymp 1
\label{18-4-9}
\end{equation}
if and only if the original operator satisfied assumption (\ref{18-2-15}) because $|\nabla (W-V/F)|\le C\mu ^{-1}$. 

Let us assume first that this non-degeneracy assumption holds. Then in the standard way we can prove that if
$Q_x=q(x_2,\mu^{-1}hD_{x_2})$ and $Q'_y=q'(y_2,\mu^{-1}hD_{y_2})$ are
$\mu^{-1}h$-pseudo-differential operators with symbols supported in $\varepsilon$-vicinity of $(\bar{x}_2,\bar{\xi}_2)$ in which condition (\ref{18-4-9}) is fulfilled, then 
\begin{equation}
F_{t\to h^{-1}\tau }\chi_T(t) \bigl( Q_x U_n(x_2,y_2,t)\,^t\!Q'_y\bigr)
\label{18-4-10}
\end{equation}
is negligible provided $T\le T_1$ and
$\mu ^{-1}T\times \varepsilon \ge C\mu ^{-1}h|\log h|$ or equivalently
\begin{equation}
T_* \Def C\varepsilon ^{-1} h|\log h|\le T \le T^* =\epsilon \mu ;
\label{18-4-11}
\end{equation}
recall that $\chi $ is an admissible function supported in
$[-1,-\frac{1}{2}]\cup [\frac{1}{2},1]$. Then
\begin{equation}
F_{t\to h^{-1}\tau }\chi_T(t) \Gamma \bigl(\psi (x) U(x,y,t)\bigr)
\label{18-4-12}
\end{equation}
is negligible as well.

This statement does not look strong enough because $T_*$ seems to be much larger than we need due to the factor $\varepsilon^{-1}$: the direct assault shows that $|\Gamma \psi K U|$ does not exceed $Ch^{-2}$ 
(where $K= k\bigl(x,\mu p (x,\xi)\bigr)^\w)$ is $\mu^{-1}h$-pseudo-differential operator cut-off; then
\begin{equation}
|F_{t\to h^{-1}\tau}\bigl(\bar{\chi}_T(t)\Gamma (\psi KU)\bigr)|\le Ch^{-2}T_*
\label{18-4-13}
\end{equation}
for $\bar{\chi}$ supported in $[-1,1]$ and equal $1$ at $[-\frac{1}{2},\frac{1}{2}]$ and $T\in [T_*,T^*]$; here we can take $K=I$ because $\tau \le c$ and we can apply ellipticity arguments in zone 
$\{|p| \ge C'\mu ^{-1}\}$; then (\ref{18-4-13}) and
the standard Tauberian arguments imply Tauberian spectral remainder estimate
\begin{multline}
\R^\T\Def |\Gamma (\psi Q_x e)(\tau ) -
h^{-1}\int _{-\infty}^\tau \Bigl(F_{t\to h^{-1}\tau'}
\bigl( \bar{\chi}_T(t)\Gamma (\psi Q U)\bigr)\Bigr)\,d\tau' | \le\\
Ch^{-2}{\frac{T_*} {T^*}} \asymp Ch^{-1} \mu ^{-1}\varepsilon^{-1}|\log h|
\label{18-4-14}
\end{multline}
 and it is close to
what we are looking for (but still is not that good) only for $\varepsilon=1$. Sure by rescaling arguments we can make $T_*= h\varepsilon^{-1}$ in the right-hand expressions of the estimates only and thus get rid off logarithmic factor (and we will do it) but it alone does not solve our problem.

The source of this trouble is that \emph{we have $\varepsilon $-scale in both $x_2$ and $\xi_2$\/},
while in the differential case treated in Section~\ref{book_new-sect-4-5} we had scale $1$ in $\xi_2$.

However, there is a way to improve our analysis: what we need is an estimate
\begin{equation}
|F_{t\to h^{-1}\tau} \bigl( \bar{\chi}_T(t)\Gamma (\psi U)\bigr)|\le
Ch^{-1}+C \mu h^{-2} \vartheta (\varepsilon)
\label{18-4-15}
\end{equation}
which we are going to prove; then dividing the right-hand expression of (\ref{18-4-15}) by $T^*$ we get the required remainder estimate
$O\bigl(\mu ^{-1}h^{-1}+h^{-2}\vartheta (\varepsilon)\bigr)$.

\medskip\noindent
(i) First of all, we want to deal with $U_n$ directly. Note that
\begin{multline*}
\tilde{\psi} =
\psi _0 (x_2,\mu ^{-1}hD_2)+\\
\psi _1(x_2,\mu ^{-1}hD_2)x_1+
\psi _2(x_2,\mu ^{-1}hD_2)\mu^{-1}hD_1+
\tilde{\psi}'(x,\mu^{-1}hD)
\end{multline*}
with symbol of $\tilde{\psi}'$ belonging to
$\mu ^{-l} |\log h |^{-\sigma}\sF^{0,0}$; for $l=2,\sigma>0$ one needs to
include also the second order terms
\begin{multline*}
\psi_{11}(x_2,\mu^{-1}hD_2) \,x_1^2+
\psi_{12}(x_2,\mu^{-1}hD_2)\,(x_1\mu^{-1}hD_1 +\mu^{-1}hD_1x_1)+\\
\psi_{11}(x_2,\mu^{-1}hD_2)\,(\mu ^{-1}hD_1)^2.
\end{multline*}

Let us consider
$|F_{t\mapsto h^{-1}\tau} \bigl(\bar{\chi}_T(t)\tilde{\psi}' U\bigr)|$ first;
due to analysis above it does not exceed
$Ch^{-2}T_* \times \mu ^{-l}|\log h|^{-\sigma}$; one can see easily that it does
not exceed the right-hand expression of (\ref{18-4-15}).

Further, note that plugging terms linear with respect to $(x_1,\mu^{-1}hD_1)$ kills the trace while plugging quadratic terms is equivalent to plugging
$(n+{\frac{1}{2}})\mu^{-1}h\psi ' (x_2,\mu^{-1}hD_2)$ with
$\psi '=\psi_{11}+\psi_{22}$.

Now we have pure $\mu^{-1}h$-pseudo-differential operators, acting on $x_2$ alone; we consider $u_n$ and replace $\Gamma$ by $\Gamma '$ acting with respect to $x_2$ alone.

\medskip\noindent
(ii) Second, we need to improve $T_*$ slightly, removing $|\log h|$ factor and to get better estimates for Fourier transform. Let us consider
$\varepsilon$-admissible partition in $(x_2,\xi_2)$:
$1 = \sum_\nu q_{(\iota)} $ and consider element $q_{(\iota )}$ of this partition.

One can easily prove by the standard rescaling technique the following estimate
\begin{equation}
|\Gamma '\bigl(\chi_T(t) q_{(\iota)} u_n\bigr) |\le
C\varepsilon ^2 \mu h^{-1}\bigl( \frac{T} {\bar{T}_*}\bigr)^s
\label{18-4-16}
\end{equation}
for arbitrarily large $s$ and $T\in [\bar{T}_*,T^*]$ where
\begin{equation}
\bar{T}_* \Def C_0\varepsilon ^{-1} h = T_*/|\log h|.
\label{18-4-17}
\end{equation}

Then
\begin{gather}
|F_{t\to h^{-1}\tau} \bigl(\chi _T(t) \Gamma' (q_{(\iota)}U_n) \bigr)|\le
C\varepsilon ^2 \mu h^{-1}\bigl(\frac{T}{\bar{T}_*}\bigr)^s T 
\label{18-4-18}\\
\shortintertext{and thus}
\quad |F_{t\to h^{-1}\tau } \bigl( \bar{\chi }_T(t) \Gamma' (q_{(\iota)}U_n)
\bigr)|\le
C\varepsilon ^2 \mu h^{-1}\bar{T}_*
\label{18-4-19}
\end{gather}
for any $T\in [\bar{T}_*,T^*]$.

Further, note that due to the standard ellipticity arguments
\begin{equation*}
F_{t\to h^{-1}\tau } \bigl(\bar{\chi}_T(t) \Gamma' (q_{(\iota)}U_n)\bigr)
\end{equation*}
is negligible for $|\tau|\le Ch|\log h| T^{-1}$ unless
\begin{equation}
|\cA_n(x_2,\xi _2)|\le C\bigl(\varepsilon + Ch|\log h|T^{-1} \bigr) \qquad \text{on\ \ } \supp q_{(\iota)},
\label{18-4-20}
\end{equation}
and that
\begin{claim}\label{18-4-21}
Under assumption (\ref{18-2-15}) the set
\begin{equation}
\Omega _{\tau,n} \Def \{(x_2,\xi_2): |\cA_n(x_2,\xi _2)-\tau |\le C\varepsilon\}
\label{18-4-22}
\end{equation}
has measure $O(\varepsilon)$ for each $\tau $, $n$. 
\end{claim}
Therefore as $T\ge T_*$ (and therefore $Ch|\log h|T^{-1}\lesssim \varepsilon$ we arrive to the estimate
\begin{equation}
|F_{t\to \mu ^{-1}h} \bar{\chi}_T(t) \Gamma' (QU_n) |\le
C \mu h^{-1} \bar{T}_* \times \varepsilon = C\mu
\label{18-4-23}
\end{equation}
which implies (\ref{18-4-15}) as we need to multiply by the number of admissible indices (which is $\asymp (\mu h)^{-1}$) resulting in $Ch^{-1}$ and take in account contribution of neglected terms. Then in virtue of the standard Tauberian arguments (with $T^*\asymp \mu$) we arrive to the estimate 
\begin{equation}
\R^\T \le C\mu^{-1}h^{-1} + C h^{-2}\vartheta (\varepsilon)
\label{18-4-24}
\end{equation}
with $\R^\T$ defined by (\ref{18-4-14}) $\bar{\chi}=1$ on 
$[-\frac{1}{2} , \frac{1}{2} ]$, $T\in [T_0,T_1]$ which is exactly estimate we want. 

\section{Estimate under assumption (\ref{18-2-16})}
\label{book_new-sect-18-4-3}

Now let $(l,\sigma)\succeq (2,0)$ and at some point 
$\bar{z}\Def (\bar{x}_2,\bar{\xi}_2)$ the following generalization of (\ref{18-4-9}) holds:
\begin{gather}
|\nabla \cA_n |\le \epsilon_0\nu \implies 
|\nabla \cA_n (x_2,\xi_2)|\asymp \nu 
\tag*{$\textup{(\ref*{18-4-9})}'$}\label{18-4-9-'}\\
\shortintertext{where}
\nu \ge C_0\mu^{-1}
\label{18-4-25}
\end{gather}
which in turn is larger than $C_0\varepsilon$ as 
$\mu \le \epsilon (h|\log h|)^{-1}$. Then \ref{18-4-9-'} holds in $B(\bar{z},\epsilon_0\nu)$.

Then we can take 
\begin{equation}
T_*\Def C \varepsilon^{-1}\nu^{-1}h|\log h| \le T^*=\epsilon \mu
\tag*{$\textup{(\ref{18-4-17})}'$}\label{18-4-17-'}
\end{equation}
and (\ref{18-4-18}), (\ref{18-4-19}) hold. On the other hand, one needs to replace (\ref{18-4-20}) by 
\begin{equation}
|\cA_n(x_2,\xi _2)|\le C\bigl(\varepsilon \nu + Ch|\log h|T^{-1} \bigr) \qquad \text{on\ \ } \supp q_{(\iota)},
\tag*{$\textup{(\ref*{18-4-20})}'$}\label{18-4-20-'}
\end{equation}
and instead of (\ref{18-4-21}) observe that 
\begin{claim}\label{18-4-26} 
Under assumption \ref{18-4-9-'} 
\begin{equation}
\Omega _{\tau,n,\nu } \Def \{(x_2,\xi_2) \in B(\bar{z},\epsilon_0\nu): |\cA_n(x_2,\xi _2)-\tau |\le C\varepsilon\nu \}
\tag*{$\textup{(\ref*{18-4-22})}'$}\label{18-4-22-'}
\end{equation}
has measure $O(\varepsilon \nu )$ for each $\tau $, $n$. 
\end{claim}

Then, as $Q$ is supported in $B(\bar{z},\epsilon_0\nu)$, (\ref{18-4-23}) is replaced by 
\begin{equation}
|F_{t\to \mu ^{-1}h} \bar{\chi}_T(t) \Gamma' (Qu_n) |\le
C \mu h^{-1} \bar{T}_* \times \nu \varepsilon = C\mu 
\tag*{$\textup{(\ref*{18-4-23})}'$}\label{18-4-23-'}
\end{equation}
and summation over $n$ returns 
$C\mu (\nu^2/(\mu h) +1) $ while summation over $z$ returns
\begin{equation}
C\int \mu (\nu(z)^2/(\mu h) +1) \nu^{-2}(z)\,dz \asymp Ch^{-1}+ C\mu |\log \mu|
\label{18-4-27}
\end{equation}
due to assumption (\ref{18-2-16}). As $\mu \lesssim (h|\log h|)^{-1}$ the second term in the right-hand expression is less than the first one. 

Therefore if $Q$ is supported in the ball of radius $\asymp 1$ and is negligible on $\cX_\inter\Def \{z: |\nabla W|\le C_0\mu^{-1}\}$ inequality (\ref{18-4-23}) remains true and the contribution of $Q$ to Tauberian remainder does not exceed $C\mu^{-1}h^{-1}$. 

Meanwhile contribution of $\cX_\inter$ (intersected with the ball of radius 
$\asymp 1$) to the Tauberian remainder does not exceed 
$C\mu h^{-1}\times \mes \cX_\inter=O(\mu^{-1}h^{-1})$ in virtue of (\ref{18-2-16}) again. Thus we arrive to Tauberian estimate (\ref{18-4-24}) again. Therefore we arrive to

\begin{proposition}\label{prop-18-4-1}
Let $d=2$, $(\bar{l},\bar{\sigma})\succeq (2,0)$, 
$h^{-\frac{1}{3}}\lesssim (h|\log h|)^{-\frac{1}{3}}\le \mu \le \epsilon (h|\log h|)^{-1}$ and \underline{either} 

\medskip\noindent
(i) $(l,\sigma)\succeq (1,2)$ and assumption \textup{(\ref{18-2-15})} be fulfilled \underline{or} 

\medskip\noindent
(ii) $(l,\sigma)\succeq (2,0)$ and assumption \textup{(\ref{18-2-16})} be fulfilled.

\medskip
Then Tauberian estimate \textup{(\ref{18-4-24})} holds where we need to calculate Tauberian expression 
\begin{equation}
h^{-2}\cN^\T\Def h^{-1}\int _{-\infty}^\tau \Bigl(F_{t\to h^{-1}\tau'}
\bigl(\bar{\chi}_T(t)\Gamma (\psi u)\bigr)\Bigr)\,d\tau'
\label{18-4-28}
\end{equation}
with $T\ge T_*= C\varepsilon^{-1} h|\log h|$ in (i) and take partition and take $T\ge T_{*\iota}= C\varepsilon^{-1}\nu_\iota^{-1} h|\log h|$. 
\end{proposition}

\section{Case $\mu \ge \epsilon (h|\log h|)^{-1}$}
\label{book_new-sect-18-4-4}

As $\mu \ge \epsilon (h|\log h|)^{-1}$ we need to remember that now $\varepsilon\ge \mu^{-1}$; therefore estimates 
\begin{gather}
|x_1|\le C_0\mu^{-1}\qquad |\mu^{-1}hD_1|\le C_0\mu^{-1}
\label{18-4-29}\\
\intertext{should be replaced by} 
|x_1|\le C_0\varepsilon,\qquad |\mu^{-1}hD_1|\le C_0\varepsilon 
\label{18-4-30}
\end{gather}
in the ``domain of interest''; however as $\mu \lesssim h^{-1}$ estimate (\ref{18-4-28}) still holds in the operator sense.

Further, as $\mu \ge c h^{-1}$ we need to consider Schr\"odinger-Pauli operator 
\begin{equation}
A=\sum_{j,k} P_jg^{jk}P_k+V-\fz \mu h F,\qquad \text{with\ \ } P_j=hD_j-\mu V_j
\label{18-4-31}
\end{equation}
with $\fz \in \bR$ and to avoid hitting the spectral gap we need to consider only $\fz=1,3,5,\ldots$. Then (\ref{18-4-30}) holds in microlocal sense and also in operator sense (but in the latter case one can take 
$\varepsilon= C(\mu^{-1} h)^{\frac{1}{2}}$).

Furthermore, one need to remember that $n\le C_0/\mu h$ should be replaced by $n\le C_0$ as $\mu \ge c h^{-1}$. This takes care of the case of assumption (\ref{18-2-15}) and we arrive to the estimate (\ref{18-4-24}) as 
$\mu \le h^{-1}$ and estimate 
\begin{equation}
\R^\T \le C + C \mu h \vartheta (\varepsilon)
\label{18-4-32}
\end{equation}
as $\mu \ge h^{-1}$.

Under assumption (\ref{18-2-16}) additional analysis is required as for
$\{\nu \le \varepsilon\}$ uncertainty principle breaks but contribution of such zone to the remainder is $C h^{-2} (\mu h +1) \varepsilon^2$ (i.e. estimates (\ref{18-4-24}) and (\ref{18-4-32}) acquire additional terms 
$C\mu^{-1}h^{-1}|\log \mu|$ and $C|\log \mu|$ respectively. 

Thus we arrive to 

\begin{proposition}\label{prop-18-4-2}
Let $d=2$ and $(\bar{l},\bar{\sigma})\succeq (2,0)$. 

\medskip\noindent
(i) Let $(l,\sigma)\succeq (1,2)$ and assumption \textup{(\ref{18-2-15})} be fulfilled. Then for $(h|\log h|)^{-1} \le \mu \le h^{-1}$ Tauberian estimate \textup{(\ref{18-4-24})} holds for Schr\"odinger operator and for 
$\mu \ge h^{-1}$ Tauberian estimate \textup{(\ref{18-4-32})} holds for Schr\"odinger-Pauli operator where we need to calculate Tauberian expression \textup{(\ref{18-4-32})} with $T\ge T_*= C\varepsilon^{-1} h|\log h|$;

\medskip\noindent
(ii) Let $(l,\sigma)\succeq (2,0)$ and assumption \textup{(\ref{18-2-16})} be fulfilled. Then for $(h|\log h|)^{-1} \le \mu \le h^{-1}$ Tauberian estimate 
\begin{equation}
\R^\T \le C\mu^{-1}h^{-1}+C |\log \mu| + C h^{-2}\vartheta (\varepsilon)
\label{18-4-33}
\end{equation}
holds for Schr\"odinger operator and for $\mu \ge h^{-1}$ Tauberian estimate
\begin{equation}
\R^\T \le C|\log \mu| + C \mu h^{-1}\vartheta (\varepsilon)
\label{18-4-34}
\end{equation}
holds for Schr\"odinger-Pauli operator where we need to take partition and to calculate Tauberian expression \textup{(\ref{18-4-32})} with 
$T\ge T_{*\iota}= C\varepsilon^{-1}\nu_\iota^{-1} h|\log h|$. 
\end{proposition}

Repeating arguments of the smooth case (in section~\ref{book_new-sect-13-5}) one can prove easily the following

\begin{remark}\label{rem-18-4-3}
(i) As $\mu \gtrsim h^{-1}$ as usual conditions (\ref{18-2-15}), (\ref{18-2-16}) could be weaken to
\begin{equation}
|(2n+1-\fz )\mu h +V|+ |\nabla V/F|\ge \epsilon_0 \qquad \forall n\in \bZ^+
\tag*{$\textup{(\ref*{18-2-15})}^*$}\label{18-2-15-*}
\end{equation}
and 
\begin{multline}
|(2n+1-\fz )\mu h +V|+ |\nabla V/F|+|\det \Hess V/F| \ge \epsilon_0\\
\forall n\in \bZ^+
\tag*{$\textup{(\ref*{18-2-16})}^*$}\label{18-2-16-*}
\end{multline} 
respectively.

\medskip\noindent
(ii) As $\mu \lesssim h^{-1}$ under condition 
\begin{equation}
 |\nabla V/F|\le \epsilon_0 \implies \det \Hess V/F\ge \epsilon_0 \\
\tag*{$\textup{(\ref*{18-2-16})}^{+}$}\label{18-2-16-+}
\end{equation} 
estimate (\ref{18-4-24}) holds and as $\mu \gtrsim h^{-1}$ under condition
\begin{multline}
|(2n+1-\fz )\mu h +V|+ |\nabla V/F|\le \epsilon_0 \implies \det \Hess V/F\ge \epsilon_0 \\
\forall n\in \bZ^+
\tag*{$\textup{(\ref*{18-2-16})}^{*+}$}\label{18-2-16-*+}
\end{multline} 
estimate (\ref{18-4-32}) holds.

\medskip\noindent
(iii) As $\mu \gtrsim h^{-1}$ in the case of the spectral gap 
\begin{gather}
|(2n+1-\fz )\mu h +V|\ge \epsilon_0 \qquad \forall n\in \bZ^+
\label{18-4-35}
\end{gather}
$\R^\T = O(\mu^{-s})$ as in the smooth case.
\end{remark}

Here in (ii) we use the fact that under condition \ref{18-2-16-+} or \ref{18-2-16-*+} (in the corresponding settings) $\cA_n$ is is \underline{either} elliptic unless $\nu \asymp \bar{\nu}$ for some 
$\bar{\nu}\ge C_0\max \bigl(\mu^{-1},(\mu^{-1}h)^{\frac{1}{2}}\bigr)$ \underline{or} unless $\nu \le \bar{\nu}= C_0\max \bigl(\mu^{-1},(\mu^{-1}h)^{\frac{1}{2}}\bigr)$. 

As $\bar{\nu}\in[ (\mu^{-1}h)^{\frac{1}{2}}, (\mu^{-1}h|\log h|)^{\frac{1}{2}}]$ we use a rescaling technique.

\chapter{$d=2$: Calculations and  main theorems}
\label{book_new-sect-18-5}

Now our purpose is to prove main results as $\mu \ge\bar{\mu}_1$, i.e. theorems \ref{thm-18-5-4} and \ref{thm-18-5-7}. In calculations we can assume that \ref{18-1-6-*} holds, i.e. $F=1$.

\section{Calculations}
\label{book_new-sect-18-5-1}

\subsection{Step I under assumption (\ref{18-2-15})}
\label{book_new-sect-18-5-1-1}

To calculate Tauberian expression (\ref{18-4-28}) we rewrite it modulo negligible terms as the sum with respect to $n,\iota$ of
\begin{multline}
\cS_{n(\iota)}\Def \\
h^{-1}\int _{-\infty}^\tau \biggl(F_{t\to h^{-1}\tau'}\Bigl(
\bar{\chi}_T(t)\Gamma \bigl(\tilde{\psi}(x,\mu^{-1}hD) q_{(\iota)}^\w
\bigl(u_n \Upsilon _n(x_1)\Upsilon _n(x_1)\bigr)\bigr)\Bigr)\biggr)\,d\tau'
\label{18-5-1}
\end{multline}
with $q_{(\iota)}=q_{(\iota)}^\w (x_2,\mu^{-1}hD_2)$.

One can see easily that on the \emph{elliptic elements\/} (with failed (\ref{18-4-20}) or \ref{18-4-20-'})
\begin{equation}
\cS_{n(\iota)} \equiv \frac{1}{2\pi} \mu h^{-1}\iint
\uptheta \bigl(\tau -\cA_n (x_2,\xi_2)\bigr) \tilde{\psi}(x,\xi)
q_{(\iota)}(x_2,\xi_2)\,dx d\xi
\label{18-5-2}
\end{equation}
modulo negligible expression; here
$\uptheta \bigl(\tau -\cA_n (x_2,\xi_2) \bigr)=\const$ on $\supp q_{(\iota)}$
and we need to calculate the difference between these two expressions on the
\emph{non-elliptic elements\/} (i.e. satisfying (\ref{18-4-20}) or \ref{18-4-20-'}). 

The total contribution of such terms
is not very large: as $\mu h\lesssim 1$ and $\nu=1$ it does not exceed
$C h^{-1} \times\mu \times \varepsilon \times 1/(\mu h)$
where the first factor is just present in (\ref{18-5-1}), the second one is an estimate (\ref{18-4-23}) of the Fourier transform
($F_{t\to \mu ^{-1}h} \bar{\chi }_T(t) \Gamma' (Qu_n)$), the third is here
because this Fourier transform is negligible for 
$|\tau ' -\tau |\ge C\varepsilon $ and the last one is the number of $n$ for which (\ref{18-4-20}) holds on $\supp q_{(\iota)}$.
Still resulting expression $C\varepsilon h^{-2}$ is larger than what we need.

We apply the method of successive approximations to calculate these terms; as
unperturbed operator we take
\begin{gather}
\cB _n=\cA _n(y,\mu ^{-1}hD_x), \label{18-5-3}\\
\shortintertext{then}
\cR _n \Def \cA _n-\cB _n =(x_2-y_2) \cQ _n(x,y,\mu ^{-1}hD_x)
\label{18-5-4}
\end{gather}
and then
\begin{multline}
U_n^\pm (x_2,y_2,t)\Def \uptheta (\pm t) U(x_2,y_2,t) =
\mp i h \cG _n^\pm \updelta (x_2-y_2)\updelta (t)=\\[3pt]
\mp i h \bar{\cG}_n^\pm \updelta (x_2-y_2)\updelta (t)
\mp i h \bar{\cG}_n^\pm \cR _n \cG _n^\pm \updelta (x_2-y_2)\updelta(t)
\label{18-5-5}
\end{multline}
where $\bar{\cG}^\pm$ and $\cG ^\pm$ are forward/backward parametrices for the Cauchy problem for operators $hD_t-\cB _n$ and $hD_t-\cA _n$ respectively. Similarly
\begin{multline}
U_n^\pm (x_2,y_2,t)\,^t\!Q_y =
\mp i h \cG _n^\pm \updelta (x_2-y_2)\,^t\!Q_y \updelta (t)=\\[3pt]
\mp i h \bar{\cG}_n^\pm \updelta (x_2-y_2)\,^t\!Q_y \updelta (t)
\mp i h \bar{\cG}_n^\pm\cR _n\cG _n^\pm\delta (x_2-y_2)Q_y^t
\updelta (t)
\label{18-5-6}
\end{multline}
for $Q_y= Q(y_2,\mu^{-1}hD_{y_2})$ where as usual we write pseudo-differential operator acting on $y$ on the right in accordance to matrix theory and $\,^t\!Q$ denotes the transposed operator.

So, $U_n = U_n^0 +U_n^1$ where $U^0$ is a Schwartz kernel of the propagator for 
$\cB _n$ and
\begin{equation}
U_n^{1\pm} = \bar{\cG}_n^\pm \cR _n U_n^\pm.
\label{18-5-7}
\end{equation}

As usual we take $(x_2-y_2)$ and carry it forward to the left using commuting relation
\begin{equation}
[\bar{\cG}_n^\pm ,(x_2-y_2)] = \bar{\cG}_n^\pm [\cB _n, x_2-y_2] \bar{\cG}_n^\pm
\label{18-5-8}
\end{equation}
and in the very front this factor will be killed by $\Gamma '$. Then for our purpose we can replace $U_n^{1\pm}$ by
\begin{equation}
\bar{\cG}_n^\pm [\cB _n, x_2-y_2] \bar{\cG}_n^\pm \cQ _n U_n^{0\pm}.
\label{18-5-9}
\end{equation}
Let us dump it in (\ref{18-4-28}). Consider first
\begin{equation}
F_{t \to h^{-1}\tau} \biggl(\varphi_T(t)\Gamma '
\Bigl(\bigl(U^1_n (\tilde{\psi} \,^t\!q_{(\iota)}^\w )_y\bigr)\Bigr)\biggr)
\label{18-5-10}
\end{equation}
with \underline{either} $\varphi_T=\bar{\chi}_T$, $T=\bar{T}_*$ \underline{or}
$\varphi_T= \chi_T$, $T\ge \bar{T}_0$; recall that
$\bar{T}_*= Ch/\varepsilon$ as $\nu=1$.

In the first case ($\varphi_T=\bar{\chi}_T$, $T=\bar{T}_*$) our estimate is rather straightforward. Note that the norm of $\bar{\cG }_n^\pm $ in the strip $\{|t|\le T\}$ does not exceed $h^{-1}T$ and the norm of $[\cB _n, x_2-y_2]$ does not exceed $C\mu ^{-1}h$. Then the factor $(x_2-y_2)$ translates into an extra factor $C \mu^{-1}h^{-1}T^2$ in the norm estimate. For $T=\bar{T}_*$ this factor is $C\mu ^{-1}h \varepsilon ^{-2} = C|\log h|^{-1}$ and in the general case it will be $C|\log h|^{-1}T^2/\bar{T}_*^2$.

Further, we can assume that 

\begin{claim}\label{18-5-11}
The operator norm of $\cR _n $ does not exceed
$C\bigl(\vartheta (\varepsilon) + \varepsilon^2\bigr)$.
\end{claim}
Really, we always can make $\partial _{x_2} W$ vanish at some point of 
$\supp q_{(\iota)}$ by an appropriate Fourier integral operator corresponding to some linear symplectomorphism at $(x_2,\xi_2)$-plane. 

Then in this case expression (\ref{18-5-10}) does not exceed
$C\mu \varepsilon^{-1}\bigl(\vartheta (\varepsilon) + \varepsilon^2\bigr)$; recall that for $U^0$ it was just $C \mu $.

In the second case ($\varphi_T= \chi_T$, $T\ge \bar{T}_*$) we, using the standard rescaling arguments, gain an extra factor $(\bar{T}_*/ T)^s$ with arbitrarily large exponent $s$ which takes care of all $(T/ \bar{T}_*)^k$ factors and summation with respect to partition of unity with respect to $t$. Now we can conclude that expression (\ref{18-5-10}) does not exceed 
$C\mu \varepsilon^{-1}\bigl(\vartheta (\varepsilon) + \varepsilon^2\bigr)$ for 
$\varphi =\bar{\chi}$, $T=T_*$. 

Then applying arguments we used to estimate the total contribution of our ``special'' elements $q_{(\iota)}$ to (\ref{18-4-28}), we see that 

\begin{claim}\label{18-5-12}
The total contribution of these ``special'' elements $q_{(\iota)}$ to (\ref{18-4-28}) with $U$ replaced by $U^1$ does not exceed $Ch^{-2}\vartheta(\varepsilon)+Ch^{-2}\varepsilon^2$
\end{claim}
what is exactly what we expect mollification error would be as $(l,\sigma)\preceq (2,0)$.

To deal with the case $(l,\sigma)\succ (2,0)$ we note that due to propagation speed $O(\mu^{-1})$ we can estimate $(x_2-y_2)$ by $C\mu^{-1}T$ and as effectively $T$ is estimated by $\bar{T}_*$ we conclude that we can assume that 

\begin{claim}\label{18-5-13}
As $(l,\sigma)\succ (2,0)$ the operator norm of $\cR _n $ does not exceed
$C\varepsilon^{\prime 2}$ with
$\varepsilon' \Def C\mu^{-1}h/\varepsilon$
\end{claim}
and replace (\ref{18-5-12}) by
\begin{claim}\label{18-5-14}
As $(l,\sigma)\succeq (2,0)$ the total contribution of these ``special'' elements $q_{(\iota)}$ to (\ref{18-4-28}) with $U$ replaced by $U^1$ does not exceed $Ch^{-2}\varepsilon^{\prime\,2}=O(\mu^{-1}h^{-1})$.
\end{claim}

Now let us consider contribution of $U^0$. One can see easily that it is exactly
equal to the left-hand expression of
\begin{multline}
\Gamma'\Bigl(
\uptheta \bigl(\tau - \cA _n(y_2,\mu ^{-1}hD_{x_2})\bigr)K_{(\iota)} \Bigr) =\\
(2\pi)^{-1}\mu h^{-1} \iint
\uptheta \bigl(\tau - \cA _n(x_2,\xi_2) \bigr) k_{(\iota)}(x_2,\xi_2)
\,dx_2 d\xi_2
\label{18-5-15}
\end{multline}
where $K_{(\iota)} = \tilde{\psi}_n q_{(\iota)}^\w$, $k_{(\iota)}$ is its
$pq$-symbol\footnote{\label{foot-18-25} Which means that 
$K= k({\overset 1 {\mathstrut{x}}},{\overset 2 {\mathstrut{\mu ^{-1}hD_x}}} )$.} and
\begin{equation}
\tilde{\psi}_n =
\langle \tilde{\psi} \Upsilon _n, \Upsilon _n\brangle
\label{18-5-16}
\end{equation}
where $\Upsilon _n=\Upsilon _n(x_1)$ and $\langle .,. \brangle$ is an inner
product in $\sL^2(\bR_{x_1})$.

The problem is that $pq$-symbols are not invariant with respect to Fourier integral operators and the difference between $pq$- and Weyl symbols is not small because of $q_{(\iota)}$. However, there is a simple walk-around. Cover our total $(x_2,\xi_2)$-domain by two, with $\partial _{\xi_2}W$ and 
$\partial _{x_2}W$ disjoint from $0$ respectively. We consider the first domain; the second case will be reduced to the first one easily, by means of  $\mu^{-1}h$-Fourier transform.

We can consider partition of unity (depending on $n$)
$q'_{(\iota_1)}(x_2)q'_{(\iota_2)}(\xi _2)$ such that for each
$\iota_1$ (and each $n$) there is no more than one number $\iota_2$ such that
condition (\ref{18-4-20}) holds on 
$\supp q'_{(\iota_1)}\times \supp q''_{(\iota_2)}$. Then
all the contributions of all other elements do not depend on this Fourier integral operator and summing with respect to $\iota_2$ we get just $q'_{(\iota_1)}(x_2)$. 

But we can take
the required linear symplectomorphism
$(x_2,\xi_2)\mapsto (x_2,\xi_2 - r_{(\iota_1)} x_2)$ and the corresponding Fourier integral operator will be just multiplication by 
$e^{-i\mu h^{-1}r_{(\iota)}x_2^2/2}$ which does not affect $q'_{(\iota_1)}(x_2)$. On the other hand, for $\tilde{\psi}$ the difference between different symbols does not exceed
$C\mu ^{-1}h$ and therefore without increasing our remainder estimate we can
rewrite our final answer as
\begin{equation}
(2\pi)^{-1}\mu h^{-1} \sum_{n\ge 0}
\iint \uptheta \bigl( \tau -\cA _n(x_2,\xi_2)\bigr)
\tilde{\psi} _n(x_2,\xi_2) \,dx_2 d\xi_2.
\label{18-5-17}
\end{equation}
Recall, that so far $\mu \le (h|\log h|)^{-1}$ and assumption (\ref{18-2-15}) was supposed to be fulfilled.

\subsection{Step I under assumption (\ref{18-2-16})}
\label{book_new-sect-18-5-1-2}

If assumption (\ref{18-2-15}) is not fulfilled we must observe that instead of claim (\ref{18-5-12}) we now have that 

\begin{claim}\label{18-5-18}
The total contribution of these ``special'' $\nu$-elements $q_{(\iota)}$ to (\ref{18-4-28}) with $U$ replaced by $U^1$ does not exceed 
$Ch^{-2}\bigl(\vartheta(\varepsilon)+\varepsilon^2\bigr)\times 
\bigl(\nu^2+ \mu h \bigr)$.
\end{claim}

Really, in comparison with the case $\nu=1$, factor $T_*$ (or $\bar{T}_*$) brings a factor $\nu^{-1}$, and volume brings factor $\nu$ but the number of indices ``$n$'' is now $\asymp(\nu^2/(\mu h) +1 )$ rather than 
$\asymp 1/(\mu h)$. 

Further, claim (\ref{18-5-13}) is preserved as factor $\nu$ in the propagation speed is compensated by factor $\nu^{-1}$ in $\bar{T}_*$ and therefore (\ref{18-5-14}) should be replaced by 

\begin{claim}\label{18-5-19}
As $(l,\sigma)\succeq (2,0)$ the total contribution of these ``special'' $\nu$-elements $q_{(\iota)}$ to (\ref{18-4-28}) with $U$ replaced by $U^1$ does not exceed $Ch^{-2} \varepsilon^{\prime\,2}\times \bigl(\nu^2+ \mu h )$.
\end{claim}

Then we have an extra term
\begin{equation}
C\mu h^{-1} \bigl(\vartheta(\varepsilon) + \varepsilon ^{\prime\,2}\bigr)
\label{18-5-20}
\end{equation}
which after summation with respect to $\nu$ under assumption (\ref{18-2-16}) just acquires factor $|\log h|$; further, under assumption \ref{18-2-16-+} only $\nu$ of the same magnitude should be taken in account near each critical point and no logarithmic factor is acquired. Thus we arrive to

\begin{claim}\label{18-5-21}
Under assumption (\ref{18-2-16}) the total contribution of these ``special'' elements $q_{(\iota)}$ to (\ref{18-4-28}) with $U$ replaced by $U^1$ does not exceed 
\begin{multline*}
C\bigl(h^{-2}\vartheta(\varepsilon)+ \mu^{-1}h^{-1}\bigr)
\bigl(1+\mu h |\log h|\bigr)\asymp\\
C\bigl(h^{-2}\vartheta(\varepsilon)+ \mu^{-1}h^{-1}\bigr) +
C\bigl(\mu h^{-1}\vartheta(\varepsilon)+ 1\bigr)|\log \mu|
\end{multline*}
\end{claim}\vglue-12pt\noindent
and
\begin{claim}\label{18-5-22}
Under assumption \ref{18-2-16-+} the total contribution of these ``special'' elements $q_{(\iota)}$ to (\ref{18-4-28}) with $U$ replaced by $U^1$ does not exceed 
$C\bigl(h^{-2}\vartheta(\varepsilon)+ \mu^{-1}h^{-1}\bigr)$. 
\end{claim}

Sure, in our case $\mu \lesssim (h|\log h|)^{-1}$ these extra terms with logarithmic factor do not matter but they will be important otherwise.

\subsection{Step II under assumption (\ref{18-2-15})}
\label{book_new-sect-18-5-1-3}

To finish the proof of main theorems we need to bring expression (\ref{18-5-17}) with $\tau =0$ to the more explicit form. 

First of all, we need to use original $(x_1,x_2)$-coordinates rather than $(x_2,\xi_2)$-coordinates after reduction. The natural expression of the area form in $(x_1,x_2)$-coordinates is $\omega = \sqrt g dx_1\wedge dx_2$ which can be considered as the area form on $\Sigma =\{p_1=p_2=0\}$ such that 
$\omega \wedge dp_1\wedge dp_2 = dx\wedge d\xi$.

Then after construction of subsubsections~\ref{book_new-sect-18-3-1-1}.1 and \ref{book_new-sect-18-3-1-2}.2 we have a symplectic map $\Phi $ which transforms 
$\Sigma $ into  $\Sigma _0=\{x_1=\xi_1=0\}$; let us denote $\Phi |_\Sigma $ by $\phi $. Further, $\Phi $ transforms $g^{-\frac{1}{2}} dp_1\wedge dp_2$ at points of $\Sigma $ exactly into $-dx_1\wedge d\xi_1$ at points of $\Sigma_0$. Therefore
\begin{equation}
\phi_*\omega =-dx_2\wedge d\xi_2,\qquad \phi=\Phi |_\Sigma.
\label{18-5-23}
\end{equation}

Consider next transformations of subsection~\ref{book_new-sect-18-3-1-3}.3. If $V=0$ these transformations will preserve $\Sigma _0$ and be identical on it. In the general case one can check easily that these transformations transform point $(0,0;x_2,\xi_2)\in\Sigma _0 $ into 
\begin{multline*}
\bigl(\zeta _1(x_2,\xi_2),\zeta _2(x_2,\xi_2);
x_1+\eta _1(x_2,\xi_2),x_2 + \eta _2(x_2,\xi_2)\bigr)\in \Sigma'_0\\[3pt]
\text{with\ \ } |\zeta _j| +|\nabla \zeta_j|\le C\mu^{-2},\ 
|\eta _j| +|\nabla \eta_j|\le C\mu^{-4}+
C\varepsilon ^l|\log \varepsilon |^{-\sigma}
\end{multline*}
and therefore after projecting
$\Sigma' _0$ onto $\Sigma_0$: $(x_1,\xi_1;x_2,\xi_2)\mapsto (0,0;x_2,\xi_2)$ we
have $(1+k)dx_2\wedge d\xi_2$ with $k=0$ unless $(l,\sigma)\succeq (2,0)$ and
$\mu $ is close to $h^{-\frac{1}{3}}$ and $k=\mu^{-4}k_1(x)$ in the latter case.

Therefore, instead of $-dx_2\wedge d\xi_2 $ in our formula we can use
$(1+k)\sqrt g\, dx_1dx_2$ and instead of $W(x_2,\xi_2,\rho _n)$ we can use
$W_n\Def W(\phi (x_1,x_2), r_n)$.

Further, obviously one can find $w(x_1,x_2)$
of the same regularity as $W$ such that $(2n+1)\mu h + W_n <0$ if and only if
$(2n+1)\mu h + w <0$ and replace $W_n(x)$ by $w(x)$. One should not be very concerned with the difference between $\psi $ and $\psi _n$ because we are interested in the calculations when $\psi =1$ and then $\psi_n=1$ as well.

Anyway, one can prove easily that
\begin{equation*}
\tilde{\psi}_n (x_2,\xi_2) =
\psi \circ \phi ^{-1} + \mu^{-2}\psi ' _n+\mu^{-4}\psi''_n
\end{equation*}
where $\psi'_n$, $\psi''_n$ are defined by (\ref{18-5-16}) with $\tilde{\psi}$ replaced by $\psi '$, $\psi''$ respectively, and with symbols $\psi '$, $\psi''$ belonging to $\sF^{l-1,\sigma}, \sF^{l-2,\sigma}$ and not depending
on $\mu,h$ and one can take $\psi''_n$, $\psi''$ equal to $0$ unless
$(l,\sigma)\succeq (2,0)$ and $\mu $ is close to $h^{-\frac{1}{3}}$.

Then we get the final formula with the correction term
\begin{multline}
h^{-2} \cN ^\MW_{2,\corr} =\\
\frac{1}{2\pi} h^{-2} \int \sum_n \Bigl( \uptheta \bigl(-w-(2n+1)\mu h\bigl) -
\uptheta \bigl( -\frac{V}{F} -(2n+1)\mu h) \Bigr)\mu h \cdot \psi (x)\, dx +\\
\mu^{-1}h^{-1} \int \sum_n 
\uptheta \bigl( -w-(2n+1)\mu h\bigl) (\psi ' _n\circ \phi )\,dx+\\
\mu^{-3}h^{-1}\int \sum_n \uptheta \bigl( -w-(2n+1)\mu h\bigl)
(\psi '' _n\circ \phi + k)\,dx.
\label{18-5-24}
\end{multline}

\medskip
Consider the second term in (\ref{18-5-24}). It is equal 
\begin{equation}
\mu ^{-1}h^{-1}\int 
\Bigl(\Tr_1 \bigl(\Pi_{-w} \psi '(x_1,x_2,\mu^{-1}hD_1, \xi_2)\bigr) \circ \phi\Bigr) \,dx
\label{18-5-25}
\end{equation}
where $\Pi _\tau$ is the spectral projector in $\sL^2(\bR)$ ($\bR=\bR_{x_1}$) for harmonic oscillator $\mu ^2 x_1^2 + h^2 D_1^2$ and $\Tr_1$ is an operator trace in $\sL^2(\bR)$, $\tau=-w$ is the spectral parameter.

One can apply Weyl formula to calculate this expression (\ref{18-5-25}).  The error will not exceed
$C\mu ^{-1}h^{-1} \times \bigl(1+ \varepsilon^{l-2}|\log h|^{-\sigma}\bigr)$
where the last factor is the remainder estimate for Weyl formula for
$\Tr_1 \bigl(\Pi_{-w} \psi '(x_1,hD_1)\bigr)$\,\footnote{\label{foot-18-26} Actually, $C$ is the remainder estimate and
$C\varepsilon^{l-2}|\log h|^{-\sigma}$ is an estimate for error in the
calculations.}. Therefore this error does not exceed desired
$Ch^{-2}\bigl(\mu ^{-1}h + \vartheta (\varepsilon)\bigr)=C\mu^{-1}h^{-1}+ Ch^{-2}\vartheta (\varepsilon)$.

On the other hand, the Weyl answer (for the same term) is
\begin{gather}
\mu h^{-1} \iint \uptheta \bigl(-w -r^2\bigl)
\psi' (x,\mu^{-1}r\cos \alpha, \mu^{-1}r \sin\alpha)\,r dr d\alpha ;
\label{18-5-26}\\
\intertext{plugging into (\ref{18-5-25}) we get}
h^{-2} \iiint \uptheta \bigl(-w -r^2\bigl)
\psi' (x,\mu^{-1}r\cos \alpha, \mu^{-1}r \sin\alpha)\,r dr d\alpha dx.
\label{18-5-27}
\end{gather}
Recall that 
\begin{equation}
w -V =O\bigl(\vartheta (\mu ^{-1})\bigr)
\label{18-5-28}
\end{equation}
and then due to assumption (\ref{18-2-15}) replacing $w$ by $V$ we make an error not exceeding $C\mu ^{-2-l}|\log h|^{-\sigma}h^{-2}$; this expression does not exceed $Ch^{-2}\vartheta (\varepsilon)$ as $l\le 2$ and 
$C\mu ^{-1}h^{-1}$ as $l>2$; recall that $\mu \ge (h|\log h|)^{-\frac{1}{3}}$.

\medskip
We can treat the third term in (\ref{18-5-24}) in the same way. Therefore

\begin{claim}\label{18-5-29}
Without deterioration of the remainder estimate, the second and third
terms in (\ref{18-5-24}) can be rewritten as $K_1\mu^{-2}h^{-2}$ and
$K_2\mu^{-4}h^{-2}$ respectively with some constants $K_1$, $K_2$.
\end{claim}

We also claim that

\begin{claim}\label{18-5-30}
With an error not exceeding $C\mu^{-1}h^{-1}+ Ch^{-2}\vartheta (\mu h)$
one can replace in the first term in (\ref{18-5-24}) the Riemann sum by the corresponding integral.
\end{claim}

Really, without any loss of the generality one can assume that
$|\partial_{x_j}V|\ge \epsilon_0$ with either $j=1$ or $j=2$ and the same is
true for $w$. Then the first part of the first term in (\ref{18-5-24}) is equal to
\begin{multline}
h^{-2} \int \sum_n \mu h \uptheta \bigl(\tau -w-(2n+1)\mu h\bigl)
(\partial_{x_j}w ) (\partial_{x_j}w' )^{-1}\psi (x)\, dx +\\
h^{-2} \mu h \int \sum_n \uptheta \bigl(\tau -w-(2n+1)\mu h\bigl)
(\partial_{x_j} w )
\bigl((\partial_{x_j} w )^{-1}- (\partial_{x_j}w' )^{-1}\bigr)\psi (x)\, dx
\label{18-5-31}
\end{multline}
where $w'$ is $\mu h$-mollification of $w$. Note that
$\bigl((\partial_{x_j} W_n )^{-1}- (\partial_{x_j}W' )^{-1}\bigr)$
does not exceed $C(\mu h)^{l-1}|\log \mu h|^{-\sigma}$ and therefore
passing in the corresponding term from the Riemann sum to the corresponding
integral brings an error not exceeding
$Ch^{-2}\times \mu h \times (\mu h)^{l-1}|\log \mu h|^{-\sigma}$.
On the other hand, the first term in (\ref{18-5-31}) is equal to
\begin{equation*}
h^{-2} \int \sum_n \bigl(\tau -w-(2n+1)\mu h\bigl)_+
\bigl(\partial_{x_j} (\partial_{x_j}w' )^{-1}\psi (x)\bigr)\mu h \, dx
\end{equation*}
with $\bigl(\partial_{x_j} (\partial_{x_j}w' )^{-1}\psi (x)\bigr)$ not exceeding
$(\mu h)^{l-2}|\log \mu h|^{-\sigma}$ and passing from the Riemann sum to the
corresponding integral brings an error not exceeding
$Ch^{-2}\times (\mu h )^2\times
(\mu h)^{l-2}|\log \mu h|^{-\sigma}$ for $(l,\sigma)\preceq (2,0)$; for
$(l,\sigma)\succeq (2,0)$ we need to repeat this procedure one more time to
recover the same estimate.

Thus, claim (\ref{18-5-30}) has been proven and with an error not exceeding
$C\bigl(\mu^{-1}h^{-1}+ h^{-2}\vartheta (\mu h)\bigr)$ the first term in (\ref{18-5-24}) is equal to
\begin{equation*}
\frac{1}{4\pi} h^{-2}
\int \bigl(  w _- - V _-\bigr) \psi (x)\, dx =
\frac{1}{4\pi} h^{-2} \int (V-w)\psi (x)\, dx
\end{equation*}
due to condition (\ref{18-2-27}).

So we get that with the indicated error the correction term is
\begin{equation}
\frac{1}{4\pi}h^{-2} \int (V-w)\psi (x)\, dx +
h^{-2}\bigl(K_1\mu^{-2}+K_2\mu^{-4}\bigr)
\label{18-5-32}
\end{equation}
with $V$, $K_1$, $K_2$ not depending on $\mu^{-1}$ and $w$ depending on it.

Comparing with the results for the weak magnetic field (section~\ref{book_new-sect-18-2}) when correction term is $0$\,\footnote{\label{foot-18-27} Surely for
$\mu \asymp(h|\log h|)^{- \frac{1}{3}}$ we can take the same mollification but we also removed operator with $A_0$-bound not exceeding $\vartheta(\varepsilon)$ and one can see easily that it brings error not exceeding 
$Ch^{-2}\vartheta (\varepsilon)$.} we see that expression (\ref{18-5-32})
vanishes modulo $O\bigl(\mu^{-1}h^{-1}+ \vartheta (\varepsilon)h^{-2}\bigr)$,
with $\varepsilon = (\mu ^{-1}h|\log h|)^{\frac{1}{2}}$ and then (\ref{18-5-24}) transforms into
\begin{multline}
h^{-2}\cN_{2,\corr}\Def \\
h^{-2} \int \sum_n \Bigl( \uptheta \bigl(\tau -W -(2n+1)\mu h\bigl) -
\uptheta \bigl(\tau -V -(2n+1)\mu h) \Bigr)\mu h\cdot  \psi (x)\, dx -\\
\frac{1}{2} h^{-2} \int (V-W)\psi (x)\, dx
\label{18-5-33}
\end{multline}
where we replaced notation $w$ by $W$. Thus we arrive under assumption (\ref{18-2-15}) to proposition~\ref{prop-18-5-1} below.

\subsection{Step II under assumption (\ref{18-2-16})}
\label{book_new-sect-18-5-1-4}

Note that up to (\ref{18-5-31}) we have not used any non-degeneracy assumption. Next step is done on each $\nu$-ball with 
$\nu \ge C_0 \max\bigl( (\mu h)^{\frac{1}{2}},\mu^{-1}\bigr)$ and then the factor $\nu^{-2}$ appears (in comparison with the calculations of the previous subsubsection) and it is compensated by factor $\nu^2$ from the ball volume; summation over all balls brings either factor $\log \mu$  under assumption (\ref{18-2-16}).

However the same factor $\log \mu$ or $1$ appears under corresponding assumption in the remainder estimate as well and therefore without deterioration of the remainder we arrive to the same expression (\ref{18-5-33}) for the correction term, arriving under assumption (\ref{18-2-16}) to proposition~\ref{prop-18-5-1}:

\begin{proposition}\label{prop-18-5-1}
Let $h^{-\frac{1}{3}}\lesssim \mu \lesssim (h|\log h|)^{-1}$ and condition \textup{(\ref{18-2-27})} be fulfilled. Then 

\medskip\noindent
(i) Under  non-degeneracy assumption \textup{(\ref{18-2-15})} Tauberian expression \textup{(\ref{18-4-28})}  equals 
\begin{equation}
h^{-2} \int \bigl( \cN^\MW _x(0) +  \cN^\MW _{x,\corr}\big)\psi (x)\,dx
\label{18-5-34}
\end{equation}
modulo $O\bigl(\mu^{-1}h^{-1}+ h^{-2} \vartheta (\varepsilon)\bigr)$ with $h^{-2}\cN^\MW _{x,\corr}$ defined according to \textup{(\ref{18-5-33})}.

\medskip\noindent
(ii) Under  non-degeneracy assumption \textup{(\ref{18-2-16})} Tauberian expression \textup{(\ref{18-4-28})}  equals \textup{(\ref{18-2-34})} modulo $O\bigl(\mu^{-1}h^{-1}+ h^{-2} \vartheta (\varepsilon)|\log h|\bigr)$.
\end{proposition}

\subsection{Calculations: $\mu \gtrsim h^{-1}|\log h|^{-1}$}
\label{book_new-sect-18-5-1-5}

In this case calculations are easy. We just should remember that as $\mu \gtrsim h^{-1}$ the correction $h^2W'$ with $W'$ defined by (\ref{18-5-42})--(\ref{18-5-43}) below to potential.

\subsection{Comparisons}
\label{book_new-sect-18-5-1-6}

We need to compare 
$\int \tilde{\cN}^\MW(x,\tau)\psi (x)\, dx$ with 
$\int \cN ^\MW(x,\tau)\psi (x) \, dx$:

\begin{proposition} \label{prop-18-5-2}
(i) Under non-degeneracy assumptions \textup{(\ref{18-2-15})} or \ref{18-2-16-+} as $\mu \lesssim h^{-1}$ and
\ref{18-2-15-*} or \ref{18-2-16-*+} as $\mu\gtrsim h^{-1}$ 
\begin{equation}
h^{-2}|\int \Bigl(\tilde{\cN}^\MW(x,\tau )-
\cN^\MW(x,\tau)\Bigr)\psi (x)\, dx |\le 
C(1+\mu h) h^ {-2}\vartheta (\varepsilon )
\label{18-5-35}
\end{equation}
for $|\tau |\le \epsilon$;

\medskip\noindent
(ii) Under non-degeneracy assumptions \textup{(\ref{18-2-16})} as 
$\mu \lesssim h^{-1}$ and \ref{18-2-16-*} as $\mu\gtrsim h^{-1}$ 
\begin{multline}
h^{-2}|\int \Bigl(\tilde{\cN}^\MW(x,\tau )-
\cN^\MW(x,\tau)\Bigr)\psi (x)\, dx |\le \\
C(1+\mu h) h^ {-2}\vartheta (\varepsilon )|\log \varepsilon|
\label{18-5-36}
\end{multline}
for $|\tau |\le \epsilon$.
\end{proposition}

\begin{proof} Without any loss of the generality we can assume that $F=1$.

\medskip\noindent
(i) Assume first that $\vartheta (\varepsilon )\le \epsilon \mu h$ and assumption (\ref{18-2-15}) or \ref{18-2-16-+} is fulfilled.

Let us introduce $\gamma $-admissible partition of unity
with $\gamma =\epsilon \vartheta (\varepsilon )$. If
$|V+(2j+1)\mu h F-\tau |\ge C\gamma $ for all $j\in \bZ^+$ on some element of
the partition, then on this element $\tilde{\cN}^\MW (x, \tau )$ and
$\cN ^\MW (x, \tau )$ coincide while on all other elements their difference
does not exceed $C\mu h^{-1}$. However, due to \textup{(\ref{18-2-15})} or \ref{18-2-15-*} the total area of the latter elements does not exceed 
$C\gamma (1+\mu ^{-1}h^{-1})$ which yields the necessary estimate.

\medskip\noindent
(ii) Assume now that $\vartheta (\varepsilon )\ge \epsilon \mu h$ and assumption (\ref{18-2-15}) is fulfilled. Then $\uptheta \bigl(\tau -V -(2n+1)\mu h\bigr)$ and  $\uptheta \bigl(\tau -\tilde{V}-(2n+1)\mu h\bigr)$ differ for each $x$ for no more than $C_0\vartheta (\varepsilon )/\mu h+1$ number of indices ``$n$'' which yields the necessary estimate.

\medskip\noindent
Similar analysis under one of the assumptions \textup{(\ref{18-2-16})}--\ref{18-2-16-*+} we leave to the reader.
\end{proof}
  
We also need to compare $W$ and $V/F$ as $\mu \lesssim (h|\log h|)^{-1}$ (i.e. as $\mu ^{-1}\gtrsim \varepsilon$). One can see easily that it is $O\bigl(\vartheta (\mu^{-1})+\mu^{-2}\bigr)$ and we arrive to 

\begin{proposition} \label{prop-18-5-3}
Let $\mu \lesssim (h|\log h|)^{-1}$. Then

\medskip\noindent
(i) Under assumption \textup{(\ref{18-2-15})} or \ref{18-2-16-+} 
\begin{equation}
|\cN^\MW _{2,\corr}|\le 
C\left\{\begin{aligned}
&\vartheta (\mu h)\qquad &&\text{as\ \ }\mu \le h^{-\frac{1}{2}},\\
&\vartheta (\mu ^{-1})\qquad &&\text{as\ \ }\mu \ge h^{-\frac{1}{2}},\ (l,\sigma)\preceq (2,0),\\
&\mu^{-4}h^{-2}\vartheta (\mu h)\qquad&&\text{as\ \ }\mu \ge h^{-\frac{1}{2}},\ (l,\sigma)\succ (2,0).
\end{aligned}\right.
\label{18-5-37}
\end{equation}
In particular, $h^{-2}\cN^\MW _{2,\corr}=O(\mu^{-1}h^{-1})$ for sure as 
$(l,\sigma) \succeq (3,0)$.

\medskip\noindent
(ii) Under assumption \textup{(\ref{18-2-16})} 
\begin{multline}
|\cN^\MW _{2,\corr}|\le \\
C|\log h| \left\{\begin{aligned}
&\vartheta (\mu h)\qquad &&\text{as\ \ }\mu \le h^{-\frac{1}{2}},\\
&\vartheta (\mu ^{-1})\qquad &&\text{as\ \ }\mu \ge h^{-\frac{1}{2}},\ (l,\sigma)\preceq (2,0),\\
&\mu^{-4}h^{-2}\vartheta (\mu h)\qquad&&\text{as\ \ }\mu \ge h^{-\frac{1}{2}},\ (l,\sigma)\succ (2,0).
\end{aligned}\right.
\label{18-5-38}
\end{multline}
In particular, $h^{-2}\cN^\MW _{2,\corr}=O(\mu^{-1}h^{-1})$ for sure as 
$(l,\sigma) \succeq (3,1)$.
\end{proposition}

\begin{proof}
An easy proof is left to the reader; as $(l,\sigma)\succ (2,0)$ we need to take into account some extra correction terms of the type $K\mu^{-2}h^{-2}$ because then $W-V/F = k \mu^{-2} +O\bigl(\vartheta\mu^{-1})\bigr)$.
\end{proof}

\section{Main theorems}
\label{book_new-sect-18-5-2}

Now, assembling results of this section we arrive to our main

\begin{theorem}\label{thm-18-5-4} 
Let $d=2$ and $A$ be a self-adjoint in $\sL^2(X)$ Schr\"odinger operator defined by \textup{(\ref{18-1-1})}.

Let  $g^{jk}, F,\psi \in \sC^{\bar{l},\bar{\sigma}}$, 
$V\in \sC^{l,\sigma}$ with $(\bar{l},\bar{\sigma})\succeq (2,0)$, 
$(\bar{l},\bar{\sigma})\succeq (l,\sigma)\succeq (1,2)$ and let conditions \textup{(\ref{18-1-6})} and \textup{(\ref{18-2-27})} be fulfilled in
$B(0,1)\subset X\subset \bR ^2$.

Let $ (h|\log h|)^{-{\frac{1}{3}}}\le \mu \le Ch^{-1}$. Then there are two framing  approximations $\tilde{A}=\tilde{A}^\pm$  as in footnote~\footref{foot-18-16} such that for each of them 

\medskip\noindent
(i) Under assumption \textup{(\ref{18-2-15})} estimate
\begin{multline}
\R^\MW_*   \Def \\ |\int \Bigl( \tilde{e}(x,x,0)-h^{-2}\cN_{2,x} ^\MW (0)-
h^{-2}\cN_{2,x\,\corr } ^\MW (0) \Bigr)\psi (x)\, dx|\le \\
C\mu ^{-1}h ^{-1}+ C\vartheta (\varepsilon) h^{-2} 
\label{18-5-39}
\end{multline}
holds with the  standard magnetic Weyl expression $h^{-2}\cN_{2,x}^\MW (\tau)$ and with the correction term  $h^{-2}\cN _{2,x\,\corr } ^\MW(\tau)$ defined by \textup{(\ref{18-5-31})} with the appropriate corrected potential  
$W=V/F + O\bigl(\vartheta (\mu^{-1})+ \mu^{-2}\bigr)$; here and below 
$\varepsilon = (\mu^{-1}h|\log h|)^{\frac{1}{2}}$;

\medskip\noindent
(ii) Let  $(l,\sigma)\succeq (2,0)$; then under assumption  \ref{18-2-16-+} the same estimate \textup{(\ref{18-5-39})} holds;

\medskip\noindent
(iii) Let  $(l,\sigma)\succeq (2,0)$; then under assumption \textup{(\ref{18-2-16})} 
\begin{equation}
\R^\MW_* \le \\
C\mu ^{-1}h ^{-1}+ C\vartheta (\varepsilon) h^{-2}|\log h| + C|\log h|
\label{18-5-40}
\end{equation}
\end{theorem}

\begin{remark}\label{rem-18-5-5}
Recall that the correction term has been estimated in proposition~\ref{prop-18-5-3}. In particular it implies estimates for $\R^\MW$ without correction term.
\end{remark}

\begin{Problem}\label{problem-18-5-6}
Either improve estimates (\ref{18-5-37}), (\ref{18-5-38}) of proposition~\ref{prop-18-5-3} or construct counter-examples showing that the improvement is impossible.
\end{Problem}

\begin{theorem}\label{thm-18-5-7}
Let $d=2$ and $A$ be a self-adjoint in $\sL^2(X)$ Schr\"odinger-Pauli operator defined by \textup{(\ref{18-1-29})}.

Let  $g^{jk}, F,\psi \in \sC^{\bar{l},\bar{\sigma}}$, 
$V\in \sC^{l,\sigma}$ with $(\bar{l},\bar{\sigma})\succeq (2,0)$, 
$(\bar{l},\bar{\sigma})\succeq (l,\sigma)\succeq (1,1)$ and condition \textup{(\ref{18-1-6})} be fulfilled in $B(0,1)\subset X\subset \bR ^2$.

Then there are two framing  approximations $\tilde{A}=\tilde{A}^\pm$ as in footnote~\footref{foot-18-16} such that for each of them 

\medskip\noindent
(i) Under assumption \ref{18-2-15-*} estimate
\begin{equation}
\R^\MW_* \le 
C+ C\vartheta (\varepsilon) \mu h^{-1} + C\bar{\vartheta}(\varepsilon)\mu^2
\label{18-5-41}
\end{equation}
holds with the correction term $h^{-2}\cN_{2,x\,\corr }^\MW$ defined again by \textup{(\ref{18-5-24})}\footnote{\label{foot-18-28} However now $W=V+W'h^2$ with
\begin{equation}
W'= - \frac{1}{12} \kappa \bigl(\fz^2+2 \bigr) +
F^{\frac{1}{2}}\sum_{j,k}\partial_j\bigl(g^{jk}\partial_kF^{-\frac{1}{2}}\bigr)
\label{18-5-42}
\end{equation}
and $\kappa $ a curvature associated with metrics $g^{jk}F^{-1}$ (and this correction is not due to the lack of the smoothness); 
\begin{equation}
g^{jk}F^{-1}=\alpha ^2 \updelta_{jk}\implies \kappa = 2 \Delta \log \alpha .
\label{18-5-43}
\end{equation}
Further, the sum in the definition contains only a finite number of terms.}.

\medskip\noindent
(ii) Let  $(l,\sigma)\succeq (2,0)$; then under assumption \ref{18-2-16-*+} the same estimate \textup{(\ref{18-5-40})} holds;

\medskip\noindent
(iii) Let  $(l,\sigma)\succeq (2,0)$; then under assumption \ref{18-2-16-*} 
\begin{equation}
\R^\MW_* \le \\
C\bigl(1+ \vartheta (\varepsilon) \mu h^{-1} + \bar{\vartheta}(\varepsilon)\mu^2\bigr)|\log \mu|
\label{18-5-44}
\end{equation}
(iv) Finally under assumption \textup{(\ref{18-4-35})} estimate
\begin{equation}
|\tilde{e}(x,x,0 ) - h^{-2} \cN_{2,x}^\MW (x,0) |\le C_s \mu ^{-s}
\qquad \text{for\ \ } x \in B(0,{\frac{1}{2}})
\label{18-5-45}
\end{equation}
holds.
\end{theorem}

We leave to the reader an easy

\begin{problem}\label{problem-18-5-8} 
Investigate conditions to the regularity (i.e. to $(l,\sigma)$ and $(\bar{l},\bar{\sigma})$) depending on $\mu$, $h$ and reversely conditions to $\mu$, $h$ depending on the regularity required to have the same results as in the smooth case. 
\end{problem}

\section{Generalizations}
\label{book_new-sect-18-5-3}

\subsection{Vanishing $V$}
\label{book_new-sect-18-5-3-1}

We want to get rid of assumption (\ref{18-2-27}) $V\le -\epsilon _0$. We need to deal only with the case  $\mu \le \epsilon _1 h^{-1}$. We should use the same rescaling arguments as in subsection~\ref{book_new-sect-13-7-1} but we need to take care of $\varepsilon$. 

Let us introduce admissible scaling functions
$\gamma (x) =\epsilon |V(x)|$ and $\rho (x)= \gamma (x)^{\frac{1}{2}}$. Then
rescaling ball $B\bigl(x,\gamma (x)\bigr)$ to $B(0,1)$ and dividing operator by
$\rho ^2$ we find ourselves again in our settings but with $h$, $\mu$ replaced
by $h_\new = h/\rho \gamma$, $\mu _\new =\mu \gamma / \rho $ respectively.
However we need a cut-off: in fact $\gamma (x) =\epsilon |V(x)|+ \bar{\gamma}_0$
with an appropriate parameter $\bar{\gamma }_0$. Modifications in some cases
are needed.

Consider first the case of theorem \ref{thm-18-5-4} when
$\bar{\mu}_1=(h |\log h|)^{- \frac{1}{3}}\le \mu \le h^{-1}$.
Then as $\mu_\new h_\new = \mu h /\rho^2=\mu h/\gamma$ and since we do not need assumption (\ref{18-2-27}) as $\mu_\new h_\new \gtrsim 1$ we will cut off $\gamma$ by picking up $\bar{\gamma}_0= \mu h $. 

Plugging $\mu_\new$, $\gamma_\new$ and corresponding mollification parameter 
\begin{equation*}
\varepsilon_\new  = (\mu_\new ^{-1}h_\new |\log h_\new|)^{\frac{1}{2}}= (\mu^{-1}h)^{\frac{1}{2}} \gamma^{-1}|\log (h/\gamma^{\frac{3}{2}})|
\end{equation*}
to the remainder estimate (\ref{18-5-39}) we see that the contribution of each partition element to the remainder does not exceed
\begin{equation*}
C\mu ^{-1}h ^{-1}\gamma + C(\mu ^{-1}h)^{\frac{l}{2}}
|\log (h/\gamma ^{\frac{3}{2}})|^{{\frac{l} {2}}-\sigma}
h^{-2}\gamma ^{3-l}
\end{equation*}
and the total remainder estimate will be then due to assumption (\ref{18-2-15})
\begin{equation*}
C\int \mu ^{-1}h \gamma ^{-1}\,d\gamma +
C(\mu ^{-1}h)^{\frac{l}{2}}
\int |\log (h/\gamma ^{\frac{3}{2}})|^{{\frac{l} {2}}-\sigma}
h^{-2}\gamma ^{3-l}\,d\gamma.
\end{equation*}
The second integral is equal to the second term in the right-hand expression of
(\ref{18-5-39}) for $l<2$; to handle the case $l=2$ we need to notice that the second term in (\ref{18-5-39}) is due to mollification error only but in the rescaled coordinates mollification parameter is $\varepsilon_\new$ which in the original coordinates becomes $\varepsilon_\new \gamma = 
(\mu^{-1}h)^{\frac{1}{2}} |\log (h/\gamma^{\frac{3}{2}})|\le (\mu ^{-1}h |\log h|)^{\frac{1}{2}}$ and thus mollification error does not exceed the second term in (\ref{18-5-39}) anyway. Therefore

\begin{claim}\label{18-5-46}
Theorem \ref{thm-18-5-4}(i) without condition (\ref{18-2-27}) remains true with the remainder estimate $C\mu ^{-1}h^{-1}|\log (\mu h)|+R_2$ where $R_2$ is the second term in the right-hand expression in (\ref{18-5-39}).
\end{claim}

The analysis of the case of theorem \ref{thm-18-2-9} is even simpler: one need to pick $\bar{\gamma}_0=\mu ^{-2}$; then on each element after rescaling
\underbar{either} condition (\ref{18-2-27}) is fulfilled \underbar{or}
$\gamma\asymp \bar{\gamma}_0$ and $\mu_\new \asymp 1$ and condition (\ref{18-2-27}) is not needed while $h_\new \asymp \mu ^3 h\le 1$.

Note that in this case 
$\varepsilon_\new = \mu_\new h_\new |\log h|=
\mu h \gamma^{-1} |\log (h/\gamma^{\frac{3}{2}})|$ and after rescaling to original coordinates it is $\varepsilon _\new \gamma \le \mu h|\log h|$.

Then one can prove easily that

\begin{claim}\label{18-5-47}
Theorem \ref{thm-18-2-9}(i) without condition (\ref{18-2-27}) remains true with the remainder estimate $C\mu ^{-1}h^{-1}|\log (\mu h)|+R_1$ where $R_1$ is the second term in the right-hand expression in (\ref{18-2-35}).
\end{claim}

Actually one can get rid of this logarithmic factor considering
propagation of singularities which is going in direction orthogonal to
$\nabla V$ with the speed $\asymp \mu ^{-1}$. Details in the smooth case can be found in subsection~\ref{book_new-sect-13-7-1}.

Thus we can replace $T^* =\epsilon \mu $ by 
$T^* = \epsilon\mu\min\bigl(|\log \gamma |,|\log (\gamma/\bar{\gamma}_0|\bigr)^\alpha $
in the Tauberian theorem provided $(1,\alpha )\preceq (l,\sigma)$.
Also, one can calculate ``the second term'' and to prove that it vanishes.

After summation we get an extra term in the remainder estimate
\begin{equation*}
C\mu ^{-1}h^{-1} \int \gamma ^{-1} \Bigl(|\log \gamma |^{-\sigma }+
|\log (\gamma/\bar{\gamma}_0|^{-\sigma}\Bigr)\,dx \asymp C\mu ^{-1}h^{-1}.
\end{equation*}
Therefore

\begin{claim}\label{18-5-48}
Theorems \ref{thm-18-5-4}(i), \ref{thm-18-2-9}(i) remain true without assumption (\ref{18-2-27}).
\end{claim}

Combining with the other statements of theorems \ref{thm-18-5-4}, \ref{thm-18-2-9} we conclude that these theorems remain true under assumption (\ref{18-2-27}) replaced by
\begin{equation}
|V|+|\nabla V|\ge \epsilon_0.
\label{18-5-49}
\end{equation}

This conclusion enables us to run rescaling again but now with the scaling function
$\gamma=\epsilon \bigl(|V|+|\nabla V|^2\bigr)^{\frac{1}{2}}$ and 
$\rho \asymp \gamma$; then $\mu_\new^{-1}h_\new^{-1} \asymp \mu^{-1}h^{-1}\gamma^2$. Then in the framework of theorem~\ref{thm-18-5-4} $\varepsilon_\new =
\bigl(\mu^{-1} h |\log (h/\gamma^2)| \bigr)^{\frac{1}{2}}\gamma^{-1}$, which in the original coordinates is 
$\varepsilon_\new\gamma \lesssim \bigl(\mu^{-1} h |\log h| \bigr)^{\frac{1}{2}}$ again. 

Meanwhile in the framework of theorem~\ref{thm-18-2-9} 
$\varepsilon_\new = \mu h |\log (h/\gamma^2)| \gamma^{-2}$ and a mollification error term in the virtue of assumption (\ref{18-2-16}) is
\begin{equation*}
C h^{-2}\int \gamma ^4 \bigl(\mu h |\log (h/\gamma^2)|\bigr)^l \gamma^{-2l} 
|\log \bigl(\mu h |\log (h/\gamma^2)|\bigr)|^{-\sigma}\, \gamma^{-1}d\gamma
\end{equation*}
which as one can see easily resets to the value of such term in theorem~\ref{thm-18-2-9}. Therefore

\begin{theorem}\label{thm-18-4-9}
Theorems \ref{thm-18-5-4}, \ref{thm-18-2-9} without assumption \textup{(\ref{18-2-27})}. 
\end{theorem}

\subsection{Other generalizations}
\label{book_new-sect-18-4-5-2}

We leave to the reader an easy parts (i) of the problems~\ref{problem-18-5-10} and \ref{problem-18-5-11} below; on the contrary, parts (ii) look extremely challenging.

\begin{Problem}\label{problem-18-5-10}
To get rid of condition (\ref{18-1-6}) (while (\ref{18-2-27}) holds) under assumption 
\begin{equation}
|\nabla (F/V)|\ge \epsilon.
\label{18-5-50}
\end{equation}
(i) Use rescaling method and recover remainder estimate $O(h^{-1})$ as $\mu \lesssim h^{-2}$.

\medskip\noindent
(ii) Generalize\emph{ Chapter~\ref{book_new-sect-14}. \nameref{book_new-sect-14}} to the non-smooth case.
\end{Problem}

\begin{Problem}\label{problem-18-5-11}
To get rid of condition $B(0,1)\subset X$.

\medskip\noindent
(i) Use rescaling method and recover remainder estimate $O(h^{-1})$ as $\mu \lesssim h^{-2}$.

\medskip\noindent
(ii) Generalize \emph{Chapter~\ref{book_new-sect-15}. \nameref{book_new-sect-15}} to the non-smooth case.
\end{Problem}

Anther challenging problem:

\begin{Problem}\label{problem-18-5-12}
Generalize \emph{Chapter~\ref{book_new-sect-17}. \nameref{book_new-sect-17}} to the non-smooth case.
\end{Problem}

\begin{remark}\label{rem-18-5-13}
One could think about generalizing \emph{Chapter~\ref{book_new-sect-16}. \nameref{book_new-sect-16}}  to the non-smooth case but as we do not have monotonicity of $e(x,x,\tau)$ with respect to operator it does not make much sense as we need to deal with framing approximations.. 

Sure we do not have such monotonicity of 
$\int e(x,x,\tau)\psi(x)\,dx$ either but this kind of expressions are building blocks for $\N(\tau)=\int e(x,x,\tau)\,dx$.
\end{remark}

\chapter{$d=3$: Weak magnetic field}
\label{book_new-sect-18-6}

In this section we consider the case
\begin{equation}
\mu \le \bar{\mu}_1\Def C(h |\log h|)^{-\frac{1}{3}}.
\label{18-6-1}
\end{equation}

In the smooth case then one needs to reduce operator to canonical form only to
analyze relatively long (comparing with $\epsilon _0\mu ^{-1}$) propagation of
singularities but not to calculate asymptotics. The same situation preserves in
non-smooth case and to recover sharp asymptotics one needs neither non-degeneracy nor extra smoothness. This is the same threshold  (\ref{18-2-6})  as in section~\ref{book_new-sect-18-2} for $d=2$ despite completely different mechanism; in some sense this is just coincidence.

\section{Preliminary remarks}
\label{book_new-sect-18-6-1}

We consider microlocally mollified Schr\"odinger operator in the strong magnetic
field given by expression (\ref{18-1-1}) with symmetric positive matrix 
$(g ^{jk})$, real-valued $V_j$, $V$ and small parameters $h$, $\varepsilon $ and large parameter $\mu $. We assume that the corresponding magnetic field $\mathbf{F}=(F^1,F^2,F^3)$ with 
\begin{equation*}
F ^l(x)\Def \frac{1}{2}\sum _{jk}\upepsilon ^{lkj}F_{kj}, \qquad
F_{kj}\Def (\partial _{x_j}V_k-\partial _{x_k}V_j)
\end{equation*}
does not vanish in $B(0,1)$; then we can straighten it and direct along
$x_3$ by an appropriate change of variables and make $V_3=0$ by an
appropriate gauge transformation. We assume that this has been done and impose
smoothness conditions to the reduced operator.

Namely we assume that for this reduced operator conditions (\ref{18-1-3}), (\ref{18-1-6}) hold and
\begin{equation}
g ^{jk}, V, F^3\in \sF ^{1,\sigma }(B(0,1)),\qquad F^1=F^2=0.
\label{18-6-2}
\end{equation}

We start from the case (further restrictions to follow)
\begin{equation}
\varepsilon \ge Ch |\log h|,\quad \mu \le C^{-1}h |\log h|^{-1}.
\label{18-6-3}
\end{equation}
Then\begin{phantomequation}\label{18-6-4}\end{phantomequation}
\begin{align}
& h ^{-1}[A,x_j]=\sum_k \bigl(g ^{jk}P_k+P_kg ^{jk}\bigr)
\tag*{$\textup{(\ref*{18-6-4})}_1$}\label{18-6-4-1}\\
& h ^{-1}[A,P_3]= h ^{-1}\sum_{j,k}P_j[g ^{jk},P_3] P_k+h^{-1}[V,P_3] \tag*{$\textup{(\ref*{18-6-4})}_2$}\label{18-6-4-2}\\
&\mu ^{-1}h ^{-1}[A,P_j]=
i (-1) ^{j+1}F_{12}\sum_k (g ^{53-j,k}P_k+P_k g ^{3-j,k})+
\tag*{$\textup{(\ref*{18-6-4})}_3$}\label{18-6-4-3}\\
&\qquad\qquad \mu ^{-1}h ^{-1}\sum_{jk}P_j[g ^{jk},P_i] P_k
+\mu ^{-1}h ^{-1}[V,P_i]\quad\text{for\ \ } i=1,2\notag
\end{align}
have symbols of class $\sF ^{1,\sigma }$ in any domain 
$\{(x,\xi):\ a(x,\xi )\le c\mu^{-2}\}$ of bounded energy where
\begin{multline}
a(x,\xi )= a_0(x,\xi )+\mu^{-2}V, \quad
a_0(x,\xi)=\sum_{j,k}g ^{jk}p_j(x,\xi )p_k(x,\xi ),\\
p_j(x,\xi )\Def(\xi _j- V_j).
\label{18-6-5}
\end{multline}
To finish the preliminary remarks let us notice that the following statement holds

\begin{proposition} \label{prop-18-6-1} 
Let $f\in \sC_0 ^\infty(\bR ^3)$, $f=1$ in $B(0,C_0)$ with large enough constant $C_0$. Let $T\ge h ^{1-\delta }$.
Then under condition \textup{(\ref{18-6-2})} with large enough constant $C=C_s$
\begin{multline}
|F_{t\to h ^{-1}\tau }
\Bigl(1-f\bigl(p_1(x,\xi ),p_2(x,\xi ), p_3(x,\xi )\bigr)\Bigr)^\w
\chi _T(t) U |\le Ch ^s\\
x,y\in B(0,{\frac{1}{2}}), \tau \le c.
\label{18-6-6}
\end{multline}
\end{proposition}

\begin{proof}
One can see easily that symbol 
$f\bigl(p_1(x,\xi ),p_2(x,\xi ), p_3(x,\xi )\bigr)$ is quantizable under condition (\ref{18-6-3}) and the proof is rather obvious. We leave it to the reader.
\end{proof}

\section{Heuristics}
\label{book_new-sect-18-6-2}

Let us consider operator (\ref{18-1-1}) with the Euclidean metrics. Then for time $T(\xi_3)\Def \epsilon _0|\xi_3|$ the shift with respect to $\xi_3 $ will be less than ${\frac{1}{2}}|\xi_3|$ and so $\xi_3$ will remain of the same magnitude and sign. Then the shift with respect to $x_3$ will be of magnitude $|\xi_3| T$ because $\xi_3$ is a speed of the propagation along $x_3$ and in order it to be observable we need logarithmic uncertainty principle
\begin{equation}
|\xi_3| T \cdot |\xi_3|\ge Ch|\log h|
\label{18-6-7}
\end{equation}
because we need to use $|\xi_3|$-sized scale along $\xi_3$. This gives us the
value $T_*$ for which we can claim that $(\Gamma_x U)(t)$ is negligible for
$ |t| \in [T_*,T^*]$: 
\begin{equation}
T_*=T_*(\xi_3)\Def Ch|\log h| \cdot |\xi_3|^{-2}.
\label{18-6-8}
\end{equation}
Since we know
$(\Gamma_x U)(t)$ for $|t|\le \bar{T}\Def \epsilon \mu^{-1}$ due to condition
(\ref{18-2-27}) and rescaling, we would like to have these intervals overlap:
$\bar{T}\ge T_*$ or, equivalently,
\begin{equation}
|\xi_3|\ge \bar{\rho}_1 \Def C(\mu h |\log h|)^{\frac{1}{2}}
\label{18-6-9}
\end{equation}
exactly as in the smooth case in section~\ref{book_new-sect-13-3}.

Then the contribution of the zone
$\cZ_1^c \Def \{(x,\xi):\  |\xi_3|\ge \bar{\rho}_1\}$ to
the remainder does not exceed 
\begin{equation}
C h^{-2} \int T^{*\,-1}(\xi_3)\, d\xi_3.
\label{18-6-10}
\end{equation}
This integral diverges logarithmically at $\xi_3=0$ and is equal to
$Ch^{-2}|\log \bar{\rho}_1|$; to improve remainder estimate to $Ch^{-2}$ we
will need to increase $T^*(\xi_3)$ to
$T^*(\xi_3)\Def \epsilon |\xi_3|\cdot \bigl|\log |\xi_3|\bigr|^{\sigma'}$ with $\sigma'>1$.

On the other hand, the contribution of the zone
$\cZ_1 \Def\{ (x,\xi):\ |\xi_3|\le \bar{\rho}_1\}$ to the remainder does not exceed $Ch^{-2}\bar{T}^{-1}\bar{\rho}_1= C\mu \bar{\rho}_1 h^{-2}$ and in order to keep it below $Ch^{-2}$ we need to have $\mu \bar{\rho} _1\le C$ which is
exactly condition (\ref{18-6-1}). 

But in this case zone $\cZ_2\Def \{(x,\xi):\ |\xi_3|\le \bar{\rho}_2=C\mu^{-1}\}$ is covered by these latter arguments; therefore we can take $C\mu^{-1}$ as a critical value and replace (\ref{18-6-9}) by a stronger restriction
\begin{equation}
|\xi_3|\ge \bar{\rho}_2  \Def C\mu^{-1}.
\label{18-6-11}
\end{equation}
Condition (\ref{18-6-9}) will be restored in sections~\ref{book_new-sect-18-7}--\ref{book_new-sect-18-9} where we analyze the case
$\mu \ge C(h|\log h|)^{-\frac{1}{3}}$ and restriction (\ref{18-6-9}) is stronger.

This analysis is basically ``smoothness-independent''. It allows us to apply weak magnetic field approach without any non-degeneracy assumptions.

\section{Variable $\varepsilon $}
\label{book_new-sect-18-6-3}

So far we did not discuss the choice of $\varepsilon$; however, since the scale
with respect to $\xi_3$ is $C|\xi_3|$ in the zone
$\cZ_2^c\Def \{(x,\xi):\ |\xi_3|\ge \bar{\rho}_2\}$ the logarithmic uncertainty
principle requires $\varepsilon \ge Ch|\log h|\cdot |\xi_3|^{-1}$. To reduce
approximation error we pick up the smallest allowed $\varepsilon$ there:
\begin{equation}
\varepsilon = Ch|\log h|\cdot |\xi_3|^{-1}.
\label{18-6-12}
\end{equation}
Then 
\begin{claim}\label{18-6-13}
The contribution of the zone $\cZ_2^c$ to the approximation error does not exceed 
$Ch^{-3}\int \vartheta \bigl(\varepsilon (\xi_3)\bigr) d\xi_3$ which is $O(h^{-2})$ as $(l,\sigma)\succeq (1,2)$.
\end{claim}
On the other hand, in the zone $\cZ_2\Def \{(x,\xi):\ |\xi_3|\le \bar{\rho}_2\}$ the scale with respect to $\xi_3$ is $\asymp \bar{\rho}_2$ and the logarithmic uncertainty principle requires there
$\varepsilon = \bar{\varepsilon}_2$ with
\begin{equation}
\bar{\varepsilon}_2 \Def Ch|\log h| \cdot \bar{\rho}_2^{-1}=C \mu h|\log h|. 
\label{18-6-14}
\end{equation}
Then the contribution of this zone to the approximation error will be
$Ch^{-3}\bar{\rho}_2 \vartheta(\bar{\varepsilon}_2)$ which is $O(h^{-2})$ as
$(l,\sigma)\succeq (1,1)$.

However, if we tried to use $\varepsilon =\bar{\varepsilon}_2$ everywhere we
would get an approximation error $Ch^{-3}\vartheta (\bar{\varepsilon}_2)$ which
is not $O(h^{-2})$ unless stronger smoothness condition is satisfied: namely, to cover $\mu \asymp (h|\log h|)^{\frac{1}{3}}$ we would need $(l,\sigma)\succeq (\frac{3}{2},1)$.

Thus, for $d=3$ we will need $\varepsilon $ to be a function of $\xi _3$ and may
be of $x$. We already met such situation in \ref{book_new-sect-4-5}. We will pick $\varepsilon $ to
be a temperate function:
\begin{equation}
\rho |\nabla _\xi \varepsilon | + \gamma |\nabla _x \varepsilon |\le
c\varepsilon
\label{18-6-15}
\end{equation}
where $(\gamma,\rho)$ is $(x,\xi)$-scale. In this case we introduce
\begin{equation}
\tilde{A}= \sum \psi _\alpha ^* \tilde{A}_\alpha \psi _\alpha
\label{18-6-16}
\end{equation}
where $\psi _\alpha $ are $(\gamma , \rho )$-admissible functions,
$\psi _\alpha ^*\psi _\alpha $ is a partition of unity and $\tilde{A}_\alpha $
are mollifications of $A$ with $\varepsilon = \varepsilon_\alpha $ calculated at
some points of the corresponding elements.

\section{Rigorous analysis}
\label{book_new-sect-18-6-4}

As we already mentioned we consider now the case of the weak magnetic field:
$\mu \le \bar{\mu}_1 \Def (h | \log h|) ^{-\frac{1}{3}}$ assuming that
\begin{equation}
l=1,\quad \sigma >1
\label{18-6-17}
\end{equation}
and we will not press very hard to reduce $\sigma $. Let us consider
point $(\bar{x},\bar{\xi} )$ with
\begin{equation}
\rho \Def |\bar{\xi} _3| \ge \bar{\rho}_2\Def C\mu ^{-1},
\qquad \varepsilon = Ch\rho ^{-1}|\log h|.
\label{18-6-18}
\end{equation}
Recall that the second condition (to $\varepsilon $) means the logarithmic
uncertainty principle; thus quantization is possible; we will take 
$\varepsilon $ variable but depending on $\xi_3 $ only.

Let us consider first a classical propagation starting from this point. Obviously $|d\xi_3/dt |\le C_0$ due to condition \ref{18-6-4-3} which implies
that $\xi_3$ keeps it magnitude during time $|t|\le \epsilon_0\rho$ even in the general (not Euclidean) case but we need a bit better. Let us consider
\begin{equation*}
\frac{d\ }{dt} \xi_3= 
\{\mu ^2a(x,\xi), \xi_3\} = \beta (x,\xi) \xi_3 +
\mu \sum _{j,k=1}^2\beta _{jk}p_jp_k - \partial _{x_3}V
\end{equation*}
with $\beta _*\in \sF^{0,\sigma}$. In the smooth case we would be able to reduce the second term in the right-hand expression to the form $\beta'' a_0(x,\xi)$ by replacing $\xi_3$ by
\begin{equation}
\ell =\xi _3 + \sum _{j,k=1}^2\beta ' _{jk}p_jp_k
\label{18-6-19}
\end{equation}
with appropriate coefficients $\beta '_*$ found from the linear algebraic
system which we would get from equality
\begin{multline}
\{\mu ^2a(x,\xi), \ell (x,\xi)\} \equiv \beta ' (x,\xi) \xi_3 +
\mu \beta '' a_0(x,\xi)- \partial _{x_3}V \\
\mod O\bigl(\mu ^{-1}\bigr)
\label{18-6-20}
\end{multline}
(since $|p_j|\le C\mu ^{-1}$). However, in on our non-smooth case we replace 
$\beta ' _{jk}$ by their $h ^{-\delta}$-mollifications with the small enough exponent $\delta>0$; then we would get the above equations (\ref{18-6-20}) modulo $O(|\log \rho |^{-\sigma})$. Finally, one can find that
$\mu ^2 \beta '' = \partial_3 F$
and therefore
\begin{multline}
\{\mu ^2a(x,\xi), \ell (x,\xi)\} \equiv \beta ' (x,\xi) \ell (x,\xi) +
\mu ^2 F a (x,\xi)- F\partial _{x_3}(V/F) \\[2pt]
\mod O(|\log \varrho |^{-\sigma}).
\label{18-6-21}
\end{multline}
Then in an appropriate time direction\footnote{\label{foot-18-29} That means for $t$ of an appropriate sign; namely the sign of $-\xi_3 \partial_{x_3}(V/F)$.}
inequality $\pm \ell \ge \frac{1}{2}\rho $ holds for $|t|\in [0,T^*]$ with
\begin{equation}
T^*=T^*(\rho ) \Def \epsilon \rho |\log \rho | ^ \sigma .
\label{18-6-22}
\end{equation}
This additional logarithmic factor is a small but a crucial progress in comparison with $T^*(\rho) =\epsilon_0\varrho$ we had before.

Then using the following function with an appropriate sign $\varsigma$ in the
standard (see Section~\ref{book_new-sect-2-3}) propagation arguments
\begin{equation}
\upchi\Bigl(\pm\frac{1}{\rho}\bigl((\ell (x,\xi )-
\ell(\bar{x},\bar{\xi} )\bigr) + c\varsigma \frac{t} {T}\Bigr)
\label{18-6-23}
\end{equation}
one can prove the similar statement on microlocal level (due to conditions
(\ref{18-6-2}), $\textup{(\ref{18-6-4})}_{1-3}$, (\ref{18-6-17})--(\ref{18-6-18}) all symbols will be quantizable):

\begin{proposition}\label{prop-18-6-2} 
Let conditions \textup{(\ref{18-6-2})}, \textup{(\ref{18-6-3})}, \textup{(\ref{18-6-17})}--\textup{(\ref{18-6-18})} be fulfilled. 

Let $Q$ be $h$-pseudo-differential operator with the symbol supported in 
$\rho $-vicinity of $(\bar{x},\bar{\xi})$ with 
$|\bar{\xi}_3|=\rho \ge C\mu ^{-1}$ and let $Q'$ be $h$-pseudo-differential operator with the symbol equal to $1$ in $\rho'$-vicinity of $(\bar{x},\bar{\xi})$ (with  $\rho '=C\rho |\log \rho |^\sigma$), intersected with $\{|\xi_3 | \ge \frac{1}{2}\rho\}$.

Then $(I-Q')e^{-ih^{-1}At} Q$ is negligible for $\varsigma t\in [0,T^*]$ with $T^*=T^*(\rho)$ defined by \textup{(\ref{18-6-22})} and appropriate $\varsigma=\pm 1$.
\end{proposition}

\begin{proof}
The standard proof is left to the reader.
\end{proof}

This statement shows that for $\varsigma t\in [0,T^*]$ singularity of
$e^{-ih^{-1}At}Q$ is confined to $\rho'$-vicinity of $(\bar{x} ,\bar{\xi})$ intersected with $\{|\xi_3 | \ge \frac{1}{2} \rho\}$.

Now let us consider classical propagation along $x_3$; then
\begin{equation*}
dx_3/dt =\{ \mu^2a, x_3\}= 2g^{33}\xi_3 + \mu \sum_{j=1}^2 2g^{3j}p_j.
\end{equation*}
Replacing $x_3$ by
\begin{equation}
\phi = x_3+\sum _{j=1}^2 \alpha_j p_j
\label{18-6-24}
\end{equation}
with appropriate coefficients $\alpha_j$ found from an appropriate linear
algebraic system we would able in the smooth case get
\begin{equation}
\{ \mu^2a, \phi\}\equiv 2g^{33}\xi_3\quad \mod O(\mu ^{-1}).
\label{18-6-25}
\end{equation}
In the non-smooth case we do the same and no mollification is needed since
$\alpha_j\in \sF^{1,\sigma}$. So, the classical shift with respect to
$\phi $ for time $t$ with $\varsigma t\in [0,T^*]$ will be at least
$\epsilon_0\rho t$ which is ``observable'' if
$\epsilon_0\rho |t |\times \rho \ge Ch|\log h|$.

On the microlocal level using the same arguments of Section~\ref{book_new-sect-2-3} with function
\begin{equation}
\upchi \Bigl(\varsigma_1 \frac{1} {\rho T}
\bigl(\phi (x,\xi)-\phi (y,\xi)\bigr)\pm c\frac{t} {T}\Bigr)
\label{18-6-26}
\end{equation}
(which again is quantizable due to our conditions) with an appropriate sign
$\varsigma_1=\pm1$ one can easily prove that for $\varsigma t\in [T_*,T^*]$,
\begin{equation}
T_*= C\rho ^{-2}h|\log h|
\label{18-6-27}
\end{equation}
and $\psi (x_3)$, supported in $\epsilon $-vicinity of $\bar{x}_3$,
$\psi (x_3) e^{-ih^{-1}At} Q\psi (x_3)\equiv 0$.

Then $\Tr \bigl(e^{-ih^{-1}At} Q \bigr)\equiv 0$.
Furthermore this is true for $|t|\in [T_*,T^*]$ \emph{with no restriction to the sign because of the trace\/}. Therefore we arrive

\begin{proposition} \label{prop-18-6-3} 
Let the above assumptions be fulfilled. Then for any $T\in [T_*,T^*]$
\begin{equation}
F_{t\to h^{-1}\tau } \Bigl( \bigl(\bar{\chi}_{T^*}(t) -
\bar{\chi}_T(t)\bigr) \Gamma (U \,^t\!Q_y) \Bigr) \equiv 0.
\label{18-6-28}
\end{equation}
\end{proposition}

\begin{proof}
Again, the standard details are left to the reader.
\end{proof}

Recall that $\chi $, $\bar{\chi}$ are admissible functions, supported in
$[-1,1]$ and equal to $0,1$ respectively on $[-\frac{1}{2}, \frac{1}{2} ]$.

Then the same statement holds for $Q=\psi (x) \chi_\rho (hD_3)$ (we can sum
with respect to the partition of unity in $(x,\xi_1,\xi_2)$). \emph{From now on
operator $Q$ is given by this expression\/}.

On the other hand, from the standard theory rescaled we already know 
$F_{t\to h^{-1}\tau } \bigl( \bar{\chi}_{\bar{T}}(t)\Gamma (u \,^t\!Q_y) \bigr)$
with $\bar{T}=\epsilon_0\mu ^{-1}$ and in order to combine these two results we
need to have $T_*\ge \bar{T}$ which is equivalent to (\ref{18-6-9});
this condition is wider than our current framework; right now we have
$(\mu h |\log h|)^{\frac{1}{2}} \lesssim \mu^{-1}$ which exactly means that
$\mu \lesssim \bar{\mu}_1 = (h|\log h|)^{-\frac{1}{3}}$.

Then from the standard theory rescaled we conclude that
\begin{equation}
|F_{t\to h^{-1}\tau } \bigl( \bar{\chi}_T(t)\Gamma (U\,^t\!Q_y) \bigr) |
\le C\rho h^{-2}
\label{18-6-29}
\end{equation}
for any $T\in [T_*,T^*]$. With this estimate the standard Tauberian arguments
imply that
\begin{multline}
\R^\T _Q\Def
|\Gamma (\tilde{e}\,^t\!Q_y)(\tau) - h^{-1}\int_{-\infty}^\tau
\Bigl(F_{t\to h^{-1}\tau'} \bigl(\bar{\chi}_T(t)\Gamma (U \,^t\!Q_y) \bigr)\Bigr) \, d\tau'|\le \\[3pt]
C\rho h^{-2}T ^{*\,-1}= C h^{-2}|\log \rho |^{-\sigma}.
\label{18-6-30}
\end{multline}

Let us sum with respect to partition on $\xi _3$ in the zone
$\cZ_2^c = \{(x,\xi):\ |\xi_3|\ge \bar{\rho}_2\}$. Then the right-hand expression of this estimate transforms into 
\begin{equation*}
Ch^{-2}\int_{\rho^* _2}^1 \rho ^{-1} |\log \rho |^{-\sigma}\,d\rho 
\end{equation*}
which does not exceed $Ch^{-2}$ as long as $\sigma >1$.

On the other hand, it follows from the standard theory rescaled that estimate (\ref{18-6-29}) holds for $Q_\rho =\psi (x) \bar{\chi}_\rho (hD_3)$ and
$T=\bar{T}= \epsilon \mu ^{-1}$. Therefore for $\rho = \bar{\rho}_2$ and $\bar{Q}=Q_{\bar{\rho}_2}$ we
have 
\begin{equation*}
\R^\T _{\bar{Q}}\le C\bar{\rho}_2 h^{-2}/\bar{T} = Ch^{-2}.
\end{equation*}

Combining this estimate with the result of the previous analysis we get the
final inequality 
\begin{equation}
\R^\T _\psi \le Ch^{-2}.
\label{18-6-31}
\end{equation}
Further, from the standard theory rescaled we know
$F_{t\to h^{-1}\tau } \bigl( \bar{\chi}_T(t)\Gamma (U \,^t\!Q_y) \bigr)$
for $T={\bar T}$ and we get the final estimate\footnote{\label{foot-18-30} There should be also a term $C\mu ^2h^{-1}$ in the right-hand expression but
for $\mu \le h^{-\frac{1}{2}}$ it is less than $Ch^{-2}$. We will discuss this
expression with more details later.}
\begin{gather}
|\int \Bigl( \tilde{e}(x,x,\tau) - h^{-3} \tilde{\cN}_{3,x} ^\W (\tau)\Bigr)
\psi (x)\, dx|\le Ch^{-2}
\label{18-6-32}
\shortintertext{where}
h^{-3}\tilde{\cN}^\W_{3,x}(\tau) =
(2\pi )^{-3} \int \uptheta \bigl(\tau - \tilde{A}(x,\xi)\bigr)\, d\xi
\label{18-6-33}
\end{gather}
and $\tilde{A}$ is a symbol of the mollified operator and this is $h$-pseudo-differential operator because we take $\varepsilon$ depending on $\xi_3$\,\footnote{\label{foot-18-31} Which means that we take appropriate partition of unity, on each element take
its own $\varepsilon$, mollify and then add. One can check easily that
perturbation is small enough to preserve all our arguments.}.

Now let us consider an approximation error in $h^{-3}\tilde{\cN}^\W$. It does
not exceed
\begin{equation}
C h^{-3}\int \varepsilon(\rho)^l |\log\varepsilon (\rho) | ^{-\sigma }\,d\rho 
\asymp Ch^{-2}|\log h|^{1-\sigma}\log \mu
\label{18-6-34}
\end{equation}
which is $O(h^{-2})$ as $\sigma=2$\,\footnote{\label{foot-18-32} Exactly in this place we need to have variable $\varepsilon =\varepsilon(\rho) $ rather than to
take its largest value $\varepsilon (\bar{\rho}_2)=C \mu h |\log h|$ which
would require too restrictive smoothness conditions to get approximation error
$O(h ^{-2})$.}.

This completes the proof of 

\begin{theorem}\label{thm-18-6-4}
Let $d=3$ and $A$ be a self-adjoint in $\sL^2(X)$ Schr\"odinger operator \textup{(\ref{18-1-1})}. Let conditions \textup{(\ref{18-1-3})}, \textup{(\ref{18-1-6})}, \textup{(\ref{18-2-27})}
fulfilled in $B(0,1)\subset X\subset \bR ^3$. Let $(l,\sigma )\succeq (1,1)$.

Let $\mu \le \bar{\mu}_1= (h|\log h|)^{-\frac{1}{3}}$. Then for two framing approximations as in footnote~\footref{foot-18-16} estimate
\begin{multline}
\R^\W\Def 
|\int \Bigl( \tilde{e}(x,x,0)-h^{-3}\cN_{3,x} ^\W (0) \Bigr)\psi (x)\,dx|\le \\
Ch ^{-2} + Ch^{-3}\vartheta (\mu h|\log h|)
\label{18-6-35}
\end{multline}
holds with the standard Weyl expression 
\begin{equation*}
h^{-3}\cN _{3,x} ^\W(\tau ) =
 \frac {1} {6\pi^2} h^{-3}\bigl(\tau -V(x)\bigr)_+^{\frac{3}{2}}\sqrt g.
\end{equation*}
In particular, $\R^\W =O(h^{-2})$ as $(l,\sigma)\succeq (1,2)$.
\end{theorem}

Recall, that this theorem does not require any non-degeneracy assumption and as it is it cannot be improved significantly in frames of the weak magnetic field approach without non-degeneracy assumptions even under improved smoothness assumptions. Later in Theorem~\ref{thm-18-6-10} we will  slightly improve this theorem with no additional assumptions (removing logarithmic factor)

In the rest of this chapter we are going to improve it under non-degeneracy assumption (\ref{18-2-15}) or (\ref{18-2-16}); it appears that extra smoothness will be required.  

\section{Outer zone}
\label{book_new-sect-18-6-5}
First of all, for $\mu \le c^{-1}(h|\log h|)^{-1}$ one can apply exactly the same arguments as before\footnote{\label{foot-18-33} And we need 
$(\bar{l},\bar{\sigma} )\succeq (l,\sigma )\succeq (1,2)$ only.} with
\begin{equation}
\rho \Def |\bar{\xi} _3| \ge \bar{\rho}_1 \Def
C\bigl(\mu h|\log h|\bigr) ^{\frac{1}{2}} \quad \text{and} \quad
\varepsilon = Ch\rho ^{-1}|\log h|
\label{18-6-36}
\end{equation}
and get an estimate (\ref{18-6-30}) for $Q$ with symbol supported in $\cZ_1^c$. Then to derive asymptotics of $\Gamma (\tilde{e}\,^t\!Q_y )$ (at this stage we include $\psi $ in $Q$) we need to calculate the Tauberian expression
\begin{equation}
h^{-1}\int_{-\infty}^\tau
\Bigl(F_{t\to h^{-1}\tau'}\bigl(\bar{\chi}_T(t)\Gamma (u \,^t\!Q_y) \bigr)\Bigr)
\,d\tau'
\label{18-6-37}
\end{equation}
with $T=\bar{T}= \epsilon _0\mu ^{-1}$. Note that

\begin{claim}\label{18-6-38} 
Under assumption (\ref{18-2-27}) expression (\ref{18-6-37}) with $T=\bar{T}$ equals modulo $O(h^s)$ to the same expression with $T= T_0\Def C h|\log h|$.
\end{claim}

Really, let us break $Q$ into few operators $Q_{(\iota )}$, with the symbols
supported in $\epsilon$-vicinity of points 
$(\bar{x}_{(\iota )},\bar{\xi}_{(\iota )})$. Without any loss of the generality one can assume that $g^{jk}=\updelta _{jk}$ at $\bar{x}_{(\iota)}$\,\footnote{\label{foot-18-34} We can achieve it by a linear change of coordinates preserving direction of the magnetic field $\mathbf{F}$.}.
Then if $|\bar{\xi}_{(\iota)j}| \ge \epsilon _0$ and $\epsilon>0$ is small enough then both the propagation speed and scale with respect to $x_j$ are
disjoint from $0$, $j=1,2,3$ and we can apply arguments of Section~\ref{book_new-sect-2-3}.

Now, to calculate expression (\ref{18-6-37}) with $T=T_0$ we can apply the method of successive approximations. However, it is not very convenient to take unperturbed operator in our standard way $\bar{A}=A(y,hD_x)$ because then one can estimate $(m+1)$-th term of (\ref{18-6-37}) with the approximation plugged only by
\begin{equation*}
h^{-3} |\log h| \bigl(\mu T_0 ^2 h^{-1} \bigr)^m \asymp
h^{-3} |\log h| \bigl(\mu h |\log h|^2\bigr)^m
\end{equation*}
which is not good as $\mu $ close to $(h|\log h|)^{-1}$ even if we reduce powers of logarithms (which is not very difficult).

Instead we take unperturbed operator (\ref{18-1-27})--(\ref{18-1-28}) exactly as in $2$-dimensional case.  Recall that without any loss of generality we assume $V_3=0$ and $V_j=V_j(x_1,x_2) \in \sF^{2,2}$.

With this choice of an unperturbed operator one can estimate $(m+1)$-th term of
(\ref{18-6-37}) with the approximation plugged by
\begin{equation*}
h^{-4} T_0 \bigl(\mu T_0 ^3 h^{-1} + T_0 ^2 h^{-1}\bigr)^m
\asymp h^{-3} |\log h| \cdot \bigl(h|\log h|^2\bigr)^m;
\end{equation*}
recall that $\mu \le C(h|\log h|)^{-1}$ and $T_0\asymp Ch|\log h|$. 

Using rescaling arguments one can prove rather easily that we can take effectively $T_0\asymp h$ and thus to reduce the power of the logarithmic factors but it is not worth of our efforts now; we just notice that modulo $O(h^{-2})$ expression (\ref{18-6-37}) is equal to the same expression with $U$ replaced by $U^0+u^1$ where $U^0$ and $U^1$ are first two terms of approximation:
\begin{align}
&U^0= -ih\sum_{\varsigma =\pm} \varsigma \bar{\cG}^\varsigma 
\updelta (x-y)\updelta (t),\label{18-6-39}\\
&U^1= -ih\sum_{\varsigma =\pm} \varsigma \bar{\cG}^\varsigma \cR 
\bar{\cG}^\varsigma \updelta (x-y)\updelta (t)\label{18-6-40}
\end{align}
where $\bar{\cG}^\pm$ are forward and backward parametrices respectively for
$(hD_t-\bar{A})$ and $\cR = A-\bar{A}$.

Now we want to get rid of the contribution of $U^1$ and only extra logarithmic
factors in its estimate prevent us. Then we need just notice that due to the
rescaling method the corresponding term in (\ref{18-6-37}) with $\bar{\chi}_T$ replaced by $\chi _T$ with $T\ge T_0 \Def Ch^{-1}$ does not exceed
$CTh^{-4} \times T^2h^{-1} \times \bigl(T_0 /T \bigr)^s$
with an arbitrarily large exponent $s$. Then the sum with respect to
$t$-partition does not exceed $Ch^{-2}$. 

On the other hand, this term calculated with $\bar{\chi}_T$ and $T=T_0$ does not exceed $Ch^{-2}$ as well. Therefore the total contribution of $U^1$ to the final answer does not exceed $Ch^{-2}$.

Finally, let us calculate expression (\ref{18-6-37}) with $U=U^0$. Without any loss of the generality one can assume that $g^{jk}(y)=\updelta_{jk}$; we need later to get the answer in the invariant form. Then
\begin{multline}
U^0(x,y,t)=\\ h^{-1}\iint
e^{-ih^{-1}\bigl(-t (\xi_3^2 + \tau '+V(y))+(x_3-y_3)\xi _3\bigr)}
\, d\xi_3 \, d_{\tau ''} {\bar e}_y(x_1,x_2;y_1,y_2;\tau '')
\label{18-6-41}
\end{multline}
where ${\bar e}_y(x_1,x_2;z_1,z_2;\tau )$ is the Schwartz kernel of the spectral
projector of ${\bar P}_1^2+{\bar P}_2^2$ and ``$y$'' indicates the same point as
in (\ref{18-1-27}). Direct calculations of Chapter~\ref{book_new-sect-13} show that
\begin{equation}
\bar{e}_y(x_1,x_2;y_1,y_2;\tau '')=
\sum _{n\ge 0} \uptheta \bigl(\tau ''-(2n+1)\mu h F(y)\bigr) F(y)\mu h^{-1}
\label{18-6-42}
\end{equation}
and plugging into (\ref{18-6-41}) we arrive to 
\begin{multline}
\Gamma_y \bigl(U^0 \,^t\!Q_y\bigr) = \\
\int
\sum _{n\ge 0} e^{ih^{-1}t\bigl(\xi_3^2 + (2n+1)\mu h F(y) +V(y)\bigr)}
q(\xi_3)\, d\xi_3 \times F(y) \mu h^{-2}\label{18-6-43}
\end{multline}
and consequently to
\begin{multline}
\quad F_{t\to h^{-1}\tau'} \Bigl( \bar{\chi}_T(t)\Gamma (u \,^t\!Q_y) \Bigr)=\\
 \sum _{n\ge 0} T \int \widehat{\bar{\chi}}
\Bigl(Th^{-1}\bigl(\tau ' -(2n+1)F(y) \mu h -V(y)-\xi_3^2 \bigr)\Bigr)\times \\
q(\xi _3)\, d\xi_3 \times F(y) \mu h^{-2}\label{18-6-44}
\end{multline}
where ${\widehat{\bar{\chi}}}$ is a Fourier transform of $\bar{\chi}$ and modulo $O(h^{-2})$ expression (\ref{18-6-37}) is finally equal to 
\begin{gather}
h^{-3}\int \tilde{\cN}^\MW _Q (x,\tau)\psi (x)\,dx\label{18-6-45}\\
\shortintertext{with}
\tilde{\cN} ^\MW _Q \Def \sum _{n\ge 0} \int
\uptheta \Bigl( \tau - (2n+1)F \mu h -\tilde{V}-\xi_3^2 \Bigr) q(\xi_3)\,d\xi_3 \times
F\mu h.
\label{18-6-46}
\end{gather}

One can check easily that $\varepsilon |\log h|^{-2}\le \rho ^2$ for
$\rho \ge \bar{\rho}_1$ and therefore the approximation error of
$h^{-3}\tilde{\cN}^\MW _Q$ does not exceed
$C\int h^{-3}\varepsilon |\log h|^{-2} d\rho \le Ch^{-2}$ (recall that $(l,\sigma)\succeq (1,2)$). So we arrive to

\begin{proposition}\label{prop-18-6-5} 
Let $(l,\sigma)\succeq  (1,2)$ and
$\bar{\mu}_1\le \mu \le c^{-1}(h |\log h|)^{-1}$. Let $Q=q (hD_3)$ with
$q=1-\bar{\chi}_{\bar{\rho}_1} (\xi_3)$.

Then
\begin{equation}
|\int \Bigl(\Gamma_x (\tilde{e}\,^t\!Q_y) - h^{-3}\tilde{\cN}^\MW _Q (x,\tau) \Bigr)
\psi (x)\,dx |\le Ch^{-2}
\label{18-6-47}
\end{equation}
with $\tilde{\cN}^\MW _Q (x,\tau)$ defined by \textup{(\ref{18-6-46})}.
\end{proposition}

\section{Analysis under assumption (\ref{18-2-15})}
\label{book_new-sect-18-6-6}
Let us prove that under assumption (\ref{18-2-15}) $\Gamma (U\psi \,^t\!Q_y)$ is
negligible for $\bar{T}\le |t|\le \epsilon$ as
\begin{equation}
\varepsilon = C\mu h |\log h|\quad \text{in\ \ } \cZ_1
\label{18-6-48}
\end{equation}
and $Q=q(hD_3)$ with symbol supported in $\cZ_1$. Recall that 
$\cZ_1=\{(x,\xi):\ |\xi_3|\le \bar{\rho}_1\}$.

\begin{proposition} \label{prop-18-6-6} 
Let condition \textup{(\ref{18-2-15})} be fulfilled, $Q=q(hD_3)$ with
appropriate symbol $q$ supported in $\cZ_1$. Then
\begin{equation}
|F_{t\to h^{-1}\tau }\bigl(\chi_T\Gamma u\psi \,^t\!Q_y \bigr)|\le Ch^s\qquad
\forall \tau: |\tau |\le \epsilon
\label{18-6-49}
\end{equation}
as $T\in [T_*,T^*]$ where $\chi $ an admissible function supported in
$[-1,-{\frac{1}{2}}]\cup [{\frac{1}{2}},1]$ and here
\begin{equation}
T_* =C\varepsilon ^{-1} h|\log h|, \qquad T^*=\epsilon.
\label{18-6-50}
\end{equation}

In particular \textup{(\ref{18-6-49})} holds as $T\in [\bar{T}, \epsilon]$ with
$\bar{T}=\epsilon_0\mu^{-1}$ for $\varepsilon $ given by \textup{(\ref{18-6-48})}.
\end{proposition}

We prove this proposition below, it implies

\begin{corollary}\label{cor-18-6-7} Under conditions \textup{(\ref{18-2-15})}, \textup{(\ref{18-2-27})} and \textup{(\ref{18-6-48})}
\begin{gather}
|\int \Bigl(\tilde{e}(x,x,\tau ) -
\int_{-\infty}^\tau F_{t \to h^{-1}\tau'}
\bigl(\bar{\chi}_{\bar{T}} u(x,x,t) \bigr)\Bigr)\psi (x)\,dx |\le Ch^{-2}
\label{18-6-51}\\
\shortintertext{and}
|\int \Bigl(\tilde{e}(x,x,\tau ) -
h^{-3}\tilde{\cN}^\MW (x,\tau ) \Bigr)\psi (x)\,dx |\le Ch^{-2}.
\label{18-6-52}
\end{gather}
\end{corollary}

\begin{remark}\label{rem-18-6-8} 
(i) Note that under assumption (\ref{18-2-15}) and $\varepsilon$ given by (\ref{18-6-48}) the contribution of zone $\cZ_1$ to the approximation error does not exceed 
\begin{equation}
C(\mu h|\log h|)^{l+\frac{1}{2}} |\log h|^{-\sigma} h^{-3}
\label{18-6-53}
\end{equation}
which is in turn does not exceed $Ch^{-2}$ as long as
\begin{equation}
\mu \le \mu^* _{(l,\sigma )}\Def
Ch^{-(2l-1)/ (2l +1)}|\log h|^{-1+ 2\sigma/(2l+1)}
\label{18-6-54}
\end{equation}
while remainder estimate does not exceed $Ch^{-2}$ anyway.

\medskip\noindent
(ii) Note that
$\mu^* _{(1,2)}= h^{-{\frac{1}{3}}}|\log h|^{\frac{1}{3}}$ and
$\mu^*_{(l,\sigma )} \sim h^{-1+\delta }$ in the smooth case when $l$ is large.
\end{remark}

\begin{proof}[Proof of Proposition \ref{prop-18-6-6}] 
We will use the canonical form reduction of the next section~\ref{book_new-sect-18-7}.

Making $\epsilon _1$-partition we can assume that on it 
$|\nabla V-\ell |\le \epsilon $ at each point with some fixed vector $\ell $, $|\ell | \ge 4\epsilon $.

Let us consider first the case when
$|\ell _\perp |\ge \frac{1}{2} \epsilon $, $\ell _\perp = (\ell_1,\ell_2)$ and
examine propagation of singularities with respect to $(x_2,\xi _2)$ for the
transformed operator.

To time $T$ the shift with respect to $x_2$ will be 
$\asymp \epsilon |\ell _2|\mu ^{-1}$ and the logarithmic uncertainty principle
means that 
$\epsilon_2 |\ell_2|\mu ^{-1}T\times\varepsilon \ge C\mu ^{-1}h |\log h|$ which for $|\ell _2|\ge \epsilon$ is equivalent to $T\ge T_*$ with $T_*$ defined
by (\ref{18-6-50}).

Similarly, to time $T$ the shift with respect to $\xi_2$ will be
$\asymp \epsilon |\ell _1|\mu ^{-1}$ and the logarithmic uncertainty principle
means that $\epsilon_2 |\ell _1 |\mu ^{-1}T\times \varepsilon \ge
C\mu ^{-1}h |\log h|$ which for $|\ell _1|\ge \varepsilon$ is equivalent to
$T\ge T_*$ again.

We can justify this analysis by our standard propagation arguments with functions\begin{phantomequation}\label{18-6-55}\end{phantomequation}
\begin{equation}
\upchi \bigl(\mu \frac {x_2-y_2 } {\ell _2T} \pm \epsilon \frac{t}{T}\bigr), \qquad
\upchi \bigl(\mu \frac {\xi_2-\eta_2} {\ell _1T} \pm \epsilon \frac{t}{T}\bigr)
\tag*{$\textup{(\ref*{18-6-55})}_{1,2}$}\label{18-6-55-1}
\end{equation}
respectively.

On the other hand, let us consider the case 
$|\ell _\parallel |\ge \frac{1}{2}\epsilon$ with $\ell _\parallel=\ell_3$ and examine propagation of singularities with respect to $\xi _3$. To time $T$ the shift will be $\asymp \epsilon |\ell _3 |$ and and the logarithmic uncertainly principle means that 
$\epsilon |\ell _3 |\mu ^{-1}\times \varepsilon \ge Ch |\log h|$ which is equivalent to $T\ge T_*$ again.

We can justify this analysis by standard means with a function
\begin{equation}
\upchi \bigl(\frac {\xi _3-\eta _3} {\ell _3T} \pm \epsilon \frac{t} {T}).
\tag*{$\textup{(\ref*{18-6-55})}_{3}$}\label{18-6-55-3}
\end{equation}

Also, in all upper bound for $T$ is a small constant because for this time
$(x_2,x_3,\xi_2,\xi_3)$ keeps in $\epsilon_1$-vicinity of the point of origin\footnote{\label{foot-18-35} Actually it can be longer for small $\ell_3$ but there is no benefit from this observation.}
which can be justified by standard means with functions
\begin{phantomequation}\label{18-6-56}\end{phantomequation}
\begin{align}
&\upchi \bigl(\frac{t}{T} \pm \epsilon \mu \frac {\xi _2-\eta _2} {T} \bigr), 
&&\upchi \bigl(\frac{t}{T}\pm \epsilon \mu \frac {x_2-y_2}{T}\bigr), \tag*{$\textup{(\ref*{18-6-56})}_{1-4}$}\label{18-6-56-*}\\
&\upchi \bigl(\frac{t}{T}\pm \epsilon \frac {\xi _3-\eta _3}{T}\bigr), 
&&\upchi \bigl(\frac{t}{T}\pm \epsilon \frac {x_3- y_3}{T}\bigr).\notag
\end{align}
\end{proof}

\begin{remark}\label{rem-18-6-9} 
The above proof uses canonical form (see the next section~\ref{book_new-sect-18-7}) and therefore requires $(\bar{l},\bar{\sigma} )=(2,1)$. However, one can provide an alternative proof, similar to ones in subsection~\ref{book_new-sect-18-2-3}, which uses ``precanonical'' form (i.e. variables $Q_1,Q_2$); now it is enough
$(\bar{l},\bar{\sigma} )=(l,\sigma)\succeq (1,2)$. We leave easy details to the reader.
\end{remark}

Now statement (i) of theorem \ref{thm-18-6-12} below follows from propositions \ref{prop-18-6-5}, \ref{prop-18-6-6}, corollary \ref{cor-18-6-7} and remark \ref{rem-18-6-8}.

\section{Analysis under assumption (\ref{18-2-16})}
\label{book_new-sect-18-6-7}

Let $(l,\sigma)\succeq (2,0)$. Then we can introduce 
\begin{equation}
\nu =\epsilon |\nabla (V/F)|+\bar{\nu},\qquad 
\bar{\nu} = C_1(\mu h|\log h|)^{\frac{1}{2}}\ge C_0\mu^{-1}
\label{18-6-57}
\end{equation}
where the last inequality is due to 
$\mu \ge \epsilon (h|\log h|)^{-\frac{1}{3}}$ and as 
$|\xi_3|\le \bar{\rho}_1$ and $\nu \ge 2\bar{\nu}$ we can apply arguments of the previous subsection with
\begin{equation}
\varepsilon = C_0 \mu h|\nu^{-1} \log h|
\label{18-6-58}
\end{equation}
and take $T_*= C h|\log h| \nu^{-1} \varepsilon^{-1} \le \epsilon \mu^{-1}$ and
$T^*\asymp 1$. 

Then the contribution of the zone
$\cZ_1 \cap \{ |\nabla (V/F)|\ge 2\bar{\nu}\}$ to the Tauberian remainder is $O(h^{-2})$ and in virtue of assumption (\ref{18-2-16}) its contribution to the approximation error does not exceed expression (\ref{18-6-53})  as $l<3$.

On the other hand, contribution of zone 
$\cZ_1 \cap \{ |\nabla (V/F)|\le 2\bar{\nu}\}$ to the Weyl remainder does not exceed $C\mu h^{-2} (\mu h|\log h|)^2$ where an extra factor (in comparison with the general case) $(\log \mu h|\log h)^{\frac{1}{2}}$ is due to assumption (\ref{18-2-16}). Finally, contribution of this zone to the approximation error also does not exceed (\ref{18-6-53}).

These arguments imply statement (ii) of theorem \ref{thm-18-6-12} below.

\section{Weakly-degenerate case}
\label{thm-18-6-8}

We need this subsection to improve theorem~\ref{thm-18-6-4} to theorem~\ref{thm-18-6-10} below. Let us consider point $\bar{x}$ in which
\begin{equation}
|\nabla (V/F)| \asymp \zeta \ge C'|\log h|^{-\sigma '}, \qquad
(l,\sigma)\succeq (1,\sigma').
\label{18-6-59}
\end{equation}
Then the same is true in its $\gamma$-vicinity with $\gamma = h^\delta$ with
some small $\delta >0$ provided $C'=C'(\delta)$ is large enough if
$(l,\sigma)= (1,\sigma')$; otherwise no restriction on $C'$ is needed. Thus all
the arguments of the previous subsection remain valid but now for time
$\bar{T}=\epsilon \mu^{-1}$ the shift with respect to $\xi_3$ will be
$\asymp \zeta \bar{T}$ provided $|\nabla _\parallel (V/F)|\asymp \zeta$ and the
shift with respect to $(x_2,\xi_2)$ will be $\asymp \mu ^{-1}\zeta \bar{T}$
provided $|\nabla _\perp (V/F)|\asymp \zeta$ and therefore logarithmic uncertainty principle is fulfilled if
\begin{equation}
\varepsilon \ge C \zeta^{-1}\mu h |\log h|.
\label{18-6-60}
\end{equation}
We need to assume also that $\gamma \ge C\varepsilon$ which is always the case
provided $\mu \le h^{\delta_0 -1}$ with $\delta_0>\delta$.

Then we can take
\begin{equation}
T^* = \left\{\begin{aligned}
&\gamma \qquad &&\text{for\ \ } l=1,\\
&\epsilon \zeta \qquad &&\text{for\ \ } l>1
\end{aligned}\right.
\label{18-6-61}
\end{equation}
and the contribution of this particular element
$\psi _\gamma (x)q_\rho (hD_3)$ to the remainder estimate will be
$O\bigl(\rho \gamma ^3 h^{-2}T^{*\,-1}\bigr)$ (we assume that $q$ is compactly
supported). 

We can take always $\rho = C\bar{\rho}_1$. Recall that we estimated contribution of the partition element   $\psi _\gamma (x)\bigl(1-q_\rho (hD_3)\bigr)$ to the remainder as $O\bigl(\gamma ^3 h^{-2}\bigr)$. Therefore the contribution of $\psi _\gamma (x)$ to the remainder  is $O(\gamma ^3h^{-2})$ anyway.

Now let us consider its contribution to the approximation error. In our
assumptions it does not exceed
\begin{multline}
Ch^{-2}\gamma^3 +C\rho \gamma^3 h^{-3} \vartheta (\varepsilon) \asymp\\
Ch^{-2}\gamma^3 +
C\gamma^3 h^{-3}(\mu h |\log h|)^{l+\frac{1}{2}} \zeta ^{-l}|\log h|^{-\sigma}
\label{18-6-62}
\end{multline}
and to keep it below $C\gamma^3h^{-2}$ we need to assume that
\begin{equation}
\zeta \ge \mu ^{(2l+1)/(2l)}h^{(2l-1)/(2l)}
|\log h|^{(2l-2\sigma -1)/( 2l)};
\label{18-6-63}
\end{equation}
in particular, for $(l,\sigma)=(1,2)$
\begin{equation}
\zeta \ge \mu ^{\frac{3}{2}}h^{\frac{1}{2}}|\log h|^{-\frac{1}{2}}
\label{18-6-64}
\end{equation}
which is possible only for $\mu \le h^{-\frac{1}{3}}|\log h|^{\frac{1}{3}}$.

\section{Main theorems}
\label{book_new-sect-18-6-9}

Now we prove

\begin{theorem}\label{thm-18-6-10}
As $(\bar{l},\bar{\sigma})=(2,1)$ the conclusion of theorem~\ref{thm-18-6-4} remains true for $\mu \le h^{-\frac{1}{3}}$.
\end{theorem}

\begin{proof}
To prove this theorem we need to analyze zone 
\begin{equation}
\{ \max \bigl(\mu^{-1},(\mu h)^{\frac{1}{2}}\bigr)\le |\xi_3|\le
(\mu h|\log h|)^{\frac{1}{2}}\}.
\label{18-6-65}
\end{equation}

Let us make $\gamma$-admissible partition with $\gamma=h^\delta$. Then on the elements with 
$|\nabla (V/F)|\ge \bar{\zeta}\Def C_0 |\log h|^{-\sigma'}$ we can apply arguments of the previous subsection and conclude that their contribution to Weyl remainder does not exceed $Ch^{-2}\gamma^3$.

On the other hand, consider elements with  $|\nabla (V/F)|\le \bar{\zeta}$. Consider $\xi_3$-partition in zone (\ref{18-6-65}). Using rescaling arguments one can prove easily that contribution of such elements to the remainder does not exceed $Ch^{-2} \gamma^3\bigl(\rho   + C\mu  \rho (\mu h/\rho^2)^s \bigr)$; we leave the standard arguments to the reader. 

Then summation with respect to $\rho$ returns $Ch^{-2} \gamma^3\bigl(   + C \mu (\mu h)^{\frac{1}{2}} \bigr)$ and summation with respect to $\gamma$ returns $Ch^{-2} \bigl(   + C \mu (\mu h)^{\frac{1}{2}} \bigr)$.
\end{proof}

\begin{corollary}\label{cor-18-6-11} 
In frames of theorem~\ref{thm-18-6-4} let $(l,\sigma)\succeq (1,2)$ and $\mu \le \mu_1^*\Def h^{-\frac{1}{3}}$. Then 
$\R^\MW =O(h^{-2})$. 
\end{corollary}

Consider now non-degenerate case. We have proven already

\begin{theorem}\label{thm-18-6-12} 
Let $d=3$ and $A$ be a self-adjoint in $\sL ^2(X)$ Schr\"odinger operator \textup{(\ref{18-1-1})}. Let conditions \textup{(\ref{18-1-3})}, \textup{(\ref{18-1-6})}, \textup{(\ref{18-2-27})}
fulfilled in $B(0,1)\subset X\subset \bR ^3$. Let $(l,\sigma )\succeq (1,1)$.

Let $h^{-\frac{1}{3}}\le \mu \le C(h|\log h|)^{-1}$, $(\bar{l},\bar{\sigma})=(2,1)$. Then for two framing approximations as in footnote~\footref{foot-18-16}

\medskip\noindent
(i) Under assumption \textup{(\ref{18-2-15})} asymptotics
\begin{multline}
|\int \Bigl( \tilde{e}(x,x,0)- h^{-2} \cN_3 ^\MW (x,0)\Bigr)\psi (x)\,dx|\le \\ 
Ch ^{-2}+ Ch^{-3}\vartheta ( h|\log h|)+ 
C h^{-3} (\mu h|\log h|)^{\frac{1}{2}}\vartheta ( \mu h|\log h|)
\label{18-6-66}
\end{multline}
holds;

\medskip\noindent
(ii) Let $(l,\sigma)\succeq (2,0)$ and $l<3$. Then under assumption \textup{(\ref{18-2-16})} asymptotics \textup{(\ref{18-6-66})} holds.
\end{theorem}

We leave to the reader an easy 

\begin{problem}\label{problem-18-6-13}
Prove that 

\medskip\noindent
(i) Theorems~\ref{thm-18-6-10} and~\ref{thm-18-6-12} remain  true for $\R^\W_\infty$ and, as $\mu \le h^{-\frac{1}{2}}$ for $\R^\W$;

\medskip\noindent
(ii) Corollary~\ref{cor-18-6-11} remains true for $\R^\W$.
\end{problem}

\chapter{$d=3$: Canonical form}
\label{book_new-sect-18-7} 

From now on we will consider more difficult case
\begin{equation}
\mu \ge \bar{\mu}_1\Def C ^{-1} (h|\log h|) ^{-\frac{1}{3}};
\label{18-7-1}
\end{equation}
further restrictions will be added later. Actually due to theorem~\ref{thm-18-6-10} we can consider 
$\mu \ge \mu_1^*\Def h^{-\frac{1}{3}}$.

We already know from subsection~\ref{book_new-sect-18-5} that the contribution of the zone $\cZ^c_1\Def \{|\xi _3|\ge \bar{\rho} _1\}$ with
$\bar{\rho}_1\Def C(\mu h |\log h|)^{\frac{1}{2}}$ to the remainder does not exceed $Ch^{-2}$ and therefore we can limit ourselves by zone 
$\cZ_1\Def \{|\xi _3|\le \bar{\rho} _1\}$.

\section{Canonical form. Preliminary remarks}
\label{book_new-sect-18-7-1}

To treat this zone $\cZ_1$ properly one needs to reduce operator to the canonical form as in Section~\ref{book_new-sect-13-2}. However, as for $d=2$, to follow this procedure one must assume that $\varepsilon \ge C\mu ^{-1}|\log h|$ at least and then the approximation error will be too large unless the regularity condition holds with $l\ge 2$.

On the other hand, these arguments still work, but for much larger $\mu $. One of the reasons is that even if we neglect an approximation error, the remainder estimate would be at least $\mu h ^{-1}\times \cR _1$ where $\cR _1$ is the remainder estimate for $1$-dimensional Schr\"odinger operator. 

Unless we assume that $V/F$ is non-degenerate (i.e. either $\nabla (V/F)\ne 0$ or $\nabla^2 (V/F)\ne 0$ ) we can derive nothing better than 
$\cR _1 = O\bigl(h ^{-2/(2+l)}\bigr)$ according to Section~\ref{book_new-sect-4-5} (we neglect logarithmic factors here). Therefore we should not expect anything better than estimate $O\bigl(\mu h ^{-1-2/(2+l)}\bigr)$ for our magnetic $3$-dimensional Schr\"odinger operator anyway. For $l=1$ this estimate coincides with $O(h^{-2})$ for $\mu \asymp h^{-\frac{1}{3}}$ exactly.

However for $\mu $ close to $h ^{-\frac{1}{3}}$ we need to modify
reduction arguments properly and the rest of the section is devoted to this.

We pay mainly attention to the case
\begin{equation}
C_0(h|\log h|)^{-\frac{1}{3}}\le \mu \le C_0^{-1}(h|\log h|)^{-1}
\label{18-7-2}
\end{equation}
and analyze zone $\cZ_1$.

First of all we need to decide what does \emph{non-essential\/} mean.

\subsubsection{Non-degenerate case}
Under non-degeneracy assumption (\ref{18-2-15}) the contribution of zone $\cZ_1$ to the approximation error in $h^{-3}\cN^\MW$ does not exceed
$Ch^{-3}\bar{\rho }_1R$ where $R$ is an upper bound for a perturbation. To
prove this one needs to notice that for 
$\rho\le \rho_1^*\Def (\mu h)^{\frac{1}{2}}$ due to this assumption the measure of the zone
\begin{equation}
\cY _\rho\Def
\bigl\{\rho ^2\le\min_{n\ge 0} |\tau - V-(2n+1)\mu h |\le 2\rho^2 \bigr\}
\label{18-7-3}
\end{equation}
does not exceed $C\rho ^2 /(\mu h) $ and its contribution to the approximation
error in $h^{-3}\cN ^\MW$ does not exceed
$C\bigl( \bar{\rho}_1 \mu h^{-2}\rho ^{-1} \bigr) R\times
\rho ^2 /\mu h $.
Then the total contribution of $\cZ_1$ to the approximation error does not
exceed this expression integrated by $\rho ^{-1}d\rho $ which gives exactly $Ch^{-3}\bar{\rho}_1R$.

On the other hand, under assumption (\ref{18-2-15}) the remainder estimate will be $O\bigl( h^{-2}\bigr)$ and therefore to keep approximation error below it we
need to make calculations with the precision
\begin{equation}
R\le R^* \Def C \mu ^{-\frac{1}{2}}h ^{\frac{1}{2}}|\log h|^{-{\frac{1}{2}}}.
\label{18-7-4}
\end{equation}
Therefore 
\begin{claim}\label{18-7-5}
In the non-degenerate case as $(l,\sigma)=(1,2)$ in $\cZ_1$ we can increase $\varepsilon $ up to
\begin{equation}
\varepsilon ^* = C \mu ^{-\frac{1}{2}}h ^{\frac{1}{2}}|\log h|^{\frac{3}{2}}
\label{18-7-6}
\end{equation}
\end{claim}
\vglue-6pt

Case of the non-degeneracy assumption (\ref{18-2-16}) does not warrant its own analysis as one expects to have the remainder estimate $O(h^{-2})$ and the approximation error as above but with an extra term $C\mu h^{-2} R^2$ which would be smaller.

\subsubsection{General case}
In the general case (i.e. without non-degeneracy assumption) the
contribution of the zone $\cZ_1$ to the approximation error in $h^{-3}\cN^\MW$ does not exceed $Ch^{-3}\bar{\rho}_1R + C\mu h^{-2}R^{\frac{1}{2}}$ while we expect the remainder estimate to be
$O\bigl(h^{-2}+\mu h^{-1-2/(l+2)}|\log h|^{- \sigma /(l+2)}\bigr)$;
then to keep approximation error below it we need to make calculations with
precision
\begin{equation}
R_*=C\min \bigl( \mu ^{-\frac{1}{2}} h ^{\frac{1}{2}}|\log h|^{-\frac{1}{2}},
\mu^{-2}+ h^{2l/(2+l)}|\log h|^{-2\sigma/(l+2)}\bigr)
\label{18-7-7}
\end{equation}
where $(l, \sigma)\preceq (2,1)$. Note that $R$ given by this formula is larger
than $\varepsilon ^l |\log h|^{-\sigma}$ with $\varepsilon $ defined by (\ref{18-7-1}). Therefore 

\begin{claim}\label{18-7-8}
In the general case we can increase $\varepsilon $ up to 
\begin{equation}
\varepsilon _* = C (R_*)^{1/l} |\log h|^{ \sigma /l}.
\label{18-7-9}
\end{equation}
\end{claim}
\vglue-6pt

Later further increase will be made.
\section{Canonical form. I}
\label{book_new-sect-18-7-2}

Let us assume that $g^{jk}\in \sF^{2,1}$ even if one can weaken this assumption. We call term \emph{non-essential\/} if its $A$-bound does not exceed $R_*$ given
by (\ref{18-7-7}); under assumption (\ref{18-2-15}) or (\ref{18-2-16}) we can increase it to $R^*$.

Therefore, in the zone $\cZ_1$ we have
\begin{equation}
\varepsilon = C(\mu ^{-1}h|\log h|)^{\frac{1}{2}}
\label{18-7-10}
\end{equation}
or larger. Note that for $\varepsilon $ given by (\ref{18-7-10}) mollification error $\varepsilon ^l |\log h|^{-\sigma}$ does not exceed $R_*$.

To transform operator $A$ to the canonical form we would like first to transform
its principal part to the ``diagonal form'' as it was done for $d=2$. This can
be achieved by the same transformation (\ref{18-3-4}) 
\begin{equation}
T(t)\Def e^{-i\mu ^{-1}h^{-1}tL^\w}\qquad\text{with\ \ }
L^\w={\frac{1}{2}}\sum_{j,k} P_jL ^{jk}P_k.
\tag{\ref{18-3-4}}
\end{equation}
Then (\ref{18-3-5})--(\ref{18-3-8}) are preserved but now we have $3\times 3$ commutator matrix 
\begin{gather*}
J=\begin{pmatrix} J'& 0\\ 0& 0&\end{pmatrix}\qquad\text{with\ \ }
J'=\begin{pmatrix} 0& -1\\ 1& 0\end{pmatrix}\\
\shortintertext{and}
\Lambda =f(x) \ell J = \begin{pmatrix} \Lambda ' & 0\\ \Lambda ''&0\end{pmatrix}
\qquad\text{with\ \ }
\Lambda '= \begin{pmatrix} \ell^{12}& -\ell^{21}\\ \ell ^{22}& -\ell^{21}\end{pmatrix},\quad
\Lambda ''= (\ell ^{32}\; -\ell^{31}).
\shortintertext{Then}
e^{t\Lambda }=\begin{pmatrix} e^{t\Lambda '} & 0\\
\Lambda '' K(t) & 1\end{pmatrix}
\end{gather*}
with matrix $K(t)=(\Lambda ')^{-1}(e^{t\Lambda '}-I)$ which is well defined even
if $\Lambda '$ is degenerate. 

Therefore we can achieve $e^{^t\!\Lambda} G e^\Lambda $ to be diagonal with the top-left element equal to $1$. However we cannot change $g^{33}$ by this method unless we use ``much larger'' operator $L$ (something like 
$\omega (x) \mu h D_3$) which is not a good idea at this moment in the non-smooth case.

The rest of the analysis of subsubsection~\ref{book_new-sect-18-3-1-1}.1 does not change and and therefore we get that 
\begin{claim}\label{18-7-11}
Modulo operator with $A$-bound not exceeding $\mu^{-1}h$
\begin{gather}
T(-1)A_0T(1) \equiv P_1^2 +P_2(g ')^\w P_2 + P_3(g '' )^\w P_3
+\mu ^2 M^\w, \label{18-7-12}\\[2pt]
M \Def
a_0 \circ \phi _1 -p_1^2-g^{\prime 22}p_2^2- g^{\prime 33}p_2^2 \label{18-7-13}
\end{gather}
\end{claim}
(compare with (\ref{18-3-18}),(\ref{18-3-19})).

Recall that according to (\ref{18-3-13}) $M=\sum _{i,j,k}\beta _{ijk}\, p_ip_jp_k$ with $\beta _{ijk}\in \sF^{1,1}$. Finally, modulo non-essential operator
\begin{equation}
T(-1)VT(1)\equiv (V\circ \phi_1)^\w.
\label{18-7-14}
\end{equation}

Now, we can assume that $g ',g''$ are functions of $x$ only because we can achieve it modifying $M$ but not changing its type. Further, condition $F=1$
means that $\sqrt{g '} f=1$ and we can apply the same construction as in subsubsection~\ref{book_new-sect-18-3-1-2}.2 and therefore we can assume without any loss of the generality that
\begin{gather}
V_1=V_3=0, \quad P_1=hD_1 \implies V_2=V_2(x_2,x_3),\label{18-7-15}\\
V_2,\partial_{x_1}V_2\in \sF^{2,1} \label{18-7-16}\\[2pt]
\shortintertext{and}
p_2 = \alpha (x,\xi _2 )\bigl(x_1 - \lambda (x_2,x_3,\xi _2)\bigr),
\qquad \alpha, \lambda \in \sF^{2,1}.\label{18-7-17}
\end{gather}

Repeating then construction of subsubsection~\ref{book_new-sect-18-3-1-2}.2, we take the next transformation $T'(t)=e^{-i\lambda ^\w D_1}$. As a result we will get operator
\begin{gather}
T'(-1)T(-1)A_0T(1)T'(1)\equiv
\bar{A}_0 + \bigl(\mu^2 M'\bigr)^\w,\label{18-7-18}\\[3pt]
\bar{A}_0\Def\bar{P}_1^2 + \bar{P}_2^2 +\bar{P}_3\alpha _0 ^\w \bar{P}_3,
\qquad \bar{P}_1=hD_1,\;
\bar{P}_2=-\mu x_1,\; \bar{P}_3=hD_3,\label{18-7-19}\\[3pt]
M'= a_0\circ \Phi - \bar{a}_0=
\sum _{i,j,k}\beta '_{ijk}\bar{p}_i\bar{p}_j\bar{p}_k\label{18-7-20}
\end{gather}
with $\bar{p}_1\Def \xi_1$, $\bar{p}_2\Def -x_1$, $\bar{p}_3\Def \xi_3$,
$\Phi\Def \phi _1\circ \phi'_1$, where $\phi_t,\phi'_t$ are
corresponding Hamiltonian flows,
$\alpha _0 =\alpha _0 (x_2,x_3,\xi_2) \in \sF^{2,1}$ (compare with
(\ref{18-3-28})--(\ref{18-3-30})).

Further, modulo operator with $\bar{A}_0$-bound not exceeding
$C\vartheta (\varepsilon)$
\begin{equation}
T'(-1)T(-1)A_0T(1)T'(1)\equiv (V\circ \Phi )^\w.
\label{18-7-21}
\end{equation}

So far precision was better than we needed: it was
$\varepsilon ^l |\log h|^{-{\sigma}}$, exactly as in subsubsection~\ref{book_new-sect-18-3-1-2}.2. We could achieve $\alpha _0=1$ but it would be not very useful.

\section{Canonical form. II}
\label{book_new-sect-18-7-3}

Now we will follow subsubsection~\ref{book_new-sect-18-3-1-3}.3 but with a twist. Namely, we define $S'$, $S''$, $W'$, $W''$ \emph{exactly\/} by (\ref{18-3-39})--(\ref{18-3-42}) and apply transformation
$T''(-t)=e^{-it\mu h^{-1}S^\w}$ with the phase $S=S'+S''$. Then modulo non-essential term
\begin{equation}
T''(-t)\Bigl( \bar{A}_0 + \bigl(\mu^2 M' + V'\bigr)^\w\Bigr) T''(t)
\equiv
\Bigl( \bigl( \mu ^2 \bar{a}_0 + \mu^2 M' + V'\bigr) \circ \psi_t \Bigr)^\w
\label{18-7-22}
\end{equation}
where $\psi_t$ is a corresponding Hamiltonian flow. Let us calculate symbol
\begin{multline*}
\bigl( \mu ^2 \bar{a}_0 + \mu^2 M' + V'\bigr) \circ \psi_t
=\mu ^2 \bar{a}_0 + \bigl(\mu^2 M' + V'\bigr) - t\mu ^2 \{\bar{a}_0, S\}\\
-t \{ \mu^2 M' + V', S\} -
{\frac {t^2} 2} \{\{\bar{a}_0, S\}, S\} +\dots ,
\end{multline*}
and plug $t=1$. Then
\begin{align*}
\mu ^2 \{\bar{a}_0, S\} =
2\mu ^2 (\xi _1 \partial _{x_1}-x_1 \partial _{\xi_1})S\ +\ &
2\mu ^2 \alpha_0 ^2\xi_3 (\partial _{x_3}S)\ +\ \{\alpha _0^2 ,S\} \xi_3^2
\\[3pt]
= \mu ^2 M'+ V' - W_0\ + \
&2\mu ^2 \alpha_0 ^2\xi_3 (\partial _{x_3}S)\ +\ \{\alpha _0^2 ,S\} \xi_3^2
\end{align*}
where the last equality is due to (\ref{18-3-39}), (\ref{18-3-40}) and $W_0=W'+W''$.

Further, let us notice that 
$|\xi_3|\le \mu^{-1}\bar{\rho}_1 = C(\mu ^{-1}h|\log h|)^{\frac{1}{2}}$ in our zone $\cZ_1$ where an extra factor $\mu^{-1}$ appears because in the reduction
part we consider $\mu^{-1}h$-pseudo-differential operator. Then one can see easily that modulo non-essential operator (in the general sense) our transformed operator equals to the quantization of
\begin{equation}
\mu ^2\bar{a}_0 +W_0 -2 \mu ^2 \alpha_0 ^2\xi_3 (\partial _{x_3}S)
-{\frac{1}{2}} \{\mu ^2M'+V'-W_0,S\}'
\label{18-7-23}
\end{equation}
where $\{.,.\}'$ again mean ``short'' Poisson brackets (with respect to
$x_1,\xi_1$ only).

 Obviously the last term in (\ref{18-7-23}) is $O(\mu^{-2})$ and it is essential only for $\mu \le \hat{\mu}_1\Def 
 \epsilon_0 h^{-\frac{1}{3}}|\log h|^{\frac{1}{3}}$.

Furthermore,
$\partial _{x_3}S= O\bigl(\mu ^{-1-l}|\log h|^{-\sigma}+\mu ^{-3}\bigr)$
and therefore the second term in (\ref{18-7-23}) is
$O\bigl((\mu ^{1-l}|\log h|^{-\sigma}+\mu ^{-2}) \cdot
(\mu ^{-1}h|\log h|)^{\frac{1}{2}} \bigr)$ which does not exceed $R^*$ but is not necessarily smaller than $R_*$.

Therefore, if one of the conditions (\ref{18-2-15}), (\ref{18-2-16}) is fulfilled and $\mu \ge \hat{\mu}_1$ then we have finished because only the first term in (\ref{18-7-23}) is essential. In this case we have arrived to the  operator
\begin{gather}
\bar{A}_0 +
W_0(x_2,x_3,\mu ^{-1}hD_2, \mu ^{-1}hD_3, \mu ^{-1}\bar{A}_0^{\frac{1}{2}})^\w
\label{18-7-24}\\
\shortintertext{with}
\quad W_0(x_2,x_3,\xi_2,\xi_3,\rho)=\cM_\rho (\mu^2 M'+V').
\label{18-7-25}
\end{gather}

Furthermore, under one of the assumptions (\ref{18-2-15}), (\ref{18-2-16}) we will see later that for  $\mu \le \hat{\mu}_1$ we do not need the full canonical form at all.

However, in the general case the second and the third terms in (\ref{18-7-23}) are essential and we need to continue in the same way as in subsubsection~\ref{book_new-sect-18-3-1-4}.4. Let us define $S_1$, $S'''$, $W_1$, $W'''$ in the same way as before:
\begin{gather}
W'''= \alpha_0^2 \partial _{x_3} \cM_\rho S,\qquad
W _1= -{\frac{1}{2}} \cM_\rho \{\mu ^2M'+V'-W_0,S\}',\label{18-7-26}\\[3pt]
2\mu ^2 (\xi _1 \partial _{x_1}-x_1 \partial _{\xi_1})S''' =
\alpha_0^2 \partial _{x_3}(S- \cM_\rho S),\label{18-7-27}\\[3pt]
2\mu ^2 (\xi _1 \partial _{x_1}-x_1 \partial _{\xi_1})S_1=
-{\frac{1}{2}} \{\mu ^2M'+V'-W_0,S\}'+ W_1;\label{18-7-28}
\end{gather}
Recall that $\alpha_0=\alpha_0 (x_2,x_3,\xi_2)$. Let us apply $\mu^{-1}h$-pseudo-differential operator
transformation 
\begin{equation*}
T'''(t)= \bigl(e^{-it\mu h^{-1}(S'''\xi_3 +S_1)} \bigr)^\w;
\end{equation*}
then
\begin{multline*}
T'''(-1) \Bigl(\mu ^2 \bar{a}_0 + W_0 -
2 \mu ^2 \alpha_0 ^2\xi_3 (\partial _{x_3}S) - 
\frac{1}{2} \{\mu ^2M'+V'-W_0,S\}'\Bigr)^\w T'''(1) \equiv\\[3pt]
\Bigl(\mu ^2 \bar{a}_0 + W_0 -
2 \mu ^2 \alpha_0 ^2\xi_3 (\partial _{x_3}S) - 
\frac{1}{2} \{\mu ^2M'+V'-W_0,S\}'-\\
\mu ^2 \{\bar{a}_0, S'''\xi_3 +S_1 \} +\dots \Bigr)^\w.
\end{multline*}
Plugging (\ref{18-7-26})--(\ref{18-7-27}) we see easily that
\begin{multline}
T'''(-1) \Bigl(\mu ^2 \bar{a}_0 + W_0 -
2 \mu ^2 \alpha_0 ^2\xi_3 (\partial _{x_3}S) -\\
{\frac{1}{2}} \{\mu ^2M'+V'-W_0,S\}'
\Bigr)^\w T'''(1) \equiv\\
\bar{A}_0 + W x_2,x_3,\mu ^{-1}hD_2, \mu ^{-1}hD_3, B)
\label{18-7-29}
\end{multline}
where $W=W_0+W'''\xi_3 +W_1$.

The last step is not really necessary, but we'll do it for the sake of the simple canonical form. 

First of all, transformation $T^V(1)$, with 
$T^V(t) = e^{-it\omega ^\w D_3 }$ and with
$\omega = \omega (x_2,x_3,\xi_2 )\in \sF^{2,1}$, transforms
$D_3$ into $(\partial _{x_3}\omega )^\w D_3+ \beta ^\w \mu^{-1}hD_3^2$ with
$\beta = \beta (x_2,x_3,\xi_2,\xi_3 )\in \sF^{1,1}$ and therefore one can
transform $\alpha_0^\w D_3$ into $D_3 +\beta ^\w \mu^{-1}hD_3^2$. Thus, modulo
non-essential operator
\begin{equation}
T^V(-1)\bar{A}_0 T^V(1)\equiv h^2 D_1^2+\mu ^2 x_1^2+D_3^2.
\label{18-7-30}
\end{equation}
One can also see easily that 
$T^V(-1)W^\w T^V(1)\equiv W ^\w + K^\w \mu ^{-1}hD_3$ with
$K \in \sF^{l-1,\sigma }$.

Further, modulo non-essential term we can rewrite our reduced operator in the
form $\bar{A}_0 + \bar{W}^\w + \bar{K}^\w \mu ^{-1}hD_3$ with $\alpha_0=1$ in
the expression (\ref{18-7-19}) for $\bar{A}_0$ and
$\bar{W}=\bar{W}(x_2,x_3,\xi_2, \mu^{-1}\bar{A}_0^{\frac{1}{2}})$ where
$\bar {W}\Def W|_{\xi_3=0} \in \sF^{l,\sigma}$ and 
$\bar{K}=\bar{K}(x_2,x_3,\xi_2, \mu^{-1}\bar{A}_0^{\frac{1}{2}})$ with 
$\bar{K}\Def (\partial_{\xi_3}W)|_{\xi_3=0}\in \sF^{l-1,\sigma}$.

Finally, applying transformation $T^{VI}(1)$ with
$T^{VI}(t) = e^{-it\mu^{-1}h^{-1}{\omega '}^\w }$ and
$\omega '= \omega '(x_2,x_3,\xi_2, \mu^{-1}\bar{A}_0^{\frac{1}{2}})$,
$\omega '\in \sF^{l-1,\sigma}$
we will get (for an appropriate symbol $\omega '$) modulo negligible operator
$\bar{A}_0 + \bar{W} ^\w$.

\section{Canonical form. III. $\mu \ge \epsilon h^{-1}|\log h|^{-1}$}
\label{book_new-sect-18-7-4}

In this case again an $A$-bound for $\mu^{-1} P_j$ does not exceed
$C\varepsilon $ and construction is absolutely straightforward as in the smooth
case; we just note that $R^*= Ch$ and
\begin{multline}
R_*= C\min\bigl( h, h^{{2l}/(2+l)}
|\log h|^{- {2\sigma }/(l+2)}\bigr)=\\[2pt]
C\left\{\begin{aligned}
&h \qquad &&\text{for\ \ } (l,\sigma)\preceq (2,0),\\[2pt]
&h |\log h|^{-{\frac \sigma 2}} \qquad&&\text{for\ \ } l=2,0<\sigma\le 1.
\end{aligned}\right.
\label{18-7-31}
\end{multline}

\chapter{$d = 3$: Tauberian theory}
\label{book_new-sect-18-8}

In this section we derive Tauberian estimates. First of all, we have decomposition similar to one we had in the $2$-dimensional case. 

Using this decomposition we  analyze the case $\mu \ge \bar{\mu}_1$ in two different settings: under one of 
non-degeneracy assumptions (\ref{18-2-15}), (\ref{18-2-16}) we prove that the Tauberian remainder estimate does not exceed $C(1+\mu h)h^{-2}$ and that in the general case we prove that it does not exceed expression 
$Ch^{-2}+ C\mu h^{-1}\cR_1$ where $\cR_1$ is the remainder estimate for $1$-dimensional Schr\"odinger operator.

\section{Decomposition}
\label{book_new-sect-18-8-1}

First of all, as we already mentioned, we have decomposition, similar to (\ref{18-4-3})--(\ref{18-4-4}): we reduce our operator to the family of $1$-dimensional Schr\"odinger operators with respect to $x_3$ which are also $\mu^{-1}h$-pseudo-differential operators with respect to $x_2$
\begin{equation}
\cA _n= hD_3^2+ r_n^2 + W(x_2,x_3,\mu^{-1}hD_2,\mu^{-1}r_n),\qquad
r_n = \bigl((2n+1)\mu h\bigr)^{\frac{1}{2}},
\label{18-8-1}
\end{equation}
with $x'=(x_1,x_2)$ and arrive to decomposition
\begin{equation}
U(x,y,t) \equiv \cT
\Bigl(\sum _{n\ge 0} U_n(x',y',y_3, t) \Upsilon _n (x_1)\Upsilon _n (y_1)\Bigr) \cT^{-1}
\label{18-8-2}
\end{equation}
where $U_n(x',y',t)$ are Schwartz kernels of the propagators for $\cA _n$. 

Here we assume that $\mu \le C_0 h^{-1}$; as $\mu \ge C_0h^{-1}$ we consider Schr\"odinger-Pauli operator with $\fz=1$ and $r_n=0$.

\section{Reduction of $\cZ_1$} 
\label{book_new-sect-18-8-2}

As we mentioned, in this case we have essentially $1$-dimensional Schr\"odinger
operator (or rather a family of them). Due to section~\ref{book_new-sect-18-6} we need to consider zone $\cZ_1= \{(x,\xi):\ |\xi_3|\le \bar{\rho}_1= C(\mu h|\log h|)^{\frac{1}{2}}\}$ only. 

Assuming first that $\vartheta(\eta)\ge \eta^2$ et us introduce admissible functions $\ell =\ell_{(n)}$ and
$\varrho =\varrho_{(n)}$ of $(x_2,x_3,\xi_2)$ similar to those introduced in
section~\ref{book_new-sect-5-1}:
\begin{multline}
\ell = \epsilon \min \Bigl\{\eta :\  |W+(2n +1)\mu h|\le \vartheta (\eta),
|\nabla W|\le \eta ^{-1}\vartheta (\eta )\Bigr\}+\bar{\ell},\\
\varrho = \vartheta (\ell )^{\frac{1}{2}}
\label{18-8-3}
\end{multline}
with $\bar{\ell}$ defined from equation
\begin{equation}
\bar{\ell} \vartheta (\bar{\ell} )^{\frac{1}{2}}=C_1 h|\log h|.
\label{18-8-4}
\end{equation}
Note that as $\vartheta(\eta)\le \eta |\log \eta|^{-2}$ we conclude that 
\begin{equation}
\bar{\ell}\ge (h|\log h|)^{\frac{2}{3}}\ge \bar{\varepsilon}= (\mu^{-1}h)^{\frac{1}{2}} |\log h|
\label{18-8-5}
\end{equation}
where the last inequality is due to $\mu \ge \bar{\mu}_1$. Recall that $\bar{\varepsilon}$ was the approximation parameter in the inner zone.

We need to have approximation parameter $\varepsilon_n=\varepsilon_{(n)}$ such that $\varepsilon \varrho_n \ge Ch|\log h|$ to satisfy uncertainty principle with respect to $x_3$ and therefore we replace $\varepsilon _{(n)}$ by
\begin{equation}
\varepsilon\Def \varepsilon_{(n)} + \varepsilon'_{(n)},\qquad 
\varepsilon'_{(n)}\Def Ch\varrho_{(n)}^{-1} |\log h|;
\label{18-8-6}
\end{equation}
then 
\begin{equation}
\varepsilon_{(n)}\lesssim \ell_{(n)} \quad \text{and}\quad 
\varepsilon_{(n)} \asymp \ell_{(n)} \iff \ell_{(n)}\asymp \bar{\ell}.
\label{18-8-7}
\end{equation}

Note first that

\begin{claim}\label{18-8-8}
If on the given partition element 
\begin{equation}
|W+(2n +1)\mu h|\asymp \varrho^2=\vartheta (\ell)\ge C\mu h|\log h|
\label{18-8-9}
\end{equation}
then on this element operator $\cA_n$ is elliptic in $\cZ_1$ and its contribution to $F_{t\to h^{-1}\tau} \bar{\chi}_T (t)U_n $ is negligible on energy levels $\tau:|\tau|\le \epsilon_1\varrho^2$.
\end{claim}
 
On the other hand, 
 
\begin{claim}\label{18-8-10}
If on the given partition element 
\begin{equation}
|\nabla W|\asymp \varrho^2\ell^{-1}=\vartheta (\ell)\ell^{-1},\qquad \varrho^2 \ge C\mu h|\log h|
\label{18-8-11}
\end{equation}
then on this element operator $\cA_n$ is microhyperbolic with respect to $(x_2,x_3,\xi_2)$ in $\cZ_1$ and its contribution to Tauberian estimate  with the given index $n$ does not exceed $C\mu h^{-1} \ell^2$.
\end{claim}

Really, microhyperbolicity with respect to $(x_2,\xi_2)$ and propagation with respect to $(\xi_2,x_2)$ gives us 
\begin{equation}
T_* = Ch|\log h|/(\varrho^2\ell^{-1}\varepsilon)
\label{18-8-12}
\end{equation}
and microhyperbolicity with respect to $x_3$ and propagation with respect to $\xi_3$ gives us the same answer as in these case we consider $\mu^{-1}h$-pseudo-differential and $h$-pseudo-differential operators respectively. One can see easily that 
$\varrho^2 \ell^{-1}\ge C(\mu h|\log h|)^{\frac{1}{2}}$ and therefore $T_*\le \epsilon $. 

\begin{remark}\label{rem-18-8-1}
However the actual estimate is based not on the propagation but on fine subpartition and ellipticity arguments exactly like in section~\ref{book_new-sect-18-4}. We leave details to the reader.
\end{remark}

Meanwhile one can see easily that 

\begin{claim}\label{18-8-13}
On  elements with $\varrho\ge (\mu h)^{\frac{1}{2}}$  additional approximation is not needed because  $Ch\varrho^{-1} |\log h|\le C\bar{\varepsilon}$.
\end{claim}

Since over such element with $\varrho^2 \ge C\mu h$ ellipticity is violated only for $\asymp \varrho^2 /\mu h$ indices $n$ then due to (\ref{18-8-10}) its contribution to the Tauberian remainder   does not exceed 
$C\mu h^{-1} \ell^{-1} \times \varrho^2 /\mu h \asymp 
Ch^{-2}\varrho^2  \ell^2\le Ch^{-2}\ell^3$. Then summation over all such elements results in $O(h^{-2})$.

Therefore we arrive to

\begin{proposition}\label{prop-18-8-2}
As $(l,\sigma)\preceq (1,2)$  contribution to the Tauberian remainder with $T=\bar{T}=\epsilon$ of the complement of the zone
\begin{equation}
\cZ\Def \{(x,\xi):  |\xi_3|+\varrho \le C \bar{\rho}_1\}
\label{18-8-14}
\end{equation}
with $\varrho\Def \min _n \varrho_n$, is $O(h^{-2})$.
\end{proposition}

Therefore we need to consider zone $\cZ$ only and we divide it into an \emph{inner zone}
\begin{equation}
\cZ_\inter \Def \{ (x,\xi):\ |\xi_3|+\varrho \le \bar{\rho}_0= 
\epsilon _0 (\mu h)^{\frac{1}{2}} \}
\label{18-8-15}
\end{equation}
and a \emph{transitional zone}
\begin{equation}
\cZ_\trans \Def
\{ (x,\xi):\ \bar{\rho}_0 \le |\xi_3|+\varrho \le \bar{\rho}_1  \}.
\label{18-8-16}
\end{equation}

\begin{corollary}\label{cor-18-8-3}
As $(l,\sigma)\preceq (2,0)$  the Tauberian remainder is $O(h^{-2})$ as 
$\mu \le (h|\log h|)^{-1}$.
\end{corollary}

We leave to the reader to prove using rescaling technique (with defined above $\ell, \varrho$) a bit more stronger

\begin{proposition}\label{prop-18-8-4}
As $(l,\sigma)\preceq (2,0)$ contribution to the Tauberian remainder with $T=\bar{T}=\epsilon$ of the complement of the zone $\cZ_\inter$
with $\varrho\Def \min _n \varrho_n$ is $O(h^{-2})$.
\end{proposition}

\begin{corollary}\label{cor-18-8-5}
As $(l,\sigma)\preceq (2,0)$ the Tauberian remainder is $O(h^{-2})$ as
$\mu \lesssim h^{-1}$.
\end{corollary}

\section{Analysis in $\cZ_\inter$} 
\label{book_new-sect-18-8-3}

Let us consider more difficult inner zone; as $\mu \ge h^{-1}$ and we consider Schr\"odinger-Pauli operator with $\fz=1$. Then we can apply the same arguments as above but over one partition element ellipticity assumption is violated for no more than $1$ index $n$. Further, as $\ell \asymp \bar{\ell}$ (and then $\bar{\varrho}\bar{\ell}\asymp \mu h|\log h|$ contribution of the element to the asymptotics for given index $n$ does not exceed 
$C\mu h^{-2}\bar{\varrho}\bar{\ell}^3 \asymp C\mu h^{-2}|\log h| \bar{\ell}^2$.

Therefore we arrive to

\begin{proposition}\label{prop-18-8-6}
As $(l,\sigma)\preceq (2,0)$ contribution to the Tauberian remainder with $T=\bar{T}=\epsilon$ of the zone $\cZ_\inter$ does not exceed
\begin{equation}
C\mu h^{-1}\Bigl( \int _{\{\ell\ge C_0\bar{\ell}\}}  \ell^{-1}\,dx'd\xi_2 +
\bar{\ell}^{-1} |\log h| \int _{\{\ell\le C_0\bar{\ell}\}}\,dx'd\xi_2 \Bigr)
\label{18-8-17}
\end{equation}
and an extra approximation error does not exceed
\begin{multline}
C\mu h^{-2} \int _{\{\ell\ge C_0\bar{\ell}\}}
 \vartheta \bigl(h|\log h|/\varrho\bigr) \,dx'd\xi_2 +\\
C\mu h^{-1}|\log h| \bar{\ell}^{-1}  \int _{\{\ell\le C_0\bar{\ell}\}} \,dx'd\xi_2 
\label{18-8-18}
\end{multline}
with $\ell\Def \min _n \ell_n$, $\varrho =(\vartheta(\ell))^{\frac{1}{2}}$ and $\bar{\ell}$ defined by \textup{(\ref{18-8-4})}.
\end{proposition}

Then, under non-degeneracy assumptions (\ref{18-2-15}) (as $(l,\sigma)\succeq (1,2)$) or (\ref{18-2-16}) (as $(l,\sigma)= (2,0)$) both expressions (\ref{18-8-17}) and (\ref{18-8-17}) do not exceed $C\mu h^{-1}$ while in the general case it does not exceed $C\mu h^{-1}|\log h| \bar{\ell}^{-1}$. Combining with proposition \ref{prop-18-8-4} as as $\mu h\lesssim 1$ we arrive to 

\begin{proposition}\label{prop-18-8-7}
(i) Let $\mu \lesssim h^{-1}$   and either $(l,\sigma)\succeq (1,2)$ and condition \textup{(\ref{18-2-15})} be fulfilled or $(l,\sigma)\succeq (2,0)$ and condition \textup{(\ref{18-2-16})} be fulfilled. 

Then both the Tauberian remainder and an additional approximation error do not exceed $Ch^{-2}$;

\medskip\noindent
(ii) Let $\mu \gtrsim h^{-1}$   and either $(l,\sigma)\succeq (1,2)$ and condition \ref{18-2-15-*} be fulfilled or $(l,\sigma)\succeq (2,0)$ and condition \ref{18-2-16-*}  be fulfilled. 

Then both the Tauberian remainder and an additional approximation error do not exceed $C\mu h^{-1}$;

\medskip\noindent
(iii) In the general case as $(\ell,\sigma)\preceq (2,0)$ both the Tauberian remainder and an additional approximation error do not exceed
\begin{equation}
 Ch^{-2}+ C\mu h^{-1}|\log h| \bar{\ell}^{-1}.
\label{18-8-19}
\end{equation}
\end{proposition}

We leave to the reader  

\begin{Problem}\label{problem-18-8-8}
Using rescaling arguments prove that as $(l,\sigma)\preceq (2,0)$ both contribution to the Tauberian remainder with $T=\bar{T}=\epsilon$ of the zone $\cZ_\inter$ and extra approximation error do not exceed
\begin{equation}
Ch^{-2} + C\mu h^{-1}\bar{\ell}^{-1}
\tag*{$\textup{(\ref*{18-8-19})}^*$}\label{18-8-19-*}
\end{equation}
with  $\bar{\ell}$ redefined by 
\begin{equation}
\bar{\ell} \vartheta (\bar{\ell} )^{\frac{1}{2}}=C_1 h.
\tag*{$\textup{(\ref*{18-8-4})}^*$}\label{18-8-4-*}
\end{equation}
where $\varepsilon$ is also redefined as
\begin{equation}
\varepsilon\Def \varepsilon_{(n)} + \varepsilon'_{(n)},\qquad 
\varepsilon'_{(n)}\Def Ch\varrho_{(n)}^{-1} (h/(\varrho_{(n)}\ell_{(n)}))^{-\delta'};
\tag*{$\textup{(\ref*{18-8-6})}^*$}\label{18-8-6-*}
\end{equation}
\end{Problem}

We can generalize estimate \ref{18-8-19-*} to $(l,\sigma) \succ (2,0)$. Namely, as   $(m,0)\preceq (l,\sigma)\prec (m+1,0)$  we introduce scaling functions
\begin{multline}
\ell =\ell_m\Def  \\
\epsilon \min \Bigl\{\eta :\  |\nabla^\alpha \bigl(W+(2n +1-\fz)\mu h\bigr) |\le \vartheta (\eta) \eta^{-|\alpha|}\ 
\shoveright{\forall \alpha:|\alpha|\le m\Bigr\},}\\
\varrho = \vartheta (\ell )^{\frac{1}{2}}
\label{18-8-20}
\end{multline}
and define $\bar{\ell}$ from  \ref{18-8-4-*}.

Let us consider $\ell$-admissible partition. Then we can apply the same arguments as before to cover zone $\{(x_2,x_3,\xi_2): \varrho^2 \gtrsim \mu h\}$ and to prove that its contribution to the Tauberian remainder does not exceed $Ch^{-2}$. 

Meanwhile contribution of zone  $\{(x_2,x_3,\xi_2):\ \ell \lesssim \bar{\ell}\}$ does not exceed $C\mu h^{-2}\bar{\varrho} \asymp C\mu h^{-1}|\log h| \bar{\ell}^{-1}$ with  $\bar{\varrho} =\vartheta (\bar{\ell} )^{\frac{1}{2}}$.

So, we need to cover zone $\{(x_2,x_3,\xi_2), \ell \ge C_0\bar{\ell}\}$. Consider particular element, and on this element 
\begin{equation}
\sum _{\alpha:|\alpha|\le m}
|\nabla^\alpha \bigl(W+(2n +1-\fz)\mu h\bigr) | \ell^{|\alpha|}\asymp \varrho^2.
\label{18-8-21}
\end{equation}

We will need 

\begin{proposition}\label{prop-18-8-9}
The contribution of each such element both to the Tauberian remainder  and to an approximation error for two appropriate framing approximations do not exceed 
\begin{equation}
Ch^{-2}\ell^3 + C\mu h^{-1}\ell^2 (h/\varrho \ell)^{-\delta}
\label{18-8-22}
\end{equation}
with arbitrarily small exponent $\delta>0$.
\end{proposition}

\begin{proof}
The proof is due to arguments of  Section~\ref{book_new-sect-5-1} based on the series of the scaling functions $\ell_k,\varrho_k$  with $k=m-1,\ldots, 1$ defined by (\ref{18-8-20}) with $\varrho (\eta)=\eta^{k+1}$ and  $\varepsilon_k$ defined by \ref{18-8-6-*} with $\delta'=\delta'_k>0$, with the following modifications:

\medskip\noindent
(i) When $\varrho_k^2 \ge \mu h$ (as $\mu h \le C_0$ only) we have an extra factor $\varrho_k^2/(\mu h)$ and therefore the first term in (\ref{18-8-21}) carries no extra factor (unlike the second term);

\medskip\noindent
(ii) Due to factor $(h/\varrho \ell)^{-\delta}$ in the second term we do not need sophisticated arguments of subsubsection~\ref{book_new-sect-5-1-2} - the only ones where we needed a large smoothness.

Note that $\ell (h/\varrho \ell)^\delta\ge \bar{\ell}$. Really  as $\delta$ is small enough this expression is monotone increasing function of $\ell$; recall that $\ell\ge \bar{\ell}$. 
\end{proof}

Therefore expression (\ref{18-8-22}) does not exceed
$C\bigl(h^{-2} + C\mu h^{-1}\bar{\ell}^{-1}\bigr)\ell^3$
and summation over partition returns 
$C\bigl(h^{-2} + C\mu h^{-1}\bar{\ell}^{-1}\bigr)$.

\begin{remark}\label{rem-18-8-10}
It follows from the analysis of the problem~\ref{problem-18-8-8} that under an extra non-degeneracy assumption
\begin{phantomequation}\label{18-8-23}\end{phantomequation}
\begin{equation}
\sum_{\alpha:\,1\le |\alpha |\le m}|\nabla^\alpha W|\ge \epsilon_0
\tag*{$\textup{(\ref*{18-8-23})}_m$}\label{18-8-23-m}
\end{equation}
(as $\mu h\lesssim 1$) and
\begin{equation}
\sum_{\alpha:\, |\alpha |\le m}|\nabla^\alpha \bigl(W+(2n+1)-\fz\bigr)|\ge \epsilon_0
\tag*{$\textup{(\ref*{18-8-23})}^*_m$}\label{18-8-23-m-*}
\end{equation}
(as $\mu h\gtrsim 1$) and with $(m,0)\prec (l,\sigma)$ both both contribution to the Tauberian remainder with $T = \bar{T}=\epsilon$  of the zone $\cZ_\inter$ and extra approximation error do not exceed 
\begin{equation}
Ch^{-2}+C\mu h^{-1-\delta}.
\label{18-8-24}
\end{equation}
\end{remark}

We leave to the reader the following

\begin{Problem}\label{Problem-18-8-11}
Improve (\ref{18-8-24}) to \begin{phantomequation}\label{18-8-25}\end{phantomequation}
\begin{equation}
Ch^{-2}+C\mu h |\log h|^{m-1}.
\tag*{$\textup{(\ref*{18-8-25})}_m$}\label{18-8-25-m}
\end{equation}
\end{Problem}

To do this one needs to improve (\ref{18-8-22}) in an obvious way.

\chapter{$d = 3$: Calculations and  main theorems}
\label{book_new-sect-18-9}

\section{Assembling what we got}
\label{book_new-sect-18-9-1}

After corresponding estimates were derived in the previous section, following subsection~\ref{book_new-sect-13-4-2}, we can rewrite (without increasing an error)  the Tauberian expression 
\begin{equation*}
h^{-1}\int_{-\infty}^0 F_{t\to h^{-1}\tau }
\Bigl(\bar{\chi}_T(t)\Gamma \bigl(\psi  U \bigr)\Bigr)\,d\tau
\end{equation*}
as 
\begin{equation*}
\frac {1}{4\pi^2} \mu h^{-2}\sum_n\int
\uptheta \bigl( - \xi_3^2 -(2n+1)\mu h- W (x_2,x_3,\xi_2)\bigr)
\alpha (x_2,x_3,\xi_2,\xi_3)\, dx_2 dx_3 d\xi_2 d\xi_3
\end{equation*}
which in turn we can rewrite as 
\begin{equation*}
\frac {1}{2\pi^2} \mu h^{-2}\sum_n \int
 \bigl(   -(2n+1)\mu h -W  (x_2,x_3,\xi_2)\bigr)_+^{\frac{1}{2}}
\bar{\alpha} (x_2,x_3,\xi_2)\, dx_2 dx_3 d\xi_2 
\end{equation*}
and then as 
\begin{equation*}
\frac {1}{2\pi^2} \mu h^{-2}\sum_n \int
 \bigl(  -(2n+1)\mu h-  W  \bigr)_+^{\frac{1}{2}}
(1+\alpha' ) \psi F\sqrt{g}\, dx 
\end{equation*}
where now $W=W(x ,\mu^{-1})$, which in turn we rewrite as the sum of the main part $h^{-3}\cN^\MW \psi dx$ and the correction term
\begin{multline}
\frac {1}{2\pi^2} \mu h^{-2}\sum_n \int
 \bigl( -(2n+1)\mu h - W  \bigr)_+^{\frac{1}{2}}
(1+\alpha' )\psi F\sqrt{g}\, dx -\\
\frac {1}{2\pi^2} \mu h^{-2}\sum_n \int
 \bigl( -(2n+1)\mu h - V/F \bigr)_+^{\frac{1}{2}}
\psi F\sqrt{g}\, dx.
\label{18-9-1}
\end{multline}

Note that
\begin{claim}\label{18-9-2}
As $ \mu \gtrsim  (h|\log h|)^{-1}$ the  correction term does not exceed the an approximation error and thus we can skip it.
\end{claim}

Really, in this case $\varepsilon \gtrsim \mu^{-1}$. Therefore in this case we are done and our main theorems~\ref{thm-18-9-4} and \ref{thm-18-9-6} below are proven.

\section{Correction term}
\label{book_new-sect-18-9-2}

However in the case  $h^{-\frac{1}{3}}\le \mu \lesssim C(h|\log h|)^{-1}$ we need to provide an alternative expression and the estimate for the correction term (\ref{18-9-1}).

\begin{proposition}\label{pro-18-9-1}
Let $(1,2)\preceq (l,\sigma)\preceq (2,0)$. Then modulo $O(h^{-2})$ one can rewrite correction term \textup{(\ref{18-9-1})} as
\begin{multline}
h^{-3}\cN^\MW  _{3\corr } (x) =\\ 
\frac{1}{4\pi^2}  \sum_{n\ge 0}
\Bigl( \bigl(\tau -W-(2n+1)\mu h F\bigr)_+^{\frac{1}{2}} -
\bigl(\tau -V-(2n+1)\mu h F\bigr)_+^{\frac{1}{2}}\Bigr)\mu h ^{-1}F-\\
\frac{1}{6\pi^2} \bigl(W_-^{\frac{3}{2}}-V_-^{\frac{3}{2}}\bigr)h^{-3}.
\label{18-9-3}
\end{multline}
\end{proposition}

\begin{proof}
Consider 
\begin{equation*}
\frac {1}{2\pi^2} \mu h^{-2}\sum_n \int
 \bigl( -(2n+1)\mu h - W  \bigr)_+^{\frac{1}{2}}
\alpha' \psi F\sqrt{g}\, dx
\end{equation*}
and  replace the Riemannian sum by the corresponding integral; then with an error not exceeding 
\begin{equation*}
C\vartheta (\mu^{-1}) (\mu h)^{\frac{3}{2}} h^{-3}\ll \mu^{\frac{1}{2}}h^{-\frac{3}{2}}
\end{equation*}
we can replace it by 
\begin{equation}
\frac {1}{6\pi^2}\mu h^{-2} \int W_-^{\frac{3}{2}}\alpha \psi F\sqrt{g}\,dx
\label{18-9-4}
\end{equation}
where with the same error we can replace $W$ by $V$.

Consider (\ref{18-9-1}) without $\alpha'$ and replace all Riemannian sums by integrals; then we get the second term in the right-hand expression of (\ref{18-9-1}) with the opposite sign, making an error 
$O\bigl( (\mu h)^{\frac{3}{2}}h^{-3}\bigr)$. 

For $\mu \le h^{-\frac{1}{3}}$ the sum of this two expressions should be $O(h^{-2})$ and therefore (\ref{18-9-4}) must be equal to the right-hand expression of (\ref{18-9-1}). Then it is true for $\mu \ge h^{-\frac{1}{3}}$.
\end{proof}

Now let us estimate correction term.

\begin{proposition}\label{prop-18-9-2}
(i) Under non-degeneracy assumption \textup{(\ref{18-2-15})} as $(l,\sigma)\preceq (2,0)$
\begin{equation}
h^{-3}\cN^\MW_\corr \le C(\mu h)^{\frac{1}{2}}  h^{-3} \left\{\begin{aligned}
&\vartheta (\mu h) \qquad &&\text{as\ \ } \mu \le h^{-\frac{1}{2}},\\
&\vartheta (\mu ^{-1})\qquad &&\text{as\ \ } h^{-\frac{1}{2}}\le \mu \le h^{-1};
\end{aligned}\right.
\label{18-9-5}
\end{equation}
in particular it is $O(h^{-2})$ provided $(l,\sigma)\succeq (\frac{3}{2},0)$;

\medskip\noindent
(ii) Under non-degeneracy assumption \textup{(\ref{18-2-16})} as $(l,\sigma)=(2,0)$ the correction term can be ignored;

\medskip\noindent 
(iii) In the general case
\begin{equation}
h^{-3}\cN^\MW_\corr \le C(\mu h)^{\frac{1}{2}}\vartheta (\mu^{-1})h^{-3}+
C\mu h^{-2}(\vartheta (\mu^{-1}))^{\frac{1}{2}};
\label{18-9-6}
\end{equation}
in particular it is $O(h^{-2})$ provided $(l,\sigma)\succeq (2,0)$.
\end{proposition}

\begin{proof}
To prove (i),(ii) we first restore $\xi_3$ replacing $z_+^\frac{1}{2}$ by 
$\frac{1}{2}\int \uptheta(z-\xi_3^2)d\xi_3$. Then for fixed $\xi_3$ we apply arguments of the proof of proposition~\ref{prop-18-5-3} recovering estimate (\ref{18-9-5}) albeit without factor $(\mu h)^{\frac{1}{2}}$.

Then integration over zone $\{|\xi_3| \le (\mu h)^{\frac{1}{2}}\}$ results in the right-hand expression in (\ref{18-9-5}). 

On the other hand, in the zone $\{|\xi_3|\ge (\mu h)^{\frac{1}{2}} \}$ we can ``save'' one  derivative by integrating by parts over $\xi_3$ rather $x$, so there will be an extra factor $(\mu h/\xi_3^2)$ which after integration over $\xi_3$ results in $(\mu h)^{\frac{1}{2}}$ again resulting in the right-hand expression in (\ref{18-9-5}). 

\medskip\noindent
Proof of (iii) is trivial.
\end{proof}

\section{Approximation error}
\label{book_new-sect-18-9-3}

Similarly to to proposition~\ref{prop-18-9-2} one can prove

\begin{proposition}\label{prop-18-9-3}
(i) Under non-degeneracy assumption \textup{(\ref{18-2-15})} as $(l,\sigma)\preceq (2,0)$  an approximation error does not exceed
\begin{align}
&C(\mu h)^{\frac{1}{2}}  h^{-3} \vartheta(\bar{\varepsilon})\qquad
&&\text{as\ \ }\mu \le h^{-1} \label{18-9-7}\\[3pt]
&C\mu  h^{-2} \vartheta(\bar{\varepsilon})\qquad
&&\text{as\ \ }\mu \ge h^{-1} \label{18-9-8}
\end{align}
with 
\begin{equation}
\bar{\varepsilon}= (\mu^{-1}h)^{\frac{1}{2}}|\log \mu|;
\label{18-9-9}
\end{equation}
therefore this error is $O\bigl((1+\mu h)h^{-2}\bigr)$ as $(l,\sigma)\succeq (1,2)$;

\medskip\noindent
(ii) As $(l,\sigma)=(2,0)$ an approximation error is 
$O\bigl((1+\mu h)h^{-2}\bigr)$;

\medskip\noindent 
(iii) In the general case an approximation error does not exceed (modulo remainder estimate)
\begin{equation}
C(\mu h)^{\frac{1}{2}}\vartheta (\bar{\varepsilon})  h^{-3}+
C\vartheta (\bar{\varepsilon}) ^{\frac{1}{2}}\mu h^{-2}
\label{18-9-10}
\end{equation}
where the first term is $O\bigl((1+\mu h)h^{-2}\bigr)$ and the second term does not exceed $C\mu h^{-1-2/(l+2)}|\log h|^{-\sigma/(l+2)}$ (which is the part of the remainder estimate) for sure. 
\end{proposition}

\section{Main theorems}
\label{book_new-sect-18-9-4}

Getting all our beans together and getting rid off assumption $F=1$ we arrive to

\begin{theorem}\label{thm-18-9-4} 
Let $d=3$ and $A$ be a self-adjoint in $\sL ^2(X)$ operator given by \textup{(\ref{18-1-1})}. Let conditions \textup{(\ref{18-1-3})}, \textup{(\ref{18-1-6})} be fulfilled in $B(0,1)\subset X\subset \bR ^3$ with
$(1,2)\preceq (l,\sigma)\preceq (\bar{l},\bar{\sigma})\succeq (2,1)$. Let \underline{either} non-degeneracy assumption \textup{(\ref{18-2-15})} be fulfilled \underline{or} $(l,\sigma)\succeq (2,0)$ and non-degeneracy assumption \textup{(\ref{18-2-16})} be fulfilled. 

\medskip\noindent
(i) As $h^{-\frac{1}{3}}\le \mu \le Ch^{-1}$  under additional condition \textup{(\ref{18-2-27})}   for two framing approximations as in footnote~\footref{foot-18-16}  asymptotics
\begin{equation}
|\int \Bigl( \tilde{e}(x,x,0)-h^{-3}\cN _3 ^\MW (x,0) -
h^{-3} \cN _{3,\corr} ^\MW (x,0)\Bigr)\psi (x)\,dx|\le Ch ^{-2}
\label{18-9-11}
\end{equation}
holds with magnetic Weyl expression given by $\textup{(\ref{book_new-13-1-9})}_3$ and with correction term \textup{(\ref{18-9-3})}, satisfying estimate \textup{(\ref{18-9-5})};

\medskip\noindent
(ii) As $h^{-1}\le \mu$  we consider Schr\"odinger-Pauli operator with $\fz=1$  for two framing approximations as in footnote~\footref{foot-18-16}  asymptotics
\begin{equation}
\R^\MW\Def |\int \Bigl( \tilde{e}(x,x,0)-h^{-3}\cN _3 ^\MW (x,0) \Bigr)\psi (x)\,d x|\le C\mu h^{-1}
\label{18-9-12}
\end{equation}
holds with magnetic Weyl expression given by $\textup{(\ref{book_new-13-1-9})}_3$.
\end{theorem}

\begin{remark}\label{rem-18-9-5} Using rescaling technique we will be able in the next subsection to get rid of condition (\ref{18-2-27}) in this and in the next theorems.
\end{remark}

Now we are going to drop non-degeneracy assumption  as $d=3$. 

\begin{theorem} \label{thm-18-9-6} 
Let $d=3$ and $A$ be a self-adjoint in $\sL ^2(X)$ operator given by \textup{(\ref{18-1-1})}. Let conditions \textup{(\ref{18-1-3})}, \textup{(\ref{18-1-6})} be fulfilled in $B(0,1)\subset X\subset \bR ^3$ with
$(1,2)\preceq (l,\sigma)\preceq (\bar{l},\bar{\sigma})\succeq (2,1)$. Then

\medskip\noindent
(i) As $h^{-\frac{1}{3}}\le \mu \le Ch^{-1}$  under additional condition \textup{(\ref{18-2-27})}   for two framing approximations as in footnote~\footref{foot-18-16}  asymptotics
\begin{multline}
|\int \Bigl( \tilde{e}(x,x,0)-h^{-3}\cN ^\MW (x,0)-
h^{-3} \cN ^\MW _\corr (x)\Bigr) \psi (x)\, dx|\le \\[3pt]
Ch ^{-2}+C\mu h ^{-1- 2/ (l+2)}|\log h| ^{-\sigma /(l+2)}
\label{18-9-13}
\end{multline}
holds with magnetic Weyl expression given by $\textup{(\ref{book_new-13-1-9})}_3$ and with correction term \textup{(\ref{18-9-3})}, satisfying estimate \textup{(\ref{18-9-6})};

\medskip\noindent
(ii) As $h^{-1}\le \mu$, $F=1$  we consider Schr\"odinger-Pauli operator with $\fz=1$  for two framing approximations as in footnote~\footref{foot-18-16}  asymptotics
\begin{multline}
|\int \Bigl( \tilde{e}(x,x,0)- h^{-3} \cN ^\MW (x,0)\Bigr) \psi (x)dx|\le\\
C\mu h ^{-1-2/(l+2)}|\log h| ^{- \sigma/(l+2) }
\label{18-9-14}
\end{multline}
holds  with magnetic Weyl expression given by $\textup{(\ref{book_new-13-1-9})}_3$.
\end{theorem}

\begin{remark}\label{rem-18-9-7}
(i) Note that the second term in the right-hand expression (\ref{18-9-14})
prevails as $\mu \ge h^{-{l/(l+2)}}|\log h|^{\sigma/(l+2)}$.

\medskip\noindent
(ii) Using remark~\ref{18-8-10} one can prove that under non-degeneracy assumption \ref{18-8-23-m} as $\mu h\lesssim 1$ and under non-degeneracy assumption \ref{18-8-23-m-*} as $\mu h\gtrsim 1$ with $W=V/F$, $2\le m$, 
$(m,0)\prec (l,\sigma)$ remainder estimate \ref{18-8-24} holds for $\R^\MW$:
\begin{equation}
\R^\MW \le Ch^{-2}+C\mu h ^{-1-\delta};
\label{18-9-15}
\end{equation}
(iii) Furthermore, using problem~\ref{Problem-18-8-11} one can prove that in the same framework  remainder estimate  \ref{18-8-25-m} holds for $\R^\MW$:
\begin{phantomequation}\label{18-9-16}\end{phantomequation}
\begin{equation}
\R^\MW \le Ch^{-2}+C\mu h |\log h|^{m-1}.
\tag*{$\textup{(\ref*{18-9-16})}_m$}\label{18-9-16-m}
\end{equation}
\end{remark}

\begin{Problem}\label{problem-18-9-8}
Either improve estimates (\ref{18-9-5}), (\ref{18-9-6}) of proposition~\ref{prop-18-9-2} or construct counter-examples showing that the improvement is impossible.
\end{Problem}

\section{Generalizations}
\label{book_new-sect-18-9-5}

\subsection{Vanishing $V$}
\label{book_new-sect-18-9-5-1}

For $d=3$ we need condition without condition (\ref{18-2-27}) only for 
$\mu h\le 1$. 

\medskip\noindent
(i) Consider first theorem \ref{thm-18-6-10} with $\mu \le h^{-\frac{1}{3}}$. Then applying the standard scaling with $\varrho=\gamma ^{\frac{1}{2}}$ and ${\bar\gamma}_0=\mu ^{-2}$ we will find ourselves either in the case of condition (\ref{18-2-27}) fulfilled or in the case of $\mu_\new =1$ when it is not needed. Then contribution of each partition element to the remainder is
$O\bigl(h_\new ^{-2}\bigr) = O\bigl(h^{-2}\gamma ^3\bigr)$ and then the total
remainder is $O(h^{-2})$. Therefore

\begin{claim}\label{18-9-17}
Theorem \ref{thm-18-6-10} holds without condition (\ref{18-2-27}).
\end{claim}

\noindent
(ii) Consider now $h^{-\frac{1}{3}}\le \mu \le h^{-1}$. Then in the non-degenerate case (under assumption (\ref{18-2-15}) or (\ref{18-2-16})) the remainder estimate was $O(h^{-2})$; plugging $h_\new $ we get the contribution of such ball to the remainder estimate does not exceed
$Ch^{-2}\varrho ^2\gamma ^2 \asymp Ch^{-2}\gamma ^3$ which results
in the total contribution of all such balls $O(h^{-2})$. So,

\begin{claim}\label{18-9-18}
Theorem \ref{thm-18-9-4} holds without condition  (\ref{18-2-27}).
\end{claim}

\noindent
(iii) Consider the general case; recall that then estimate contains an extra term
\begin{equation}
R_3\Def C\mu h^{-1-2/(l+2)}|\log h|^{- \sigma /(l+2)}.
\label{18-9-19}
\end{equation}

However, this term does not translates well: after we plug 
$\mu_\new$, $h_\new$ it produces
\begin{equation*}
C\mu h ^{ l /(l+2)}h^{-2}
|\log (h/\varrho \gamma )|^{- \sigma /(l+2) }\times
\gamma ^{3-(l-1)/(2(l+2))}
\end{equation*}
and after multiplication by $\gamma^{-3}$ we have $\gamma$ in the negative degree.

However, we can introduce a different kind of scaling. Namely, introduce
$\gamma$ defined by 
\begin{multline}
\gamma=\min \bigl\{\eta: \ |V|\le C\vartheta (\eta),\ 
|\nabla V|\le C\eta ^{-1}\vartheta (\eta)\bigr\}, \\
\varrho = \gamma ^{\frac{l}{2}}|\log \gamma |^{-\frac {\sigma} {2}}
\label{18-9-20}
\end{multline}
where we assume so far that $(l,\sigma )\preceq (2,0)$. 

\medskip\noindent
(a) If $|V|\asymp \gamma ^l|\log \gamma |^{-\sigma }$ we are in the degenerate
situation with $V_\new$ disjoint from $0$. 

Plugging $h_\new$ into $Ch^{-2}$ we get $Ch^{-2}\varrho ^2\gamma ^2$ which does not exceed $Ch^{-2}\gamma ^3$. Plugging $h_\new$, $\mu _\new$ into $R_3$ we get
\begin{multline}
C\mu h^{-1-2/(l+2)}|\log(\frac {h} {\varrho\gamma })|^{- \sigma/(l+2) }
\times \varrho ^{2/(l+2) }\gamma ^{-l/ (l+2)}\times \gamma ^3\asymp\\
R_3 \Bigl(|\log ({\frac h{\rho \gamma}})|\cdot
|\log\gamma |/|\log h|\Bigr)^{-\sigma/(l+2)}
\times \gamma ^3
\label{18-9-21}
\end{multline}
which does not exceed $R_3\gamma^3$ as $\sigma\ge 0$ and 
$R'_3\gamma^3$ as $\sigma< 0$ with 
\begin{equation}
R'_3\Def 
C\mu h^{-1-2/(l+2)}|\log h|^{- \sigma /(l+2)}|\log (\mu h)|^{- \sigma /(l+2)}.
\label{18-9-22}
\end{equation}

\noindent
(b) If $|\nabla V|\asymp \gamma ^{l-1}|\log \gamma |^{-\sigma }$ then after
rescaling we find ourselves in the non-degenerate situation with $V_\new$ not
necessarily disjoint from $0$. However, we already examined this before
and found that contribution of the ball to the remainder estimate would not
exceed $Ch_\new ^{-2}$ which in turn does not exceed $Ch^{-2}\gamma ^3$ and then
the total contribution of such balls does not exceed $Ch^{-2}$. So,

\begin{claim}\label{18-9-23}
Let $(l,\sigma)\preceq (2,0)$. Then theorem \ref{thm-18-9-6} holds without condition (\ref{18-2-27}) as $\sigma\ge 0$ and with the last term in the right-hand expression of (\ref{18-9-13}) replace by (\ref{18-9-22}) as $\sigma <0$.
\end{claim}

\begin{Problem}\label{problem-18-9-9}
(i) Get rid off assumption (\ref{18-2-27}) in remark~\ref{rem-18-9-7}(ii), (iii). To do this use a scaling functions 
\begin{multline}
\gamma=\min \bigl\{\eta: \ \sum_{\alpha:|\alpha|\le l-1} 
|\nabla^\alpha  V/F|\le C\eta ^{-|\alpha|}\vartheta (\eta)\bigr\}, \\
\varrho = \gamma ^{\frac{l}{2}}|\log \gamma |^{-\frac {\sigma} {2}}
\label{18-9-24}
\end{multline}

\medskip\noindent
(ii) Using this get rid off assumption $(l,\sigma)\preceq (2,0)$ in claim (\ref{18-9-23}).

\medskip\noindent
(iii) Improve our arguments and prove that one can replace $R'_3$ by $R_3$ even as $\sigma<0$.
\end{Problem}

\subsection{Vanishing $F$}
\label{book_new-sect-18-9-5-2}

We assume that $g^{jk},F^k$ have an extra smoothness $1$ in comparison which is
required in theorems we refer to. This assumption allows us to avoid problems
straightening magnetic field. Let us introduce 
\begin{equation}
\gamma = \epsilon |F| +\bar{\gamma}, \qquad \bar{\gamma}=\mu^{-\frac{1}{2}}.
\label{18-9-25}
\end{equation}
Then after recalling we have $\mu_\new=\mu \gamma^2$ and $h_\new=h/\gamma$.

\medskip\noindent
(i) Consider first case $\mu \le h^{-\frac{1}{3}}$. Then 
$\mu_\new \le h_\new^{-\frac{1}{3}}$ as well and the contribution of $B\bigl(x,\gamma (x)\bigr)$ to the remainder  is
$O(h_\new^{-2})=O(h^{-2}\gamma^2)$ and the remainder does not exceed
\begin{equation}
Ch^{-2}\int \gamma(x)^{-1}\,d x 
\label{18-9-26}
\end{equation}
which does not exceed $Ch^{-2}\mu^{\frac{1}{2}}$ in the general case, and 
$Ch^{-2}|\log \mu$, $Ch^{-2}|\log \mu$ under assumptions
\begin{phantomequation}\label{18-9-27}\end{phantomequation}
\begin{gather}
|\nabla \mathbf{F}/V|\ge \epsilon_0,
\tag*{$\textup{(\ref*{18-9-27})}_1$}\label{18-9-27-1}\\[3pt]
\sum_{j,k} |\Delta_{jk}(\nabla \mathbf{F}/V)|\ge \epsilon_0
\tag*{$\textup{(\ref*{18-9-27})}_2$}\label{18-9-27-2}
\end{gather}
respectively where $\Delta_{jk}(\nabla \mathbf{F})$ denote $2\times 2$ minors of the matrix $(\nabla \mathbf{F})=(\partial_{x_j}F^k)$.

Note that assumption \ref{18-9-27-1} implies that after rescaling assumption (\ref{18-2-15}) is fulfilled.

\medskip\noindent
(ii) Consider now case $h^{-\frac{1}{3}}\le\mu \le h^{-1}$. Under condition \ref{18-9-27-1}  contribution of $B\bigl(x,\gamma (x)\bigr)$ to the remainder is
$O(h_\new ^{-2})=O(h^{-2}\gamma ^2)$ and then the total remainder is
$O(h^{-2}\log \mu )$ and it is $O(h^{-2})$ under assumption  \ref{18-9-27-2}.

We leave to the reader rather standard and tedious
\begin{Problem}\label{Problem-18-9-10}
Consider the case when assumption \ref{18-9-27-1} fails.
\end{Problem}

\medskip\noindent
(iii) Finally, for $h^{-1}\lesssim \mu \lesssim h^{-2}$ magnitude of the remainder and the
principal part strongly depend on the behavior of $\mes \{x, F(x)\le \eta\}$
as $\eta\to +0$. In particular, in the generic case 
\begin{equation}
|\det (\nabla {\mathbf F})|\ge \epsilon_0
\tag*{$\textup{(\ref*{18-9-27})}_3$}\label{18-9-27-3}
\end{equation}
the remainder is $O(\mu ^{-2}h^{-4})$ while the principal part is 
$O(\mu ^{-3}h^{-6})$ - if we consider Schr\"odinger operator (then for 
$\mu \ge C_0h^{-2}$ principal part is $0$ and remainder is $O(\mu^{-\infty})$). 

\medskip\noindent
(iii)$^*$ On the other hand,  for Schr\"odinger-Pauli operator with $\fz=1$ under assumption \ref{18-9-27-2} the principal part is $O(\mu h^{-2})$ and the remainder is $O(\mu h^{-1})$ as $\mu \ge h^{-1}$.

\bibliographystyle{alpha}

\providecommand{\bysame}{\leavevmode\hbox to3em{\hrulefill}\thinspace}

\vglue .06truein

\begin{tabular}{rrl}
&{\hskip 200 pt} &Department of Mathematics,\cr
&&University of Toronto,\cr
&&40, St.George Str.,\cr
&&Toronto, Ontario M5S 2E4\cr
&&Canada\cr
&&ivrii@math.toronto.edu\cr
&&Fax: (416)978-4107\cr
\end{tabular}

\end{document}